\documentclass[preprint,3p,times,nopreprintline]{elsarticle} 
\pdfoutput=1

\usepackage{framed,multirow}

\usepackage{mathtools}
\usepackage{amsfonts}
\usepackage{amssymb,latexsym}
\usepackage{amsbsy}
\usepackage{subfig,float,subfloat}
\usepackage{amsthm,bm}
\usepackage[normalem]{ulem}
\usepackage{commath}
\usepackage[labelfont=bf]{caption}
\usepackage{pifont}
\usepackage{listings}
\usepackage{multirow} 
\usepackage{tikz}
\usetikzlibrary{shapes.geometric}
\usetikzlibrary{arrows.meta,arrows}
\usepackage{tabulary}
\usepackage{enumitem}
\usepackage[pagewise,displaymath, mathlines]{lineno}
\usepackage[linkcolor=blue,colorlinks=true]{hyperref}
\usepackage{cleveref}
\crefname{mythm}{Theorem}{Theorems}
\crefformat{mythm}{Theorem~#2#1#3}
\crefname{remark}{Remark}{Remarks}
\crefname{lemma}{Lemma}{Lemmas}
\crefname{table}{Table}{Tables}
\crefrangeformat{table}{Tables~#3#1#4 to~#5#2#6}
\crefname{figure}{Figure}{Figures}
\crefrangeformat{figure}{Figures~#3#1#4 to~#5#2#6}
\crefformat{equation}{#2#1#3}
\crefformat{appendix}{Appendix~#2#1#3}

\usepackage{natbib}
\usepackage{comment}
\newcommand{\red}[1]{\textcolor{black}{#1}}

\newcommand{\Cv}{\mathrm{C_{v}}}
\newcommand{\Cp}{\mathrm{C_{p}}}

\newcommand{\vecu}{\normalfont{\textbf{u}}}
\newcommand{\vecv}{\normalfont{\textbf{v}}}
\newcommand{\Id}{\normalfont{\textbf{I}}}

\newcommand{\nor}{\normalfont{\textbf{n}}}

\newcommand{\q}{\normalfont{\textbf{q}}}

\newcommand{\D}{\normalfont{\textbf{D}}}
\newcommand{\U}{\normalfont{\textbf{U}}}
\newcommand{\V}{\bm{\vartheta}}
\newcommand{\W}{\bm{\mathcal{W}}}
\newcommand{\eG}{\textbf{e}_\Gamma}
\newcommand{\eP}{\textbf{e}_\Pi}

\newcommand{\dd}{\mathrm{d}}

\newcommand{\Rus}{\scalebox{.8}{$\scriptscriptstyle\rm Rus$}}
\newcommand{\HLLC}{\scalebox{.8}{$\scriptscriptstyle\rm HLLC$}}

\newcommand{\x}{\normalfont{\textbf{x}}}
\newcommand{\nik}{\normalfont{\textbf{n}}_{(i,k)j}}
\newcommand{\nkj}{\normalfont{\textbf{n}}_{i(k,j)}}
\newcommand{\GM}{\scalebox{.8}{$\scriptscriptstyle\rm GM$}}

\newcommand{\lb}{\llbracket}
\newcommand{\rb}{\rrbracket}
\newcommand{\bpsi}{\bm{\psi}}

\newcommand{\press}{\mathrm{p}}
\newtheorem{remark}{\textbf{Remark}}[section]
\newtheorem{prop}{\textbf{Proposition}}[section]
\crefname{prop}{proposition}{propositions}
\Crefname{prop}{Proposition}{Propositions}
\newtheorem{corollary}{\textbf{Corollary}}[section]
\newtheorem{mythm}{\textbf{Theorem}}[section]
\newtheorem*{leib*}{\textit{Leibniz's rule}}
\newtheorem{lemma}{\textbf{Lemma}}[section]

\newcommand{\tD}{\Tilde{\D}}
\renewcommand{\norm}[1]{\left\lVert#1\right\rVert}

\newcommand{\phib}{\boldsymbol{\phi}}

\makeatletter
\newsavebox\myboxA
\newsavebox\myboxB
\newlength\mylenA
\newcommand*\xoverline[2][1.0]{%
    \sbox{\myboxA}{$\m@th#2$}%
    \setbox\myboxB\null
    \ht\myboxB=\ht\myboxA%
    \dp\myboxB=\dp\myboxA%
    \wd\myboxB=#1\wd\myboxA
    \sbox\myboxB{$\m@th\overline{\copy\myboxB}$}
    \setlength\mylenA{\the\wd\myboxA}
    \addtolength\mylenA{-\the\wd\myboxB}%
    \ifdim\wd\myboxB<\wd\myboxA%
       \rlap{\hskip 0.5\mylenA\usebox\myboxB}{\footnotesize{\usebox\myboxA}}%
    \else
        \hskip -0.5\mylenA\rlap{\usebox\myboxA}{\footnotesize{\hskip 0.5\mylenA\usebox\myboxB}}%
    \fi}
\makeatother

\newcommand{\equaltext}[1]{\ensuremath{\stackrel{\text{#1}}{=}}}

\usepackage{scalerel,stackengine}
\stackMath
\newcommand\widhat[1]{%
\savestack{\tmpbox}{\stretchto{%
  \scaleto{%
    \scalerel*[\widthof{\ensuremath{#1}}]{\kern-.6pt\bigwedge\kern-.6pt}%
    {\rule[-\textheight/2]{1ex}{\textheight}}
  }{\textheight}%
}{0.5ex}}%
\stackon[1pt]{#1}{\tmpbox}%
}
\usepackage{url}
\usepackage{xcolor}
\definecolor{newcolor}{rgb}{.8,.349,.1}

\title{A discontinuous Galerkin spectral element method for a nonconservative compressible multicomponent flow model}

\author[1]{R\'emi Abgrall}
\cortext[cor1]{Corresponding author's email addresses: pratik.rai@polytechnique.edu (Pratik Rai), florent.renac@onera.fr (Florent Renac)}
\author[2,3]{Pratik Rai,\corref{cor1}} \corref{cor1}
\ead{pratik.rai@onera.fr}
\author[3]{Florent Renac\corref{cor1}}

\address[1]{Institute of Mathematics, University of Zürich, Switzerland}
\address[2]{CMAP, \'Ecole Polytechnique, Route de Saclay, 91128 Palaiseau Cedex, France}
\address[3]{DAAA, ONERA, Universit\'e Paris Saclay F-92322 Ch\^atillon, France}

\providecommand{\keywords}[1]
{
  \small	
  \textbf{Keywords. } #1
}

\providecommand{\AMS}[1]
{
  \small	
  \textbf{AMS subject classifications. } #1
}
\begin{document}

\maketitle
\vspace{-2cm}

\begin{abstract}
In this work, we propose an accurate, robust (the solution remains in the set of states), and stable discretization of a nonconservative model for the simulation of compressible multicomponent flows with shocks and material interfaces. We consider the gamma-based model by Shyue [J. Comput. Phys., 142 (1998), 208--242] where each component follows a stiffened gas equation of state (EOS). We here extend the framework proposed in Renac [J. Comput. Phys., 382 (2019), 1–26] and Coquel et al. [J. Comput. Phys. 431 (2021), 110135] for the discretization of hyperbolic systems, with both fluxes and nonconservative products, to unstructured meshes with curved elements in multiple space dimensions. The framework relies on a high-order discontinuous Galerkin spectral element method (DGSEM) using collocation of quadrature and interpolation points as proposed by Gassner [SIAM J. Sci. Comput., 35 (2013)] in the case of hyperbolic conservation laws. We modify the integrals over discretization elements where we replace the physical fluxes and nonconservative products by two-point numerical fluctuations. The contributions of this work are threefold. First, we analyze the semi-discrete DGSEM discretization of general hyperbolic systems with conservative and nonconservative terms and derive the conditions to obtain a scheme that is high-order accurate, free-stream preserving, and entropy stable when excluding material interfaces. Second, we design a three-point scheme with a HLLC solver for the gamma-based model that does not require a root-finding algorithm for the approximation of the nonconservative products. The scheme is proved to be robust and entropy stable for convex entropies, to preserve uniform profiles of pressure and velocity across material interfaces (material interface preservation), and to satisfy a discrete minimum principle on the specific entropy and maximum principles on the parameters of the EOS. Third, the HLLC solver is applied at interfaces in the \red{DGSEM}, while we consider two kinds of fluctuations in the integrals over discretization elements: the former is entropy conservative (EC), while the latter preserves material interfaces (CP). Time integration is performed using high-order strong-stability preserving Runge-Kutta schemes. The fully discrete scheme is shown to preserve material interfaces with CP fluctuations. Under a given condition on the time step, both EC and CP fluctuations ensure that the cell-averaged solution remains in the set of states; satisfy a minimum principle on any convex entropy and maximum principles on the EOS parameters. These results \red{have allowed us to} use existing limiters in order to restore positivity, and discrete maximum principles of degrees-of-freedom within elements. Numerical experiments in one and two space dimensions on flows with discontinuous solutions support the conclusions of our analysis and highlight stability, robustness and accuracy of the \red{proposed DGSEM} with either CP, or EC fluctuations, while the scheme with CP fluctuations is shown to offer better resolution capabilities.
\end{abstract}



\\

\AMS{65M12, 65M70, 76T10}

\allowdisplaybreaks


\section{Introduction}\label{sec: chap 5 Introduction}

The discussion in this paper focuses on the approximation in multiple space dimensions of a compressible multicomponent flow model in nonconservative form. We consider a gamma-based model \cite{shyue1998efficient} for a mixture with a stiffened gas equation of state (EOS) approximating components including both gas and compressible liquids (hereafter referred to as the SG-gamma model). The model is written in quasi-conservative form \cite{abgrall1996prevent} to preserve velocity and pressure profiles across material interfaces separating different components. The model approximates mixture quantities and presents the main advantage of being independent of the number of components. We are here interested in high-order, robust (i.e., preserving the solution in the set of admissible states), and entropy stable simulations of flows with shocks, material interfaces, and complex interactions triggering small scale flow phenomena. Numerical approximation of multicomponent and multiphase flows based on interface capturing methods has been the subject of numerous works (see, e.g. \cite{kawai2011high,saurel_pantano_ARFM_18} and references therein).

We discretize the SG-gamma model using the discontinuous Galerkin spectral element method (DGSEM) with Gauss-Lobatto quadrature rules \cite{kopriva2010quadrature}. \red{Such schemes were first introduced in previous works \cite{black1999conservative,Despres2000}, and in \cite{gassner2013SBP-SAT} the DGSEM was reinterpreted as a diagonal norm summation-by-parts (SBP) operator that falls into the general framework of conservative elementwise flux differencing schemes \cite{fisher2013high}}. In this framework, using entropy conservative (EC) numerical fluxes from Tadmor \cite{tadmor1987numerical}, semi-discrete EC finite-difference and spectral collocation schemes have been derived in \cite{fisher2013high,gassner2013SBP-SAT}. The framework has been extended to nonconservative hyperbolic systems on Cartesian meshes in \cite{renac2019entropy,rai2021entropy} by using EC numerical fluctuations from \cite{castro2013entropy}. In both frameworks, a semi-discrete entropy inequality may be obtained by replacing physical fluxes and nonconservative products with EC numerical fluxes or fluctuations within discretization elements, while using entropy stable ones at interfaces between elements. The design of the latter relies either on adding upwind-type dissipation \cite{ismail2009affordable} to EC numerical fluxes and fluctuations \cite{carpenter2014entropy,fjordholm2012arbitrarily,derigs2017novel,gassner2016split,renac2019entropy,rai2021entropy}, or on designing approximate Riemann solvers \cite{despres98,renac2021multicomp,renac17a,renac17b}. Note that the SBP operators take into account the numerical quadrature when approximating integrals compared to other schemes that require their exact evaluation to achieve entropy stability \cite{jiang1994cell,hiltebrand2014entropy,hiltebrand2018entropy}. 

Here, we extend the framework proposed in \cite{renac2019entropy} to multidimensional unstructured meshes with curved elements by using tensor multiplication of quadrature rules and function basis \cite{kopriva2010quadrature} that satisfy geometric conservation laws (the so-called metric identities \cite{kopriva_metric_id_06}) at the discrete level. This framework has been recently applied to the approximation with a well-balanced DGSEM of balance laws with geometric source terms on multidimensional high-order meshes in \cite{waruszewski2021entropy}. We here rather focus on specific properties of discretizations of nonconservative multicomponent flows: preservation of material interfaces, discrete conservation of physical fluxes, maximum principles on purely transported quantities such as the EOS parameters, and fully discrete entropy stability. The latter property presents some difficulties due to the form of the entropy associated to SG-gamma model. First, as a model for the mixture, the  properties of the individual components, such as mass and void fractions, are not known in general which prevents the evaluation of the entropy variables necessary for the derivation of EC fluctuations from the Castro et al. condition \cite{castro2013entropy}. Then, the entropy is not convex as is often the case in phase transition models \cite{helluy_seguin_phase_trans_06}\red{,} which will restrict the entropy stability across shocks not interacting with material fronts. We here circumvent these difficulties by considering a specific entropy that allows the evaluation of the entropy variables. This entropy satisfies a Gibbs relation and thus defines a complete EOS \cite{MenikoffPlohr_RP_89} and is concave with respect to two thermodynamic intensive properties of the mixture, so entropy stability can be ensured when excluding material interfaces. Let us stress that material interfaces do not require the scheme to be entropy stable, but rather a consistent approximation of the energy equation to ensure pressure equilibrium \cite{abgrall1996prevent}. We thus also design material interface preserving (CP) fluctuations to preserve at the discrete level pressure and velocity fields across such interfaces. 

We then design a HLLC solver \cite{toro1989riemann,toro1994restoration} for the SG-gamma model that does not require a root-finding algorithm to evaluate the nonconservative product in contrast to other schemes \cite{dumbser_balsara_HLLEM_16,castro_etal_HLLC_NC_13,tokareva2010hllc}. We analyze the properties of a three-point scheme using the HLLC solver and prove that the scheme is robust and entropy stable for convex entropies defining a complete EOS, preserves uniform profiles of pressure and velocity across material interfaces, and satisfies a discrete minimum principle on the specific entropy and maximum principles on the parameters of the EOS. We then apply the HLLC solver at mesh interfaces in the \red{numerical scheme} and analyze the properties of the fully discrete scheme with an explicit first-order Euler time integration. We derive conditions on the time step so that the cell-averaged solution is a convex combination of degrees-of-freedom (DOFs) and updates of three-point schemes, so the scheme inherits the properties of the three-point scheme. In particular, the \red{proposed DGSEM} satisfies a minimum principle on the entropy irrespective of the fluctuations (EC or CP) that are used in the discretization elements. As a consequence, the DGSEM with CP fluctuations within elements and the HLLC solver at interfaces is able to handle shocks and to preserve material interfaces. Time integration is performed using high-order strong-stability preserving Runge-Kutta schemes that are convex combinations of explicit Euler schemes, while linear scaling limiters \cite{zhang2010positivity,zhang2010maximum} are applied at the end of each stage to impose positivity and maximum principles at all DOFs within discretization elements.

This paper is organized as follows. The SG-gamma model and the entropy pair are described in \cref{sec: chap 5 Choice of the model}. We then introduce the semi-discrete \red{DGSEM} on multidimensional and high-order unstructured meshes in \cref{sec: The discontinuous Galerkin spectral element method (DGSEM)}. In \cref{sec: chap 5 Numerical fluxes for the gamma model}, we derive CP and EC numerical fluxes for the SG-gamma model, and in \cref{sec: chap 5 HLLC solver for interface fluxes}, we propose the HLLC approximate Riemann solver and analyze its properties. We recall the main properties of the fully-discrete \red{DGSEM} in \cref{sec: chap 5 Properties of the DGSEM scheme}. A posteriori limiters are described in \cref{sec: chap 5 a posterioiri limiters}. The results are assessed by numerical experiments in one and two space dimensions in \cref{sec: chap 5 numerical tests} and concluding remarks about this work are given in \cref{sec: chap 5 summary}.

%
%
\section{The SG-gamma model}\label{sec: chap 5 Choice of the model}

\subsection{Governing equations and thermodynamic model}\label{sec:gamma_model_PDEs}

In this work, we consider the gamma-based compressible multicomponent flow model where each component is assumed to be a stiffened gas \cite{shyue1998efficient} and refer to it as the SG-gamma model. The main interest in this model is that the number of unknowns is independent of the number of components. We are here interested in high-order approximations of the associated Cauchy problem in $d$ space dimensions for flows with $n_s$ components:

\begin{subequations} \label{Eqn: chap 5 Cauchy prob}
\begin{alignat}{2}
    \partial_t\vecu + \nabla\cdot\textbf{f}(\vecu) + \textbf{c}(\vecu)\nabla\vecu &= 0, \quad && {\bf x}\in\mathbb{R}^d,t>0,\label{eqn: chap 5 GM} \\
    \vecu(\x,0) &= \vecu_0(\x), \label{eqn: chap 5 IC}\quad && {\bf x} \in \mathbb{R}^d,
\end{alignat}
\end{subequations}

\noindent where $\x=(x_1,\dots,x_d)$ are the spatial coordinates, $\textbf{f}(\vecu)=\big(\textbf{f}_1(\vecu),\dots,\textbf{f}_d(\vecu)\big)$, $\textbf{c}(\vecu)\nabla\vecu=\sum_{i=1}^d\textbf{c}_i(\vecu)\partial_{x_i}\vecu$, and

\begin{equation} \label{Eqn: chap 5 shorthand GM vectors}
    \vecu = 
    \begin{pmatrix}
    \rho\\
    \rho \vecv\\
    \rho E\\
    \Gamma\\
    \Pi
    \end{pmatrix}, \quad
    \textbf{f}(\vecu) = 
    \begin{pmatrix}
    \rho \vecv^\top\\
    \rho \vecv\vecv^\top + \press\Id\\
    (\rho E +\press)\vecv^\top\\
    0\\
    0
    \end{pmatrix},\quad
    \textbf{c}(\vecu)\nabla\vecu = 
    \begin{pmatrix}
    0\\
    0\\
    0\\
    \vecv\cdot\nabla\Gamma\\
    \vecv\cdot\nabla\Pi
    \end{pmatrix},
\end{equation}

\noindent represent the vector of state variables, the physical fluxes and the nonconservative products\footnote{The components of the third-order tensor in (\cref{eqn: chap 5 GM}) read ${\bf c}({\bf u})_{ijk}={\bf c}_k({\bf u})_{ij}$ with the ${\bf c}_k({\bf u})$ in $\mathbb{R}^{n_{eq}\times n_{eq}}$, so the $i$th component of ${\bf c}({\bf u})\nabla{\bf u}$ reads $\sum_{j=1}^{n_{eq}}\sum_{k=1}^d{\bf c}({\bf u})_{ijk}\partial_{x_k}{\bf u}_j$. Likewise, for $\nor=(n_1,\dots,n_d)^\top$, we have $\textbf{c}(\vecu)\nor=\sum_{i=1}^dn_i\textbf{c}_i(\vecu)$.}, respectively. The mixture density, momentum, total energy, and internal energy are defined as

\begin{equation*}
 \rho = \sum_{i=1}^{n_s}\alpha_i\rho_i,\quad \rho \vecv = \sum_{i=1}^{n_s}\alpha_i\rho_i\vecv_i,\quad \rho E = \rho e + \frac{1}{2}\rho\vecv\cdot\vecv, \quad \rho e = \sum_{i=1}^{n_s}\alpha_i\rho_ie_i,
\end{equation*}
where $\rho_i$, $\vecv_i$, and $e_i$ represent the density, the velocity vector and the specific internal energy of the $i$th component. The model assumes immiscible phases and a saturation condition on the void fractions $\alpha_i$:
\begin{equation}\label{eqn: chap 5 saturation condition}
    \sum_{i=1}^{n_s}\alpha_i =1.
\end{equation}

The partial pressures are related to partial densities and specific internal energies through the stiffened gas EOS:

\begin{equation}\label{eqn: chap 5 stiffened gas EOS}
    \press_i(\rho_i,e_i) = (\gamma_i -1)\rho_ie_i - \gamma_i\press_{\infty_i}, \quad e_i=\Cv_iT_i+\frac{\press_{\infty_i}}{\rho_i}, \quad i=1, \dots, n_s,
\end{equation}

\noindent where $\gamma_i=\Cp_i/\Cv_i >1$ is the ratio of specific heats, $T_i$ is the temperature of the species, and $\press_{\infty_i}\geqslant 0$ is a pressure-like constant. Observe that when $\press_{\infty_i}= 0$ in (\cref{eqn: chap 5 stiffened gas EOS}) we recover the polytropic EOS. Assuming thermal equilibrium of the species, $T_i(\rho_i,e_i) = T(\rho,e)$, $1\leqslant i\leqslant n_s$, the EOS for the mixture is conveniently defined by \cite{shyue1998efficient} 

\begin{equation}\label{eqn: chap 5 mixture EOS}
    \frac{\press + \gamma \press_{\infty}}{(\gamma-1)} = \press\Gamma + \Pi = \rho e, \quad \rho e = \rho\Cv T + \sum_{i=1}^{n_s}\alpha_i\press_{\infty_i},
\end{equation}

\noindent where $\Cv  = \sum_{i=1}^{n_s} Y_i\Cv_i$ denotes the specific heat at constant volume of the mixture, the $Y_i=\tfrac{\alpha_i\rho_i}{\rho}$ are the mass fractions of the species, and the EOS parameters $\Gamma$ and $\Pi$ are defined by

\begin{equation}\label{eqn: chap 5 Gamma and Pi}
    \Gamma = \frac{1}{\gamma-1}=\sum_{i=1}^{n_s}\frac{\alpha_i}{\gamma_i-1},\quad
    \Pi = \frac{\gamma\press_{\infty}}{\gamma-1} = \sum_{i=1}^{n_s}\frac{\alpha_i\gamma_i\press_{\infty_i}}{\gamma_i-1}.
\end{equation}


Hyperbolicity of the SG-gamma model requires that the solutions to (\cref{Eqn: chap 5 Cauchy prob}) belong to the set of states

\begin{equation}\label{eq:set_of_states_GM}
    \Omega_{\mathrm{GM}}=\left\{\vecu\in\mathbb{R}^{n_{eq}}:\; \rho >0,\; \vecv\in\mathbb{R}^d,\; \rho e>\press_\infty,\; \Gamma>0,\; \Pi\geqslant 0\right\},
\end{equation}

\noindent with $n_{eq}=d+4$. The matrix-valued function $\sum_{i=1}^dn_i\big(\textbf{f}_i'(\vecu) + \textbf{c}_i(\vecu)\big)$ in $\mathbb{R}^{n_{eq}\times n_{eq}}$ admits real eigenvalues

\begin{equation}\label{eqn: chap 5 eigenvalues}
    \lambda_1(\vecu) = \vecv\cdot\nor-c, \quad \lambda_2(\vecu) = \dots = \lambda_{n_{eq}-1}(\vecu) = \vecv\cdot\nor, \quad
    \lambda_{n_{eq}}(\vecu) = \vecv\cdot\nor+c,
\end{equation}

\noindent for all unit vector $\nor=(n_1,\dots,n_d)$, where $c=\sqrt{\gamma(\gamma-1)(\rho e-\press_{\infty})/ \rho}$ is the speed of sound of the mixture. Here, $\{\lambda_i\}_{2\leqslant i\leqslant n_{eq}-1}$ are associated to linearly degenerate (LD) fields, while $\lambda_1$ and $\lambda_{n_{eq}}$ are associated to genuinely nonlinear (GNL) fields. Observe that (\cref{eqn: chap 5 GM}) is not strictly hyperbolic as the eigenvalues associated to the LD fields are not distinct. 

\subsection{Entropy pair}\label{sec:gamma_model_entropy}

Solutions to (\cref{Eqn: chap 5 Cauchy prob}) may develop discontinuities and (\cref{eqn: chap 5 GM}) has to be understood in a weak sense. Weak solutions are not necessarily unique and (\cref{Eqn: chap 5 Cauchy prob}) must be supplemented with further admissibility conditions to select the physical solution. We here focus on entropy inequalities

\begin{equation}\label{eqn: chap 5 entropy criteria}
    \partial_t \eta(\vecu) + \nabla\cdot \textbf{q}(\vecu) \leqslant 0, \qquad {\bf x}\in\mathbb{R}^d,t>0,
\end{equation}

\noindent for some convex entropy -- entropy flux pair $(\eta(\vecu),\textbf{q}(\vecu))$. One common way to derive such pair \red{is by defining} partial entropies of the species $s_i(\rho_i,\theta)=-\Cv_i\left(\ln\theta+(\gamma_i-1)\ln\rho_i\right)$, where $\theta=1/T$ is the inverse of the temperature, that satisfies a Gibbs relation

\begin{equation}\label{eqn: chap 5 second law of TD}
    T\dd s_i = \dd e_i - \frac{\press_i}{\rho^2_i}\dd \rho_i,
\end{equation}

\noindent so the mixture entropy reads 

\begin{equation*}
 \sum_{i=1}^{n_s} Y_is_i=-\Cv\ln\theta - \sum_{i=1}^{n_s} Y_i(\gamma_i-1)\Cv_i\ln\rho_i.
\end{equation*}

Unfortunately, we cannot use this entropy to derive EC fluxes in \cref{sec: chap 5 Numerical fluxes for the gamma model} \red{as} $Y_i$ \red{is an unknown,} that \red{cannot be evaluated} from $\vecu$ in (\cref{Eqn: chap 5 shorthand GM vectors}). Hence, here we consider an alternative pair

\begin{equation}\label{eqn: chap 5 entropy pair}
\eta(\vecu) = -\rho s, \quad \textbf{q}(\vecu) = -\rho s\vecv,
\end{equation}

\noindent where, upon introducing $\tau=\tfrac{1}{\rho}$\red{,} the covolume of the mixture, the specific entropy reads

\begin{equation}\label{eqn: chap 5 mixture physical entropy}
s(\tau,e,\Cv,\Gamma,\Pi)=\Cv\ln\left(\frac{\press+\press_\infty}{\rho^\gamma}\right) \overset{(\cref{eqn: chap 5 Gamma and Pi})}{=} \Cv\Big(\ln\Big(e-\frac{\Pi}{\Gamma+1}\tau\Big)+\frac{1}{\Gamma}\ln\tau-\ln\Gamma\Big).
\end{equation}

The rationale for considering this entropy is as follows. First, we will see in \cref{ssec: chap 5 Entropy conservative numerical fluxes} that the EC fluctuations can be explicitly computed without knowledge of the individual mass fractions.  Then, for smooth solutions of (\cref{eqn: chap 5 GM}), we have 

\begin{equation*}
 \partial_t\Cv+\vecv\cdot\nabla\Cv =0, \quad \partial_t\tau+\vecv\cdot\nabla\tau=\tau\nabla\cdot\vecv, \quad \partial_t e+\vecv\cdot\nabla e = -\tau\press\nabla\cdot\vecv,
\end{equation*}

\noindent hence

\begin{equation*}
 \partial_ts+\vecv\cdot\nabla s = \partial_\tau s\big(\partial_t\tau+\vecv\cdot\nabla \tau\big) + \partial_es\big(\partial_te+\vecv\cdot\nabla e\big)
 = \tau \nabla\cdot\vecv( \partial_\tau s - \press\partial_e s)=\Cv\nabla\cdot\vecv\left( \frac{1}{\Gamma} - \frac{\press_\infty}{\rho e-\press_\infty} - \frac{\press}{\rho e-\press_\infty} \right)\equaltext{(\cref{eqn: chap 5 mixture EOS})} 0,
\end{equation*}

\noindent so the mixture entropy is conserved, i.e., (\cref{eqn: chap 5 entropy criteria}) is an equality. Moreover, $(\tau,e)\mapsto s(\tau,e,\Cv,\Gamma,\Pi)$ is obviously strictly concave in $\Omega_{\GM}$ and from the first equality in the above relation, we conclude that the inequality in (\cref{eqn: chap 5 entropy criteria}) makes sense even if $\eta({\bf u})$ is not convex\red{,} which is also the case of the physical entropy, $\sum_{i=1}^{n_s} Y_is_i$, or in some phase transition models \cite{helluy_seguin_phase_trans_06}. Finally, differentiating (\cref{eqn: chap 5 mixture physical entropy}) while fixing $\Cv$, $\Gamma$ and $\Pi$, gives

\begin{equation}\label{eq:Gibbs_pcple_math_entropy}
\dd s = \Cv\left(\frac{\dd\press}{\press+\press_\infty}-\gamma\frac{\dd\rho}{\rho}\right) \equaltext{(\cref{eqn: chap 5 mixture EOS})} \Cv\left(\frac{(\gamma-1)\dd(\rho e)}{\press+\press_\infty}-\gamma\frac{\dd\rho}{\rho}\right) = \frac{\rho\Cv}{\rho e -\press_\infty}\left(\dd e - \frac{\press}{\rho^2}\dd\rho\right),
\end{equation}

\noindent so the entropy satisfies a Gibbs relation similar to (\cref{eqn: chap 5 second law of TD}) but with a different temperature $\tilde{T}=\tfrac{1}{\rho\Cv}(\rho e-\press_\infty)$ instead of $T$ in (\cref{eqn: chap 5 mixture EOS}). The entropy hence defines a complete EOS with both pressure and temperature \cite{MenikoffPlohr_RP_89}. This latter observation will be important in \cref{sec: chap 5 HLLC solver for interface fluxes} to prove entropy stability of the HLLC solver through the existence of local minimum entropy principles (see proof of \cref{lemma: existence of wave speeds and invariance relations}).


We end this section by deriving the entropy variables associated to the entropy $\eta(\vecu)$ in (\cref{eqn: chap 5 entropy pair}) when considering pure phases. For pure phases, we have $\dd\Gamma=\dd\Pi=\dd\Cv=0$ so we obtain
\begin{align*}
    \dd \eta = -\rho\dd s-s\dd\rho
    \equaltext{(\cref{eqn: chap 5 mixture physical entropy})}-\rho\Cv  \frac{\dd(\press+\press_{\infty})}{\press+\press_\infty} + \gamma\Cv  \dd\rho - s\dd\rho
    &\equaltext{(\cref{eqn: chap 5 mixture EOS})} -\rho\Cv  (\gamma-1)\frac{\dd \rho e}{\press+\press_\infty}  + \gamma\Cv  \dd\rho - s\dd\rho\\
    &=-\rho\frac{(\gamma-1)\Cv  }{\press+\press_\infty}\left(\dd\rho E-\vecv\cdot\dd\rho\vecv+\frac{\vecv\cdot\vecv}{2}\dd\rho\right) + \gamma\Cv  \dd\rho - s\dd\rho\\
    &\equaltext{(\cref{eq:def_zeta})} -\zeta \dd\rho E+\zeta \vecv\cdot\dd\rho\vecv + \left(\gamma\Cv  -s-\zeta\frac{\vecv\cdot\vecv}{2}\right)\dd\rho,
\end{align*}

\noindent where 
\begin{equation}\label{eq:def_zeta}
 \zeta = \frac{(\gamma-1)\Cv  \rho}{\press+\press_{\infty}},
\end{equation}

\noindent and we get the following expression of the entropy variables 
\begin{equation}\label{eqn: chap 5 entropy variables}
    \V(\vecu)= \frac{\partial}{\partial \vecu}\eta(\vecu) = 
    \begin{pmatrix}
    \gamma\Cv  -s-\zeta\frac{\vecv\cdot\vecv}{2}\\
    \zeta \vecv\\
    -\zeta\\ 
    0\\
    0 
    \end{pmatrix}.
\end{equation}

%
%
%


%
%
\section{The discontinuous Galerkin spectral element method (DGSEM)}\label{sec: The discontinuous Galerkin spectral element method (DGSEM)}

In this section, we recall the DGSEM framework \cite{rai2021entropy,kopriva2010quadrature,renac2019entropy} which is used to discretize the Cauchy problem (\cref{Eqn: chap 5 Cauchy prob}). Here, the space domain $\Omega=\mathbb{R}^d$ is discretized using a mesh $\Omega_h$ consisting of nonoverlapping and nonempty cells $\kappa$ (quadrangles for $d=2$ and hexahedra for $d=3$) forming a partition of $\Omega$. By ${\cal E}_h$ we denote the set of interfaces in $\Omega_h$. For the sake of clarity, we introduce the DGSEM in two space dimensions $d=2$, as the extension to $d=3$ is straightforward while its derivation for $d=1$ can be found in \cite{rai2021entropy,renac2019entropy}.


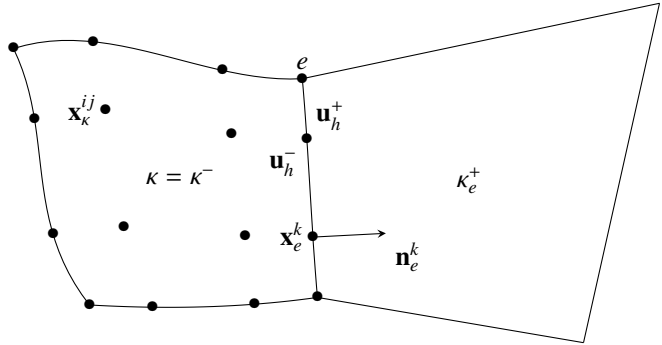
\begin{figure}[htb]
    \centering
\begin{tikzpicture}
[declare function={la(\x) = (\x+1/sqrt(5))/(-1+1/sqrt(5)) * (\x-1/sqrt(5))/(-1-1/sqrt(5)) * (\x-1)/(-2);
                   lb(\x) = (\x+1)/(-1/sqrt(5)+1) * (\x-1/sqrt(5))/(-2/sqrt(5)) * (\x-1)/(-1/sqrt(5)-1);
                   lc(\x) = (\x+1)/(1/sqrt(5)+1) * (\x+1/sqrt(5))/(2/sqrt(5)) * (\x-1)/(1/sqrt(5)-1);
                   ld(\x) = (\x+1)/(2) * (\x+1/sqrt(5))/(1+1/sqrt(5)) * (\x-1/sqrt(5))/(1-1/sqrt(5)); 
                   }]
\draw (1.2,1.725) node {$\kappa=\kappa^-$};
\draw (5.,1.65)   node {$\kappa_e^+$};
\draw (2.8,3.)    node[above] {$e$};
\draw (2.86,2.20) node[below left]  {${\bf u}_h^-$};
\draw (2.86,2.20) node[above right] {${\bf u}_h^+$};
\draw [>=stealth,->] (2.94,0.90) -- (3.9,0.95) ;
\draw (3.9,0.95) node[below right] {${\bf n}_e^k$};
\draw (2.94,0.90) node[left] {${\bf x}_e^k$};
\draw (0.21,2.59) node[left] {${\bf x}_\kappa^{ij}$};
\def\xad{-1.00};  \def\yad{3.40}; \def\xbd{0.05}; \def\ybd{3.49}; \def\xcd{1.75}; \def\ycd{3.11}; \def\xdd{2.80}; \def\ydd{3.00};
\def\xac{-0.723}; \def\yac{2.46}; \def\xbc{0.21}; \def\ybc{2.59}; \def\xcc{1.87}; \def\ycc{2.27}; \def\xdc{2.86}; \def\ydc{2.20};
\def\xab{-0.48};  \def\yab{0.94}; \def\xbb{0.45}; \def\ybb{1.03}; \def\xcb{2.05}; \def\ycb{0.91}; \def\xdb{2.94}; \def\ydb{0.90};
\def\xaa{0.00};   \def\yaa{0.00}; \def\xba{0.83}; \def\yba{-0.03}; \def\xca{2.17}; \def\yca{0.01}; \def\xda{3.00}; \def\yda{0.10};

\draw (\xad,\yad) node {$\bullet$}; \draw (\xbd,\ybd) node {$\bullet$}; \draw (\xcd,\ycd) node {$\bullet$}; \draw (\xdd,\ydd) node {$\bullet$}; 
\draw (\xac,\yac) node {$\bullet$}; \draw (\xbc,\ybc) node {$\bullet$}; \draw (\xcc,\ycc) node {$\bullet$}; \draw (\xdc,\ydc) node {$\bullet$}; 
\draw (\xab,\yab) node {$\bullet$}; \draw (\xbb,\ybb) node {$\bullet$}; \draw (\xcb,\ycb) node {$\bullet$}; \draw (\xdb,\ydb) node {$\bullet$}; 
\draw (\xaa,\yaa) node {$\bullet$}; \draw (\xba,\yba) node {$\bullet$}; \draw (\xca,\yca) node {$\bullet$}; \draw (\xda,\yda) node {$\bullet$}; 
\draw [domain=-1:1] plot ({la(\x)*\xaa+lb(\x)*\xba+lc(\x)*\xca+ld(\x)*\xda}, {la(\x)*\yaa+lb(\x)*\yba+lc(\x)*\yca+ld(\x)*\yda});
\draw [domain=-1:1] plot ({la(\x)*\xad+lb(\x)*\xbd+lc(\x)*\xcd+ld(\x)*\xdd}, {la(\x)*\yad+lb(\x)*\ybd+lc(\x)*\ycd+ld(\x)*\ydd});
\draw [domain=-1:1] plot ({la(\x)*\xda+lb(\x)*\xdb+lc(\x)*\xdc+ld(\x)*\xdd}, {la(\x)*\yda+lb(\x)*\ydb+lc(\x)*\ydc+ld(\x)*\ydd});
\draw [domain=-1:1] plot ({la(\x)*\xaa+lb(\x)*\xab+lc(\x)*\xac+ld(\x)*\xad}, {la(\x)*\yaa+lb(\x)*\yab+lc(\x)*\yac+ld(\x)*\yad});
\draw [>=stealth,-] (3.,0.1) -- (6.5,-0.5) ;
\draw [>=stealth,-] (6.5,-0.5) -- (7.5,4.) ;
\draw [>=stealth,-] (7.5,4.) -- (2.8,3.) ;
\end{tikzpicture}
\caption{Notations for the mesh in two space dimensions ($d=2$): cell $\kappa\in\Omega_h$ with quadrature points $\x_\kappa^{ij}$ (bullets $\bullet$), edge $e\in\partial\kappa$ with quadrature point $\x_e^k$, associated unit outward normal $\nor_e^k$, traces of the approximate solution $\vecu^\pm_h$ on $e$; and adjacent cell $\kappa_e^+$ sharing edge $e$.}
\label{fig: chap 3 mesh}
\end{figure}

\subsection{Numerical approximation and function space}

Let us consider the reference element $\text{I}^2=[-1,1]^2$ with coordinates $\bm{\xi} =(\xi,\eta)$ and the reference edge $\text{I}=[-1,1]$, the functions ${\bf x}_\kappa(\bm{\xi})$ and ${\bf x}_e(\xi)$ map reference to physical element and edge, respectively. The approximate solution of (\cref{Eqn: chap 5 Cauchy prob}) is sought in the function space of piecewise polynomials

\begin{equation*}
\mathcal{V}^p_h = \big\{\phi\in L^2(\Omega_h): \phi|_{\kappa}\circ{\bf x}_\kappa(\bm{\xi}) \in \mathcal{Q}^p(\text{I}^2), \forall\kappa \in \Omega_h\big\},    
\end{equation*}

\noindent where $\mathcal{Q}^p(\text{I}^2)$ denotes the space of polynomials over the reference element $\text{I}^2$ formed by the tensor product of polynomials of degree at most $p$ in each direction. The approximate solution reads

\begin{equation}\label{eqn: chap 3 approx soln}
    \vecu_h({\bf x},t) := \sum^p_{i,j=0} \phi^{ij}_\kappa({\bf x}) \U^{ij}_\kappa(t) \quad \forall {\bf x}\in \kappa, \quad \kappa\in \Omega_h, \quad t \geqslant 0,
\end{equation}

\noindent where $\{\phi^{ij}_\kappa\}_{0\leqslant i,j\leqslant p}$ constitutes a basis of $\mathcal{V}^p_h$ restricted onto $\kappa$, with dimension $(p+1)^2$, and $\{\U^{ij}_\kappa\}_{0\leqslant i,j\leqslant p}$ are the DOFs. Let $\ell_{0\leqslant k\leqslant p}$ denote the Lagrange interpolation polynomials associated to the Gauss-Lobatto nodes over $\text{I}$: $-1=\xi_0<\xi_1<\dots<\xi_p=1$. We define the basis functions as the tensor products of these polynomials:

\begin{equation}\label{eqn: chap 3 lagrange polynomial}
    \phi^{ij}_\kappa({\bf x})=\phi^{ij}({\bf x}_\kappa(\bm{\xi}))=\ell_i(\xi)\ell_j(\eta), \quad 0\leqslant i,j\leqslant p,
\end{equation}

\noindent which satisfy the following cardinality relation 

\begin{equation}\label{eq:cardinality_relation}
    \phi^{ij}_\kappa({\bf x}^{i'j'}_\kappa)=\phi^{ij}_\kappa({\bf x}_\kappa(\bm{\xi}_{i'j'}))=\ell_i(\xi_{i'})\ell_j(\eta_{j'})=\delta_{ii'}\delta_{jj'}, \quad 0\leqslant i,i',j,j'\leqslant p,
\end{equation}

\noindent where $\delta_{ii'}$ is the Kronecker delta. The DOFs are therefore point values of the solution: $\U^{ij}_\kappa(t)=\vecu_h({\bf x}^{ij}_\kappa,t)$. Likewise, the elements $\kappa$ are interpolated on the same grid of quadrature points as the numerical solution, i.e., ${\bf x}_\kappa(\bm{\xi})=\sum_{0\leqslant i,j\leqslant p}\ell_i(\xi)\ell_j(\eta){\bf x}_\kappa^{ij}$ and $\x_e(\xi)=\sum_{k=0}^p\ell_i(\xi)\x_e^k$ (see \cref{fig: chap 3 mesh}).

The integrals over elements and interfaces are approximated by using Gauss-Lobatto quadrature rules where the quadrature and interpolation points are collocated:

\begin{equation}\label{eqn: chap 3 quadrature}
    \int_{\kappa}f({\bf x})\dd V\approx\sum^p_{i,j=0} \omega_i\omega_jJ^{ij}_\kappa f({\bf x}^{ij}_\kappa), \quad \int_e f({\bf x})dS\approx \sum^p_{k=0} \omega_kJ_e^k f({\bf x}^k_e),
\end{equation}

\noindent where $\omega_i,\omega_j > 0$ are the quadrature weights and $J^{ij}_\kappa=J_\kappa(\x^{ij}_\kappa)=|{\bf x}_\kappa'(\bm{\xi}_{ij})|$, and $J_e^k=J_e(\x_e^k)=|\x_e'(\xi_k)|$. The cell-averaged solution thus reads

\begin{equation}\label{eqn: cell-averaged soln}
    \langle\vecu\rangle_\kappa (t) := \sum^p_{i,j=0}\omega_i\omega_j\frac{J^{ij}_\kappa}{|\kappa|}\U^{ij}_\kappa(t) \approx \frac{1}{|\kappa|}\int_\kappa\vecu_h(\x,t)\dd V,
\end{equation}

\noindent where $|\kappa|=\sum^p_{i,j=0}\omega_i\omega_jJ^{ij}_\kappa$ is the volume of the cell $\kappa$.

We also introduce the discrete difference matrix with entries
\begin{equation}\label{eqn: chap 3 diff matrix}
    D_{kl}=\ell'_l(\xi_k), \quad 0\leqslant k,l \leqslant p,
\end{equation}

\noindent where the property $\sum^p_{l=0} \ell_l \equiv 1$ implies

\begin{equation}\label{eqn: chap 3 sum_Dkl_vanishes}
    \sum^p_{l=0} D_{kl}= 0 \quad \forall 0 \leqslant k \leqslant p.
\end{equation}

The discrete difference matrix is known to satisfy the SBP property \cite{kopriva2010quadrature}:

\begin{equation}\label{eqn: chap 3 SBP}
    \omega_k D_{kl} + \omega_l D_{lk} =\delta_{kp}\delta_{lp}-\delta_{k0}\delta_{l0} \quad \forall 0 \leqslant k,l \leqslant p.
\end{equation}

Finally, the discretization is assumed to satisfy the following metric identities \cite{kopriva_metric_id_06}

\begin{equation}\label{eqn: metric identities}
		\sum_{k=0}^p D_{ik}J_\kappa^{kj}\nabla\xi({\bf x}_\kappa^{kj}) + D_{jk}J_\kappa^{ik}\nabla\eta({\bf x}_\kappa^{ik}) = 0 \quad \forall 0\leqslant i,j \leqslant p,
\end{equation}

\noindent and volume and edge metric terms are related by the following relations when evaluated at edges

\begin{equation}\label{eq:link_vol_surf_metric}
\begin{aligned}
    &J_\kappa^{pj}\nabla\xi(\x_\kappa^{pj})=J_e(\x_\kappa^{pj}){\bf n}_e(\x_\kappa^{pj}), \quad J_\kappa^{0j}\nabla\xi(\x_\kappa^{0j})=-J_e(\x_\kappa^{0j}){\bf n}_e(\x_\kappa^{0j}) \quad \forall 0\leqslant j\leqslant p, \\ 
		&J_\kappa^{ip}\nabla\eta(\x_\kappa^{ip})=J_e(\x_\kappa^{ip}){\bf n}_e(\x_\kappa^{ip}), \quad J_\kappa^{i0}\nabla\eta(\x_\kappa^{i0})=-J_e(\x_\kappa^{i0}){\bf n}_e(\x_\kappa^{i0}) \quad \forall 0\leqslant i\leqslant p,
\end{aligned}
\end{equation}

\noindent where ${\bf n}_e$ denotes the unit normal to $e$ in $\partial\kappa$ pointing outward from $\kappa$ (see \cref{fig: chap 3 mesh}).

\subsection{Semi-discrete form}\label{ssec: chap 3 Semi-discrete form}

We start with the following semi-discrete weak form of (\cref{eqn: chap 5 GM}) (see, e.g., \cite{franquet2012runge,rai2021entropy,renac2019entropy}): find $\vecu_h$ in $({\cal V}_h^p)^{n_{eq}}$ such that

\begin{equation*}
 \sum_{\kappa\in\Omega_h} \int_\kappa v_h\big(\partial_t\vecu_h + \nabla\cdot\textbf{f}(\vecu_h)+\textbf{c}(\vecu_h)\nabla\vecu_h\big)\dd V+\sum_{e\in\mathcal{E}_h}\int_{e} v_h^-\D^-(\vecu^-_h,\vecu^+_h,\nor_e) + v_h^+\D^+(\vecu^-_h,\vecu^+_h,\nor_e) \dd S = 0 \quad \forall v_h\in{\cal V}_h^p,
\end{equation*}

\noindent where $\vecu^\pm_h({\bf x},t)=\lim_{\varepsilon\downarrow 0}\vecu_h({\bf x}\pm\varepsilon\nor_e(\x),t)$ are the traces of $\vecu_h$ at $\x$ on a given cell interface $e\in\mathcal{E}_h$ (see \cref{fig: chap 3 mesh}) and $\D^\pm(\cdot,\cdot,\cdot)$ are the numerical fluctuations that are applied at the interfaces and are considered under the form

\begin{subequations}\label{eq:def_D+_fromD-}
\begin{align}
    \D^-(\vecu^-,\vecu^+,\nor) &= \textbf{h}(\vecu^-,\vecu^+,\nor) - \textbf{f}(\vecu^-)\cdot\nor + \textbf{d}^-(\vecu^-,\vecu^+,\nor),\\
    \D^+(\vecu^-,\vecu^+,\nor) &= \textbf{f}(\vecu^+)\cdot\nor - \textbf{h}(\vecu^-,\vecu^+,\nor) + \textbf{d}^+(\vecu^-,\vecu^+,\nor),
\end{align}
\end{subequations}

\noindent to allow proper discretizations of each term in (\cref{eqn: chap 5 GM}). They satisfy the consistency relations

\begin{equation}\label{eq:consistency_interf_fluct}
 \textbf{h}(\vecu,\vecu,\nor)=\textbf{f}(\vecu)\cdot\nor, \quad \textbf{d}^\pm(\vecu,\vecu,\nor)=0 \quad\forall \vecu \in \Omega_{\GM},
\end{equation}

\noindent and will be introduced in \cref{sec: chap 5 HLLC solver for interface fluxes}.

The integrals over mesh elements should be modified where we replace the physical fluxes and nonconservative products with numerical fluctuations of the form  \cite{castro2013entropy,castro2006high,rai2021entropy,renac2019entropy}

\begin{subequations}\label{eqn: chap 5 EC fluxes abstract form}
\begin{align}
    \D^-_{X}(\vecu^-,\vecu^+,\nor) &= \textbf{h}_{X}(\vecu^-,\vecu^+,\nor) - \textbf{f}(\vecu^-)\cdot\nor + \textbf{d}_{X}^-(\vecu^-,\vecu^+,\nor),\\
    \D^+_{X}(\vecu^-,\vecu^+,\nor) &= \textbf{f}(\vecu^+)\cdot\nor - \textbf{h}_{X}(\vecu^-,\vecu^+,\nor) + \textbf{d}_{X}^+(\vecu^-,\vecu^+,\nor),
\end{align}
\end{subequations}

\noindent where the subscript $_X$ will correspond to either entropy conservative, $_{X=ec}$, or contact preserving, $_{X=cp}$, numerical fluctuations, respectively, that will be introduced in \cref{sec: chap 5 Numerical fluxes for the gamma model}. 

We then substitute $v_h$ for the Lagrange interpolation polynomials (\cref{eqn: chap 3 lagrange polynomial}) and consider the quadrature rules (\cref{eqn: chap 3 quadrature}). Using the discrete difference matrix (\cref{eqn: chap 3 diff matrix}), the semi-discrete problem reads: find $\vecu_h$ in (\cref{eqn: chap 3 approx soln}) such that

\begin{equation}\label{eqn: chap 5 modified DG semi-discrete}
\begin{aligned}
    \omega_i\omega_jJ^{ij}_\kappa\frac{\dd}{\dd t}\U^{ij}_\kappa &+ \omega_i\omega_j\sum^p_{k=0}\left(
		 D_{ik}\tilde{\D}_X(\U^{ij}_\kappa,\U^{kj}_\kappa,\nor_{(i,k)j})
		+D_{jk}\tilde{\D}_X(\U^{ij}_\kappa,\U^{ik}_\kappa,\nor_{i(j,k)})\right)\\
		&+\sum_{e\in\partial\kappa}\sum_{k=0}^p\phi_\kappa^{ij}(\x_e^k)\omega_kJ_e^k \D^-\big(\U^{ij}_\kappa,\vecu_h^+(\x_e^k,t),\nor_e^k\big)\quad \forall t >0,\quad \kappa\in\Omega_h, \quad 0\leqslant i,j\leqslant p,
\end{aligned}
\end{equation}

\noindent where by (\cref{eq:cardinality_relation}) $\phi_\kappa^{ij}(\x_e^k)=1$ if $\x_\kappa^{ij}=\x_e^k$ and $\phi_\kappa^{ij}(\x_e^k)=0$ else. The fluctuations in the volume integrals read

\begin{subequations}\label{eqn: chap 3 tilde D}
 \begin{align}
    \tD_X(\vecu^-,\vecu^+,\nor) &:= \D^-_X(\vecu^-,\vecu^+,\nor)-\D^+_X(\vecu^+,\vecu^-,\nor), \label{eqn: chap 3 tilde Da} \\
		&\underset{(\cref{eqn: chap 3 sum_Dkl_vanishes})}{\equaltext{(\cref{eqn: chap 5 EC fluxes abstract form})}} \textbf{h}_X(\vecu^-,\vecu^+,\nor) + \textbf{h}_X(\vecu^+,\vecu^-,\nor) + {\bf d}_{X}^-(\vecu^-,\vecu^+,\nor) - {\bf d}_{X}^+(\vecu^+,\vecu^-,\nor), \label{eqn: chap 3 tilde Db}
 \end{align}
\end{subequations}

\noindent and

\begin{equation}\label{eqn: normal averages}
    \nor_{(i,k)j}=\frac{1}{2}\left(J^{ij}_\kappa\nabla\xi(\bm{\xi}_{ij})+J^{kj}_\kappa\nabla\xi(\bm{\xi}_{kj})\right), \quad \nor_{i(j,k)}=\frac{1}{2}\left(J^{ij}_\kappa\nabla\eta(\bm{\xi}_{ij})+J^{ik}_\kappa\nabla\eta(\bm{\xi}_{ik})\right)
\end{equation}

\noindent must be introduced to keep conservation of the physical fluxes \cite{wintermeyer2017entropy} and preserve uniform states. The numerical flux and fluctuations in (\cref{eqn: chap 5 EC fluxes abstract form}) satisfy the consistency conditions

\begin{equation}\label{eqn: numerical fluxes consistency condition}
    \textbf{h}_X(\vecu,\vecu,\nor)=\textbf{f}(\vecu)\cdot\nor,\quad \textbf{d}_{X}^\pm(\vecu,\vecu,\nor) =0 \quad\forall \vecu \in \Omega_{\GM},
\end{equation}

\noindent and entropy conservative (EC) fluctuations $\D_{ec}^\pm(\cdot,\cdot,\cdot)$ satisfy \cite{castro2013entropy}

\begin{equation}\label{Eqn: chap 3 entropy conserv defn nc}
    \V(\vecu^-)^\top\D^-_{ec}(\vecu^-,\vecu^+,\nor)+\V(\vecu^+)^\top\D^+_{ec}(\vecu^-,\vecu^+,\nor)=\lb\q(\vecu)\rb\cdot\nor \quad\forall \vecu^\pm \in \Omega_{\GM},
\end{equation}

\noindent where $\lb v_h\rb=v^+_h-v^-_h$ and $\V(\vecu):=\nabla_{\bf u}\eta(\vecu)$ denotes the entropy variables. Additionally, the interface fluctuations (\cref{eq:def_D+_fromD-}) in (\cref{eqn: chap 5 modified DG semi-discrete}) are assumed to be entropy stable:

\begin{equation}\label{Eqn: chap 3 entropy stable defn nc}
    \V(\vecu^-)^\top\D^-(\vecu^-,\vecu^+,\nor)+\V(\vecu^+)^\top\D^+(\vecu^-,\vecu^+,\nor)\geqslant\lb\q(\vecu)\rb\cdot\nor \quad\forall \vecu^\pm \in \Omega_{\GM}.
\end{equation}

Finally, the initial condition (\cref{eqn: chap 5 IC}) is projected onto the function space by $\U^{ij}_\kappa(0) = \vecu_0(x^{ij}_\kappa)$ for all $\kappa\in\Omega_h$ and $0\leqslant i,j\leqslant p$. The theorem below summarizes the main properties of the semi-discrete scheme (\cref{eqn: chap 5 modified DG semi-discrete}) for the discretization of general systems of the form (\cref{eqn: chap 5 GM}). Points (ii) and (iv) have been proved in \cite{waruszewski2021entropy} with slightly different volume fluctuations (see \cref{rk:vol_fluct} for comments on fluctuations in (\cref{eqn: high-order approximation proof})).

\begin{mythm}\label{thm: HO accuracy and entropy inequality}
Let $\D^\pm_{X}(\cdot,\cdot,\cdot)$ in (\cref{eqn: chap 5 EC fluxes abstract form}) be consistent fluctuations with ${\bf d}_{X}^\pm(\cdot,\cdot,\cdot)$ in (\cref{eqn: chap 3 tilde D}) satisfying
\begin{linenomath*}
\begin{subequations}\label{eqn: high-order approximation proof}
\begin{align}
\normalfont{\bf d}_{X}^{\pm}(\vecu^-,\vecu^+,\nor) &= \mathcal{C}^\pm(\vecu^-,\vecu^+,\nor)\lb\vecu\rb, \label{eqn: high-order approximation proof a}\\
\mathcal{C}(\vecu^-,\vecu^+,\nor) &:= \mathcal{C}^+(\vecu^-,\vecu^+,\nor) + \mathcal{C}^-(\vecu^-,\vecu^+,\nor), \label{eqn: high-order approximation proof b}\\
\mathcal{C}(\vecu^-,\vecu^+,\nor) + \mathcal{C}(\vecu^+,\vecu^-,\nor) &= \big(\textbf{c}(\vecu^-)+\textbf{c}(\vecu^+)\big){\bf n}, \label{eqn: high-order approximation proof c}\\
\mathcal{C}(\vecu,\vecu,\nor) &= {\bf c}(\vecu){\bf n}, \label{eqn: high-order approximation proof d}
\end{align}
\end{subequations}
\end{linenomath*}

\noindent where $\nor=(n_1,\dots,n_d)^\top$, ${\bf c}(\vecu){\bf n}=\sum_{i=1}^dn_i\textbf{c}_i(\vecu)$, and $\lb \vecu\rb = \vecu^+-\vecu^-$, and let $\D^-(\cdot,\cdot,\cdot)$ in (\cref{eqn: chap 5 modified DG semi-discrete}) be consistent (\cref{eq:consistency_interf_fluct}) fluctuations. Then, the semi-discrete scheme (\cref{eqn: chap 5 modified DG semi-discrete}) has the following properties:
\begin{enumerate}[label=(\roman*)]
 \item it is a high-order accurate approximation of smooth enough solutions to (\ref{Eqn: chap 5 Cauchy prob});
 \item it preserves uniform states (free-stream preservation);
 \item the cell-averaged solution (\cref{eqn: cell-averaged soln}) satisfies the following cell-averaged semi-discrete scheme
\begin{equation}\label{eq:cell_averaged_DGSEM}
\begin{aligned}
    |\kappa|\frac{\dd}{\dd t}\langle\vecu_h\rangle_\kappa(t) &+ \sum^p_{i,j,k=0} \omega_i\omega_j\left(D_{ik} {\bf c}(\U_\kappa^{ij})\nor_{(i,k)j}\U_\kappa^{kj} + D_{jk} {\bf c}(\U_\kappa^{ij})\nor_{i(j,k)}\U_\kappa^{ik}\right)\\ 
		&+\sum_{e\in\partial\kappa}\sum_{k=0}^p\omega_kJ_e^k \left({\bf h}\big(\vecu_h^-(\x_e^k,t),\vecu_h^+(\x_e^k,t),\nor_e^k\big) + {\bf d}^-\big(\vecu_h^-(\x_e^k,t),\vecu_h^+(\x_e^k,t),\nor_e^k\big) \right) \quad \forall t >0,\quad \kappa\in\Omega_h,
\end{aligned}
\end{equation}

\noindent which ensures a discretely conservative approximation of the physical fluxes in (\cref{eqn: chap 5 GM});

 \item if the fluctuations $\D_X(\cdot,\cdot,\cdot)$ in the volume integral are further assumed to be EC (\cref{Eqn: chap 3 entropy conserv defn nc}) and the fluctuations at interfaces $\D(\cdot,\cdot,\cdot)$ are entropy stable (\cref{Eqn: chap 3 entropy stable defn nc}) for a convex entropy pair $(\eta(\vecu),\q(\vecu))$, then the following semi-discrete entropy inequality holds:
\begin{equation}\label{eqn: semi-discrete entropy inequality}
    |\kappa|\frac{\dd}{\dd t}\langle\eta(\vecu_h)\rangle_\kappa  
		+ \sum_{e\in\partial\kappa}\sum_{k=0}^p\omega_kJ_e^k Q\big(\vecu_h^-(\x_e^k,t),\vecu_h^+(\x_e^k,t),\nor_e^k\big) \leqslant 0 \quad \forall t >0,\quad \kappa\in\Omega_h,
\end{equation}

\noindent with the consistent and conservative entropy flux 

\begin{equation}\label{eq:entropy_numer_flux}
    Q(\vecu^-,\vecu^+,\nor) = \frac{1}{2}\left(\q(\vecu^-)+\q(\vecu^+)\right)\cdot\nor+\frac{1}{2}\V(\vecu^-)^\top\D^-(\vecu^-,\vecu^+,\nor)-\frac{1}{2}\V(\vecu^+)^\top\D^+(\vecu^-,\vecu^+,\nor).
\end{equation}
\end{enumerate}
\end{mythm}

\begin{proof}
\textit{Preliminary:} By the metric identities (\cref{eqn: metric identities}) and (\cref{eqn: chap 3 sum_Dkl_vanishes}), the metric terms in (\cref{eqn: normal averages}) satisfy

\begin{equation}\label{eq:metId4nor}
 \sum^p_{k=0} D_{ik}\nor_{(i,k)j} + D_{jk}\nor_{i(j,k)} = \frac{1}{2}\sum^p_{k=0} D_{ik}\big(J_\kappa^{ij}\nabla\xi(\x_\kappa^{ij})+J_\kappa^{kj}\nabla\xi(\x_\kappa^{kj})\big) +  D_{jk}\big(J_\kappa^{ij}\nabla\eta(\x_\kappa^{ij})+J_\kappa^{kj}\nabla\eta(\x_\kappa^{ik})\big) = 0.
\end{equation}

\textit{High-order accuracy:} It is sufficient to prove that the volume integral in (\cref{eqn: chap 5 modified DG semi-discrete}) is a high-order approximation of $\nabla\cdot{\bf f}(\vecu)+{\bf c}(\vecu)\nabla\vecu$ at $(\x_\kappa^{ij},t)$ for smooth enough solutions. High-order accuracy of the conservative term $\nabla\cdot{\bf f}(\vecu)$ has been proved in \cite{chen2017entropy,ranocha_18} and the high-order accuracy of ${\bf c}(\vecu)\nabla\vecu$  in one space dimension has been proved in \cite[Th.~3.2]{renac2019entropy}. Since we are using tensor products of one-dimensional operators, the proof of accuracy follows by considering each space dimension independently.

\textit{Free-stream preservation:} Let us assume that $\U_\kappa^{ij}={\bf u}$ in the space residuals for all $\kappa\in\Omega_h$ and $0\leqslant i,j\leqslant p$, then by consistency: $\D^-({\bf u},{\bf u},{\bf n}_e^k)=0$ and $\tD_X({\bf u},{\bf u},\nor)=2{\bf f}({\bf u})\cdot{\bf n}$. Further using (\cref{eq:metId4nor}), (\cref{eqn: chap 5 modified DG semi-discrete}) becomes

\begin{equation*}
 0 = \omega_i\omega_jJ^{ij}_\kappa\frac{\dd}{\dd t}\U^{ij}_\kappa + \omega_i\omega_j\sum^p_{k=0} 2{\bf f}({\bf u})(D_{ik}\nor_{(i,k)j} + D_{jk}\nor_{i(j,k)}) = \omega_i\omega_jJ^{ij}_\kappa\frac{\dd}{\dd t}\U^{ij}_\kappa.
\end{equation*}

\textit{Cell-averaged semi-discrete scheme:} Summing up (\cref{eqn: chap 5 modified DG semi-discrete}) over $0\leqslant i,j\leqslant p$ gives

\begin{equation*}
 |\kappa|\frac{\dd}{\dd t}\langle\vecu_h\rangle_\kappa(t) + \sum^p_{j=0}\omega_jA_j + \sum^p_{i=0}\omega_iB_i + \sum_{e\in\partial\kappa}\sum_{k=0}^p\omega_kJ_e^k \D^-\big(\vecu_h^-(\x_e^k,t),\vecu_h^+(\x_e^k,t),\nor_e^k\big) = 0,
\end{equation*}

\noindent where from (\cref{eqn: chap 3 tilde D}), we have

\begin{align*}
 A_j &= \sum^p_{i,k=0} \omega_iD_{ik}\tilde{\D}_X(\U^{ij}_\kappa,\U^{kj}_\kappa,\nor_{(i,k)j})
 = \sum^p_{i,k=0} \omega_iD_{ik}\big({\bf h}_X(\U^{ij}_\kappa,\U^{kj}_\kappa,\nor_{(i,k)j}) + {\bf h}_X(\U^{kj}_\kappa,\U^{ij}_\kappa,\nor_{(i,k)j})\\
 &\hspace*{6cm}+ {\bf d}^-(\U^{ij}_\kappa,\U^{kj}_\kappa,\nor_{(i,k)j}) - {\bf d}^+(\U^{kj}_\kappa,\U^{ij}_\kappa,\nor_{(i,k)j})\big) \\
B_i &= \sum^p_{j,k=0} \omega_jD_{jk}\tilde{\D}_X(\U^{ij}_\kappa,\U^{ik}_\kappa,\nor_{i(j,k)})
 = \sum^p_{j,k=0} \omega_jD_{jk}\big({\bf h}_X(\U^{ij}_\kappa,\U^{ik}_\kappa,\nor_{i(j,k)}) + {\bf h}_X(\U^{ik}_\kappa,\U^{ij}_\kappa,\nor_{i(j,k)})\\
 &\hspace*{6cm}+ {\bf d}^-(\U^{ij}_\kappa,\U^{ik}_\kappa,\nor_{i(j,k)}) - {\bf d}^+(\U^{ik}_\kappa,\U^{ij}_\kappa,\nor_{i(j,k)})\big)
\end{align*}

Let us consider the first term. Using (\cref{eq:link_vol_surf_metric}), we have

\begin{align*}
 A_j &\underset{(\cref{eqn: numerical fluxes consistency condition})}{\equaltext{(\cref{eqn: chap 3 SBP})}} J_e(\x_\kappa^{pj}){\bf f}(\U_\kappa^{pj})\cdot{\bf n}_e(\x_\kappa^{pj}) + J_e(\x_\kappa^{0j}){\bf f}(\U_\kappa^{0j})\cdot{\bf n}_e(\x_\kappa^{0j}) \\
&+ \sum^p_{i,k=0} \omega_iD_{ik}{\bf h}_X(\U^{ij}_\kappa,\U^{kj}_\kappa,\nor_{(i,k)j}) - \omega_kD_{ki}{\bf h}_X(\U^{kj}_\kappa,\U^{ij}_\kappa,\nor_{(i,k)j}) + \omega_iD_{ik}{\bf d}^-(\U^{ij}_\kappa,\U^{kj}_\kappa,\nor_{(i,k)j}) + \omega_kD_{ki}{\bf d}^+(\U^{kj}_\kappa,\U^{ij}_\kappa,\nor_{(i,k)j}) \\
 &\underset{(\cref{eqn: high-order approximation proof}a,b)}{\overset{i\leftrightarrow k}{=}} J_e(\x_\kappa^{pj}){\bf f}(\U_\kappa^{pj})\cdot{\bf n}_e(\x_\kappa^{pj}) + J_e(\x_\kappa^{0j}){\bf f}(\U_\kappa^{0j})\cdot{\bf n}_e(\x_\kappa^{0j}) + \sum^p_{i,k=0}\omega_iD_{ik} {\cal C}(\U^{ij}_\kappa,\U^{kj}_\kappa,\nor_{(i,k)j})(\U^{kj}_\kappa-\U^{ij}_\kappa) \\
 &\underset{(\cref{eqn: high-order approximation proof d})}{\equaltext{(\cref{eqn: chap 3 SBP})}} J_e(\x_\kappa^{pj})\big({\bf f}(\U_\kappa^{pj})\cdot{\bf n}_e(\x_\kappa^{pj})-{\bf c}(\U_\kappa^{pj}){\bf n}_e(\x_\kappa^{pj})\U_\kappa^{pj}\big) + J_e(\x_\kappa^{0j})\big({\bf f}(\U_\kappa^{0j})\cdot{\bf n}_e(\x_\kappa^{0j}) -{\bf c}(\U_\kappa^{0j}){\bf n}_e(\x_\kappa^{0j})\U_\kappa^{0j}\big) \\
 &+\sum^p_{i,k=0} \omega_iD_{ik}{\cal C}(\U^{ij}_\kappa,\U^{kj}_\kappa,\nor_{(i,k)j})\U^{kj}_\kappa + \omega_kD_{ki}{\cal C}(\U^{ij}_\kappa,\U^{kj}_\kappa,\nor_{(i,k)j})\U^{ij}_\kappa \\
 &\underset{(\cref{eqn: high-order approximation proof c})}{\overset{i\leftrightarrow k}{=}} J_e(\x_\kappa^{pj})\big({\bf f}(\U_\kappa^{pj})\cdot{\bf n}_e(\x_\kappa^{pj})-{\bf c}(\U_\kappa^{pj}){\bf n}_e(\x_\kappa^{pj})\U_\kappa^{pj}\big) + J_e(\x_\kappa^{0j})\big({\bf f}(\U_\kappa^{0j})\cdot{\bf n}_e(\x_\kappa^{0j}) -{\bf c}(\U_\kappa^{0j}){\bf n}_e(\x_\kappa^{0j})\U_\kappa^{0j}\big) \\
 &+\sum^p_{i,k=0} \omega_iD_{ik}\left( {\bf c}(\U_\kappa^{ij}) + {\bf c}(\U_\kappa^{kj})\right)\nor_{(i,k)j}\U_\kappa^{kj} \\
 &\underset{(\cref{eqn: high-order approximation proof d})}{\equaltext{(\cref{eqn: chap 3 SBP})}} J_e(\x_\kappa^{pj}){\bf f}(\U_\kappa^{pj})\cdot{\bf n}_e(\x_\kappa^{pj}) + J_e(\x_\kappa^{0j}){\bf f}(\U_\kappa^{0j})\cdot{\bf n}_e(\x_\kappa^{0j}) + \sum^p_{i,k=0} \omega_iD_{ik} {\bf c}(\U_\kappa^{ij})\nor_{(i,k)j}(\U_\kappa^{kj}-\U_\kappa^{ij}),
\end{align*}

\noindent where $i\leftrightarrow k$ indicates an inversion of indices $i$ and $k$ in some of the terms. Likewise, we have
\begin{equation*}
 B_i = J_e(\x_\kappa^{ip}){\bf f}(\U_\kappa^{ip})\cdot{\bf n}_e(\x_\kappa^{ip}) + J_e(\x_\kappa^{i0}){\bf f}(\U_\kappa^{i0})\cdot{\bf n}_e(\x_\kappa^{i0}) + \sum^p_{j,k=0} \omega_jD_{jk} {\bf c}(\U_\kappa^{ij})\nor_{i(j,k)}(\U_\kappa^{ik}-\U_\kappa^{ij}).
\end{equation*}

From (\cref{eq:metId4nor}) we deduce
\begin{equation*}
 \sum^p_{j=0}\omega_jA_j + \sum^p_{i=0}\omega_iB_i = \sum_{e\in\partial\kappa}\sum_{k=0}^p\omega_kJ_e^k {\bf h}\big(\vecu_h^-(\x_e^k,t)\big)\nor_e^k + 
 \sum^p_{i,j,k=0}\omega_i\omega_j\left( D_{ik} {\bf c}(\U_\kappa^{ij})\nor_{(i,k)j}\U_\kappa^{kj} + D_{jk} {\bf c}(\U_\kappa^{ij})\nor_{i(j,k)}\U_\kappa^{ik} \right),
\end{equation*}

\noindent and we get (\cref{eq:cell_averaged_DGSEM}).

\textit{Entropy stability:} Let us now consider EC fluxes in the volume integral in the semi-discrete \red{DGSEM} (\cref{eqn: chap 5 modified DG semi-discrete}). \red{We} left multiply (\ref{eqn: chap 5 modified DG semi-discrete}) by $\V_\kappa^{ij}=\V(\vecu^{ij}_\kappa)$ and sum up over $0\leqslant i,j\leqslant p$ to get

\begin{equation*}
 |\kappa|\frac{\dd}{\dd t}\langle\eta\rangle_\kappa(t) + \sum^p_{j=0}\omega_jC_j + \sum^p_{i=0}\omega_iE_i + \sum_{e\in\partial\kappa}\sum_{k=0}^p\omega_kJ_e^k \V\big(\vecu_h^-(\x_e^k,t)\big)\cdot\D^-\big(\vecu_h^-(\x_e^k,t),\vecu_h^+(\x_e^k,t),\nor_e^k\big) = 0,
\end{equation*}

\noindent where

\begin{equation*}
 C_j = \sum^p_{i,k=0} \omega_iD_{ik}\V_\kappa^{ij}\cdot\left({\bf D}^-_{ec}(\U^{ij}_\kappa,\U^{kj}_\kappa,\nor_{(i,k)j}) - {\bf D}^+_{ec}(\U^{kj}_\kappa,\U^{ij}_\kappa,\nor_{(i,k)j})\right),
 \end{equation*}
 and
 \begin{equation*}
 E_i = \sum^p_{j,k=0} \omega_jD_{jk}\V_\kappa^{ij}\cdot\left({\bf D}^-_{ec}(\U^{ij}_\kappa,\U^{ik}_\kappa,\nor_{i(j,k)}) - {\bf D}^+_{ec}(\U^{ik}_\kappa,\U^{ij}_\kappa,\nor_{i(j,k)})\right).
\end{equation*}

Using (\cref{Eqn: chap 3 entropy conserv defn nc}), we have

\begin{align*}
 C_j &= \sum^p_{i,k=0} \omega_iD_{ik}\left(\V_\kappa^{ij}\cdot{\bf D}^-_{ec}(\U^{ij}_\kappa,\U^{kj}_\kappa,\nor_{(i,k)j}) + \V_\kappa^{kj}\cdot{\bf D}^-_{ec}(\U^{kj}_\kappa,\U^{ij}_\kappa,\nor_{(i,k)j}) - \big(\q(\U^{ij}_\kappa)-\q(\U^{kj}_\kappa)\big)\cdot\nor_{(i,k)j}\right) \\
 &\underset{(\cref{eqn: numerical fluxes consistency condition})}{\equaltext{(\cref{eqn: chap 3 SBP})}} J_e(\x_\kappa^{pj}){\bf q}(\U_\kappa^{pj})\cdot{\bf n}_e(\x_\kappa^{pj}) + J_e(\x_\kappa^{0j}){\bf q}(\U_\kappa^{0j})\cdot{\bf n}_e(\x_\kappa^{0j}) \\
&+ \sum^p_{i,k=0} \omega_iD_{ik}\V_\kappa^{ij}\cdot{\bf D}^-_{ec}(\U^{ij}_\kappa,\U^{kj}_\kappa,\nor_{(i,k)j}) - \omega_kD_{ki}\V_\kappa^{kj}\cdot{\bf D}^-_{ec}(\U^{kj}_\kappa,\U^{ij}_\kappa,\nor_{(i,k)j}) - \big(\omega_iD_{ik}\q(\U^{ij}_\kappa)+\omega_kD_{ki}\q(\U^{kj}_\kappa)\big)\cdot\nor_{(i,k)j} \\
 &\underset{(\cref{eqn: chap 3 sum_Dkl_vanishes})}{\overset{i\leftrightarrow k}{=}} J_e(\x_\kappa^{pj}){\bf q}(\U_\kappa^{pj})\cdot{\bf n}_e(\x_\kappa^{pj}) + J_e(\x_\kappa^{0j}){\bf q}(\U_\kappa^{0j})\cdot{\bf n}_e(\x_\kappa^{0j}) - \sum^p_{i,k=0} \omega_iD_{ik}J_\kappa^{kj}\q(\U^{ij}_\kappa)\nabla\xi(\x_\kappa^{kj}).
\end{align*}

Likewise

\begin{equation*}
 E_i = J_e(\x_\kappa^{ip}){\bf q}(\U_\kappa^{ip})\cdot{\bf n}_e(\x_\kappa^{ip}) + J_e(\x_\kappa^{i0}){\bf q}(\U_\kappa^{i0})\cdot{\bf n}_e(\x_\kappa^{i0}) - \sum^p_{j,k=0} \omega_jD_{jk}J_\kappa^{ik}\q(\U^{ij}_\kappa)\nabla\eta(\x_\kappa^{ik}),
\end{equation*}

\noindent and again using the metric identities (\cref{eqn: metric identities}), we get

\begin{equation*}
 |\kappa|\frac{\dd}{\dd t}\langle\eta\rangle_\kappa(t) + \sum_{e\in\partial\kappa}\sum_{k=0}^p\omega_kJ_e^k \left(\V(\x_e^k)\cdot\D^-\big(\vecu_h^-(\x_e^k,t),\vecu_h^+(\x_e^k,t),\nor_e^k\big) + \q\big(\vecu_h^-(\x_e^k,t)\big)\cdot{\bf n}_e^k \right) = 0.
\end{equation*}

Then, from (\cref{eq:entropy_numer_flux}) we have 

\begin{equation*}
 \V^-\cdot\D^-(\vecu^-,\vecu^+,\nor)+\q^-\cdot\nor=Q(\vecu^-,\vecu^+,\nor)+\frac{1}{2}(\q^--\q^+)\cdot\nor+\frac{1}{2}\V^-\cdot\D^-(\vecu^-,\vecu^+,\nor)+\frac{1}{2}\V^+\cdot\D^+(\vecu^-,\vecu^+,\nor)
\end{equation*}

\noindent and since the sum of the three last terms are non-negative by (\cref{Eqn: chap 3 entropy stable defn nc}), we obtain the desired entropy inequality (\cref{eqn: semi-discrete entropy inequality}).
\end{proof}


\begin{remark}\label{rk:vol_fluct}
Fluctuations with the following expressions fall into the category (\cref{eqn: high-order approximation proof a}):

\begin{equation*}
 \mathcal{C}^\pm({\bf u}^-,{\bf u}^+,{\bf n}) = \frac{1}{2}\big(\alpha{\bf c}({\bf u}^\pm)+(1-\alpha){\bf c}({\bf u}^\mp)\big){\bf n}, \quad 0\leqslant\alpha\leqslant 1,
\end{equation*}

\noindent which belongs to the Volpert path family of schemes \cite{volpert1967spaces} and correspond to the skew-symmetric splitting $\alpha{\bf c}\nabla{\bf u} + (1-\alpha)(\nabla\cdot({\bf c}{\bf u})-(\nabla\cdot{\bf c}){\bf u})$ of the nonconservative product. Relations (\cref{eqn: high-order approximation proof b}) and (\cref{eqn: high-order approximation proof d}) indeed correspond to the consistency condition of the Volpert path family of schemes \cite{volpert1967spaces}, while (\cref{eqn: high-order approximation proof c}) is only necessary to get (\cref{eq:cell_averaged_DGSEM}) which will be useful to prove \cref{th:properties_discr_DGSEM}. 
\end{remark}

In the next two sections, we describe the fluctuations we use in the \red{present DGSEM} (\cref{eqn: chap 5 modified DG semi-discrete}) for the discretization of the SG-gamma model (\cref{Eqn: chap 5 Cauchy prob}) and (\cref{Eqn: chap 5 shorthand GM vectors}). In \cref{sec: chap 5 Numerical fluxes for the gamma model}, we first focus on designing CP and EC fluctuations that are applied in the volume integral. We then design in \cref{sec: chap 5 HLLC solver for interface fluxes} a HLLC approximate Riemann solver for (\cref{eqn: chap 5 GM}) that is applied at interfaces. 

\section{Numerical fluctuations for the volume integrals}\label{sec: chap 5 Numerical fluxes for the gamma model}

In this section we focus on numerical fluctuations (\cref{eqn: chap 5 EC fluxes abstract form}) for the scheme (\cref{eqn: chap 5 modified DG semi-discrete}) that will be applied in the volume integrals. We will make use of the Leibniz identities, which we recall here: let $a^+,a^-,b^+,b^-,c^+,c^-$ in $\mathbb{R}$ and have finite values, then we have

\begin{linenomath*}
\begin{equation}\label{eqn: chap 5 Leibniz rule}
    \displaystyle \lb ab \rb = \xoverline{a}\lb b\rb + \xoverline{b}\lb a\rb, \quad \lb abc \rb = \xoverline{a}\big(\xoverline{b}\lb c\rb + \xoverline{c} \lb b\rb\big) + \xoverline{bc}\lb a\rb,
\end{equation}
\end{linenomath*}

\noindent where $\displaystyle\xoverline{a}= \frac{a^++a^-}{2}$ is the arithmetic mean and $\lb a\rb = a^+-a^-$ the jump.

\subsection{Contact preserving numerical fluxes}\label{ssec: chap 5 contact preserving numerical fluxes}

\begin{figure}[htb]
    \centering
\begin{tikzpicture}
[declare function={la(\x) = (\x+1/sqrt(5))/(-1+1/sqrt(5)) * (\x-1/sqrt(5))/(-1-1/sqrt(5)) * (\x-1)/(-2);
                   lb(\x) = (\x+1)/(-1/sqrt(5)+1) * (\x-1/sqrt(5))/(-2/sqrt(5)) * (\x-1)/(-1/sqrt(5)-1);
                   lc(\x) = (\x+1)/(1/sqrt(5)+1) * (\x+1/sqrt(5))/(2/sqrt(5)) * (\x-1)/(1/sqrt(5)-1);
                   ld(\x) = (\x+1)/(2) * (\x+1/sqrt(5))/(1+1/sqrt(5)) * (\x-1/sqrt(5))/(1-1/sqrt(5)); 
                   }]
%
%
\def\xad{-1.00};  \def\yad{3.40}; \def\xbd{0.05}; \def\ybd{3.49}; \def\xcd{1.75}; \def\ycd{3.11}; \def\xdd{2.80}; \def\ydd{3.00};
\def\xac{-0.723}; \def\yac{2.46}; \def\xbc{0.21}; \def\ybc{2.59}; \def\xcc{1.87}; \def\ycc{2.27}; \def\xdc{2.86}; \def\ydc{2.20};
\def\xab{-0.48};  \def\yab{0.94}; \def\xbb{0.45}; \def\ybb{1.03}; \def\xcb{2.05}; \def\ycb{0.91}; \def\xdb{2.94}; \def\ydb{0.90};
\def\xaa{0.00};   \def\yaa{0.00}; \def\xba{0.83}; \def\yba{-0.03}; \def\xca{2.17}; \def\yca{0.01}; \def\xda{3.00}; \def\yda{0.10};

\draw (\xad,\yad) node {$\bullet$}; \draw (\xbd,\ybd) node {$\bullet$}; \draw (\xcd,\ycd) node {$\bullet$}; \draw (\xdd,\ydd) node {$\bullet$}; 
\draw (\xac,\yac) node {$\bullet$}; \draw (\xbc,\ybc) node {$\bullet$}; \draw (\xcc,\ycc) node {$\bullet$}; \draw (\xdc,\ydc) node {$\bullet$}; 
\draw (\xab,\yab) node {$\bullet$}; \draw (\xbb,\ybb) node {$\bullet$}; \draw (\xcb,\ycb) node {$\bullet$}; \draw (\xdb,\ydb) node {$\bullet$}; 
\draw (\xaa,\yaa) node {$\bullet$}; \draw (\xba,\yba) node {$\bullet$}; \draw (\xca,\yca) node {$\bullet$}; \draw (\xda,\yda) node {$\bullet$}; 
\draw [domain=-1:1] plot ({la(\x)*\xaa+lb(\x)*\xba+lc(\x)*\xca+ld(\x)*\xda}, {la(\x)*\yaa+lb(\x)*\yba+lc(\x)*\yca+ld(\x)*\yda});
\draw [domain=-1:1] plot ({la(\x)*\xad+lb(\x)*\xbd+lc(\x)*\xcd+ld(\x)*\xdd}, {la(\x)*\yad+lb(\x)*\ybd+lc(\x)*\ycd+ld(\x)*\ydd});
\draw [domain=-1:1] plot ({la(\x)*\xda+lb(\x)*\xdb+lc(\x)*\xdc+ld(\x)*\xdd}, {la(\x)*\yda+lb(\x)*\ydb+lc(\x)*\ydc+ld(\x)*\ydd});
\draw [domain=-1:1] plot ({la(\x)*\xaa+lb(\x)*\xab+lc(\x)*\xac+ld(\x)*\xad}, {la(\x)*\yaa+lb(\x)*\yab+lc(\x)*\yac+ld(\x)*\yad});
\draw [>=stealth,-] (3.,0.1) -- (6.5,-0.5) ;
\draw [>=stealth,-] (6.5,-0.5) -- (7.5,4.) ;
\draw [>=stealth,-] (7.5,4.) -- (2.8,3.) ;
\def \cl {purple}
\draw (0.8,1.625) node {\textcolor{\cl}{$\Omega_L$}};
\draw (2.1,1.625) node {\textcolor{\cl}{$\Omega_R$}};
\draw (5.0,1.625) node {\textcolor{\cl}{$\Omega_R$}};
\draw [line width=1.5pt,\cl] (1.5,-0.5) .. controls (1.35,1.55) .. (1.55,3.8);
\end{tikzpicture}
\caption{\red{The above representation shows two cells, where the left cell contains a material interface (line in color) which separates the states  $\Omega_L = \{\alpha_{1L},\rho_L,\Gamma_L,\pi_L\}$ and $\Omega_R = \{\alpha_{1R},\rho_R,\Gamma_R,\pi_R\}$ with discontinuous density and EOS parameters and uniform velocity $\vecv$ and pressure $\press$. Note that the present scheme does not require the material discontinuity to be aligned with mesh interfaces to preserve uniform pressure and velocity profiles across it.}}
\label{fig: mesh_interface}
\end{figure}
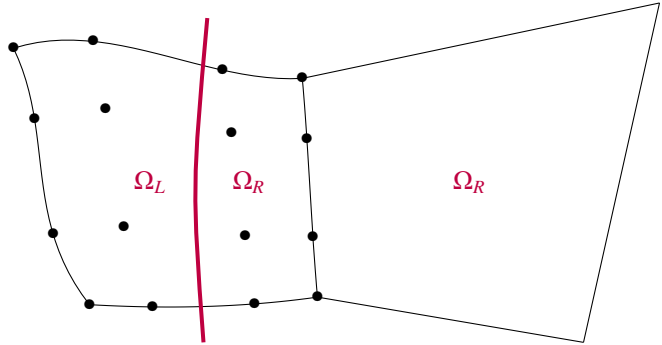

Here, we focus on deriving conditions that will ensure that the numerical fluxes maintain uniform pressure and velocity profiles across an isolated material interface. To this purpose, we introduce CP fluctuations in (\cref{eqn: chap 5 EC fluxes abstract form}), with $_{X=cp}$, where $\textbf{h}_{cp}=(h^{\rho}_{cp}, \textbf{h}^{\rho\vecv\top}_{cp}, h^{\rho E}_{cp},0,0)^\top$ is the vector of numerical fluxes for the conservative equations of mass, momentum and energy, and $\textbf{d}_{cp}^{\pm}=(0,0,0,d^{\pm}_{\Gamma},d^{\pm}_{\Pi})^\top$ is the vector for the fluctuations for the nonconservative products in (\cref{eqn: chap 5 GM}). 

For the sake of brevity, we here focus on volume fluctuations only, though the same relations may be derived for the interface fluctuations by following the same lines (see \cref{ssec: chap 5 preservation of uniform states}). Instead of using the symmetrizer in (\cref{eqn: chap 3 tilde D}), we assume that the numerical fluxes are symmetric: ${\bf h}_{cp}(\vecu^-,\vecu^+,\nor)={\bf h}_{cp}(\vecu^+,\vecu^-,\nor)$ without loss of generality (they will be, see \cref{prop: chap 5 material-interface preserving fluctuation flux}). Let us represent the conserved and nonconserved quantities in (\cref{eqn: chap 5 GM}) using $A\in\{\rho,\rho \vecv,\rho E\}$ and $B\in\{\Gamma,\Pi\}$, respectively. The \red{discretization under the DGSEM} for (\cref{eqn: chap 5 GM}) now reads

\begin{subequations}\label{eqn: chap 5 semi-disrete GM}
\begin{align}
\omega_i\omega_j J^{ij}_{\kappa}\frac{d}{dt}A^{ij}_{\kappa} &+ 2\omega_i\omega_j\sum^p_{k=0}\left(D_{ik}\textbf{h}^A_{cp}(\U^{ij}_{ \kappa},\U^{kj}_{ \kappa},\nik)+D_{jk}\textbf{h}^A_{cp}(\U^{ij}_{\kappa},\U^{ik}_{\kappa},\nkj)\right) \nonumber\\
+&\sum_{e\in\partial\kappa}\sum_{k=0}^p\phi_\kappa^{ij}(\x_e^k)\omega_kJ_e^k \D_A^-\big(\U^{ij}_{ \kappa},\vecu_h^+(\x_e^k,t),\nor_e^k\big)=0\label{eqn: chap 5 semi-discrete A}\\
\omega_i\omega_j J^{ij}_{\kappa}\frac{d}{dt}B^{ij}_{\kappa} &+ \omega_i\omega_j\sum^p_{k=0}\bigg(D_{ik}\big(\textbf{d}_{cp}^{B,-}(\U^{ij}_{\kappa},\U^{kj}_{\kappa},\textbf{n}_{(i,k)j})-\textbf{d}_{cp}^{B,+}(\U^{kj}_{\kappa},\U^{ij}_{\kappa},\textbf{n}_{(i,k)j})\big) \nonumber\\
+& D_{jk}\big(\textbf{d}_{cp}^{B,-}(\U^{ij}_{\kappa},\U^{ik}_{\kappa},\textbf{n}_{i(k,j)})-\textbf{d}_{cp}^{B,+}(\U^{ik}_{\kappa},\U^{ij}_{\kappa},\textbf{n}_{i(k,j)})\big)\bigg) +\sum_{e\in\partial\kappa}\sum_{k=0}^p\phi_\kappa^{ij}(\x_e^k)\omega_kJ_e^k \D_B^-\big(\U^{ij}_{ \kappa},\vecu_h^+(\x_e^k,t),\nor_e^k\big)=0\label{eqn: chap 5 semi-discrete B}
\end{align}
\end{subequations}


Now let us suppose that the initial condition consists of a material interface with uniform velocity, $\vecv$, and pressure, $\press$, and states $\rho_L$, $\Gamma_L$ and $\pi_L$ in $\Omega_L$ and $\rho_R$, $\Gamma_R$ and $\pi_R$ in $\Omega_R$ with $\overline{\Omega_L}\cup\overline{\Omega_R}=\overline{\Omega}$, then so do the DOFs, \red{see \cref{fig: mesh_interface}}. We now derive conditions for the numerical fluxes (\cref{eqn: chap 5 EC fluxes abstract form}) to preserve the uniform states in time. 

We, first, focus on the velocity state and impose the semi-discrete scheme (\cref{eqn: chap 5 semi-disrete GM}) to satisfy a discrete counterpart to the differential relation $\rho\dd \vecv = \dd\rho \vecv - \vecv\dd\rho = 0$. Ignoring the interface fluxes, $\omega_i\omega_jJ_\kappa^{ij}\rho_\kappa^{ij}\dd_t\vecv_\kappa^{ij}=0$ requires

\begin{equation*}
\omega_i\omega_j\sum^p_{k=0}\Big(\big(\textbf{h}^{\rho\vecv}_{cp}(\U^{ij}_\kappa,\U^{kj}_\kappa,\textbf{n}_{(i,k)j})-\vecv h^\rho_{cp}(\U^{ij}_\kappa,\U^{kj}_\kappa,\textbf{n}_{(i,k)j})\big)D_{ik}+\big(\textbf{h}^{\rho\vecv}_{cp}(\U^{ij}_\kappa,\U^{ik}_\kappa,\textbf{n}_{i(k,j)})-\vecv h^\rho_{cp}(\U^{ij}_\kappa,\U^{ik}_\kappa,\textbf{n}_{i(k,j)})\big)D_{jk}\Big) = 0,
\end{equation*}

\noindent and a sufficient condition reads

\begin{equation}\label{eqn: chap 5 sufficient condn unif. vel.}
    \textbf{h}^{\rho\vecv}_{cp}(\vecu^-,\vecu^+,\nor)=\tilde{\vecv}(\vecu^-,\vecu^+) h^{\rho}_{cp}(\vecu^-,\vecu^+,\nor)+ \tilde{\press}(\vecu^-,\vecu^+,\nor),
\end{equation}

\noindent where $\tilde{\vecv}$ and $\tilde{\press}$ are any consistent discretizations of the velocity vector $\vecv$ and pressure $\press{\bf n}$. Similarly, a semi-discrete equation for the pressure (\cref{eqn: chap 5 mixture EOS}) can be obtained by using
$\Gamma\dd\press =\dd\rho E - (\tfrac{1}{2}\vecv\cdot\vecv)\dd\rho - \press\dd\Gamma - \dd\Pi$ from (\cref{eqn: chap 5 Gamma and Pi}), and again ignoring surface contributions in (\cref{eqn: chap 5 semi-disrete GM}), $\omega_i\omega_jJ_\kappa^{ij}\Gamma_\kappa^{ij}\dd_t\press(\U_\kappa^{ij})=0$ requires

\begin{equation*}
\begin{aligned}
\omega_i\omega_j\sum^p_{k=0}&D_{ik}\bigg(2h^{\rho E}_{cp}(\textbf{U}^{ij}_\kappa,\textbf{U}^{kj}_\kappa,\textbf{n}_{(i,k)j})-\vecv\cdot\vecv h^\rho_{cp}(\textbf{U}^{ij}_\kappa,\textbf{U}^{kj}_\kappa,\nor_{(i,k)j}) -\press\left(d^-_\Gamma(\textbf{U}^{ij}_\kappa,\textbf{U}^{kj}_\kappa,\textbf{n}_{(i,k)j})-d^+_\Gamma(\textbf{U}^{kj}_\kappa,\textbf{U}^{ij}_\kappa,\textbf{n}_{(i,k)j})\right)\\
&- d^-_\Pi(\textbf{U}^{ij}_\kappa,\textbf{U}^{kj}_\kappa,\textbf{n}_{(i,k)j})+d^+_\Pi(\textbf{U}^{kj}_\kappa,\textbf{U}^{ij}_\kappa,\textbf{n}_{(i,k)j}) \bigg)\\
+&D_{jk}\bigg(2h^{\rho E}_{cp}(\textbf{U}^{ij}_\kappa,\textbf{U}^{ik}_\kappa,\textbf{n}_{i(k,j)})-\vecv\cdot\vecv h^\rho_{cp}(\textbf{U}^{ij}_\kappa,\textbf{U}^{ik}_\kappa,\nor_{i(j,k)}) -\press\left(d^-_\Gamma(\textbf{U}^{ij}_\kappa,\textbf{U}^{ik}_\kappa,\textbf{n}_{i(k,j)})-d^+_\Gamma(\textbf{U}^{ik}_\kappa,\textbf{U}^{ij}_\kappa,\textbf{n}_{i(k,j)})\right)\\
&- d^-_\Pi(\textbf{U}^{ij}_\kappa,\textbf{U}^{ik}_\kappa,\textbf{n}_{i(k,j)})+d^+_\Pi(\textbf{U}^{ik}_\kappa,\textbf{U}^{ij}_\kappa,\textbf{n}_{i(k,j)}) \bigg) = 0,
\end{aligned}
\end{equation*}

\noindent and subtracting the trivial quantity $2\omega_i\omega_j\sum^p_{k=0}\big({\bf f}^{\rho E}(\textbf{U}^{ij}_\kappa)-\tfrac{\vecv\cdot\vecv}{2}{\bf f}^{\rho}(\textbf{U}^{ij}_\kappa)\big)\cdot(D_{ik}\nor_{(i,k)j}+D_{jk}\nor_{i(j,k)})=0$, from (\cref{eq:metId4nor}), a sufficient condition reads

\begin{multline}\label{eqn: chap 5 sufficient condn unif. press.}
    h^{\rho E}_{cp}(\vecu^-,\vecu^+,\nor) - (\rho E^- + \press)\vecv\cdot\nor = \frac{\vecv\cdot\vecv}{2}\big(h^{\rho}_{cp}(\vecu^-,\vecu^+,\nor)-\rho^-\vecv\cdot\nor\big) \\ +\frac{\press}{2} \left(d^{-}_{\Gamma}(\vecu^-,\vecu^+,\nor) - d^{+}_{\Gamma}(\vecu^+,\vecu^-,\nor)\right) + \frac{1}{2}\big(d^{-}_\Pi(\vecu^-,\vecu^+,\nor) - d^{+}_\Pi(\vecu^+,\vecu^-,\nor)\big).
\end{multline}

We can now propose CP fluxes for the volume integral.

\begin{prop}\label{prop: chap 5 material-interface preserving fluctuation flux}
Numerical fluxes of the form (\cref{eqn: chap 5 EC fluxes abstract form}) where
\begin{equation}\label{eqn: chap 5 material-interface preserving flux b}
\normalfont{\textbf{h}}_{cp}(\vecu^-,\vecu^+,\nor)=
\begin{pmatrix}
\xoverline{\rho}\, \xoverline{\vecv}\cdot\nor\\
\xoverline{\rho}\,\xoverline{\vecv}\,\xoverline{\vecv}\cdot\nor+\xoverline{\press}\nor\\
(\xoverline{\rho E}+\xoverline{\press})\xoverline{\vecv}\cdot\nor\\
0\\
0
\end{pmatrix},
\quad
\textbf{d}_{cp}^{\pm}(\vecu^-,\vecu^+,\nor)=\frac{1}{2}\vecv^\pm\cdot\nor
\begin{pmatrix}
0\\
0\\
0\\
\lb\Gamma\rb\\
\lb\Pi\rb
\end{pmatrix},
\end{equation}

\noindent preserve the uniform pressure and velocity fields across contact discontinuities and material interfaces for the SG-gamma model (\cref{eqn: chap 5 GM}) and (\cref{Eqn: chap 5 shorthand GM vectors}), with the mixture EOS (\cref{eqn: chap 5 mixture EOS}).
\end{prop}

\begin{proof}
Checking that condition (\cref{eqn: chap 5 sufficient condn unif. vel.}) holds for (\cref{eqn: chap 5 material-interface preserving flux b}) is direct, wile using the mixture EOS (\cref{eqn: chap 5 mixture EOS}), we get

\begin{equation*}
    \frac{\press}{2} \left(d^{-}_{\Gamma}(\vecu^-,\vecu^+,\nor) - d^{+}_{\Gamma}(\vecu^+,\vecu^-,\nor)\right) + \frac{1}{2}\big(d^{-}_\Pi(\vecu^-,\vecu^+,\nor) - d^{+}_\Pi(\vecu^+,\vecu^-,\nor)\big) = \frac{\vecv\cdot\nor}{2}\big(\press\lb\Gamma\rb + \lb\Pi\rb\big) \overset{(\cref{eqn: chap 5 Leibniz rule})}{=} \frac{\vecv\cdot\nor}{2}\left(\lb\rho E\rb -\frac{\vecv\cdot\vecv}{2}\lb\rho\rb\right),
\end{equation*}

\noindent so (\cref{eqn: chap 5 sufficient condn unif. press.}) holds as well.
\end{proof}

\begin{remark}
The CP numerical fluxes (\cref{eqn: chap 5 material-interface preserving flux b}) are similar to the one proposed in \cite{kennedy2008reduced}. Here we have modified the contributions towards the energy equation to satisfy (\cref{eqn: chap 5 sufficient condn unif. press.}).
\end{remark}

\subsection{Entropy conservative numerical fluxes}\label{ssec: chap 5 Entropy conservative numerical fluxes}

We, now, propose EC fluxes for the SG-gamma model that are applied to the modified volume integral in (\cref{eqn: chap 5 modified DG semi-discrete}) and, according to \cref{thm: HO accuracy and entropy inequality}, these numerical fluxes will contribute to the entropy stability of the numerical scheme. 

\begin{prop}\label{prop: chap 5 EC fluctuation flux}
Consider the entropy pair (\cref{eqn: chap 5 entropy pair}) with (\cref{eqn: chap 5 mixture physical entropy}), then fluctuations of the form (\cref{eqn: chap 5 EC fluxes abstract form}) with
\begin{equation}\label{eqn: chap 5 defn of EC fluxes}
    \normalfont{\textbf{h}}_{ec}(\vecu^-,\vecu^+,\nor) =
    \begin{pmatrix}
    \hat{\rho}\xoverline{\vecv}\cdot\nor\\
    {\bf h}^{\rho \vecv}_{ec}\\
    h^{\rho E}_{ec}\\
    0\\
    0
    \end{pmatrix},\quad
    \textbf{d}_{ec}^\pm(\vecu^-,\vecu^+,\nor) = \frac{1}{2}\vecv^\pm\cdot\nor
    \begin{pmatrix}
    0\\
    0\\
    0\\
    \lb \Gamma\rb\\
    \lb \Pi\rb
    \end{pmatrix},
\end{equation}

\noindent where
\begin{equation*}
{\bf h}^{\rho \vecv}_{ec}(\vecu^-,\vecu^+,\nor)=\hat{\rho}(\xoverline{\vecv}\cdot\nor)\xoverline{\vecv}+\left(\xoverline{\rho}\xoverline{(\gamma-1)}\xoverline{\left(\displaystyle\frac{\Cv}{\zeta}\right)} - \xoverline{\press_{\infty}}\right)\nor,\quad
h^{\rho E}_{ec}(\vecu^-,\vecu^+,\nor) = \left(\widhat{\left(\displaystyle\frac{\Cv  }{\zeta}\right)}+\frac{\vecv^-\cdot \vecv^+}{2}\right)\hat{\rho}\xoverline{\vecv}\cdot\nor + \left(\xoverline{\rho}\xoverline{(\gamma-1)}\xoverline{\left(\displaystyle\frac{\Cv  }{\zeta}\right)}\right)\xoverline{\vecv}\cdot\nor,
\end{equation*}

\noindent and $\tfrac{\Cv}{\zeta}=\tfrac{\Gamma}{\rho}(\press+\tfrac{\Pi}{\Gamma+1})$ from (\cref{eq:def_zeta}), are EC in the sense (\cref{Eqn: chap 3 entropy conserv defn nc}) when excluding material interfaces, i.e., $\lb\Gamma\rb\equiv\lb\Pi\rb\equiv 0$. 
\end{prop}

\begin{proof}
As we consider pure phases, the system (\cref{eqn: chap 5 GM}) is conservative and (\cref{Eqn: chap 3 entropy conserv defn nc}) reduces to the Tadmor condition \cite{tadmor1987numerical}
\begin{equation*}
    \textbf{h}_{ec}(\vecu^-,\vecu^+,\nor)\cdot\lb\V\rb - \lb\bpsi(\vecu)\rb\cdot\nor =0 \quad \forall \vecu^\pm \in \Omega_{\GM},
\end{equation*}

\noindent where the entropy variables $\V$ are defined in (\cref{eqn: chap 5 entropy variables}), and $\bpsi \equiv \textbf{f}^\top\V-\textbf{q}$ is the entropy potential and reads
\begin{equation*}
    \bpsi(\vecu) \equaltext{(\cref{eqn: chap 5 entropy variables})} -\zeta (\rho E+\press)\vecv + \zeta \vecv(\rho \vecv\cdot\vecv+\press)+\left(\gamma\Cv  -s-\zeta\frac{\vecv\cdot\vecv}{2}\right)\rho \vecv + \rho s\vecv
    =(\gamma\Cv  -\zeta e)\rho \vecv
    =\left((\gamma-1)\Cv  \rho  -\press_\infty \zeta \right)\vecv,
\end{equation*}

\noindent so 
\begin{equation}\label{eqn: chap 5 lb psi rb}
    \lb\bpsi(\vecu)\rb\cdot\nor \equaltext{(\cref{eqn: chap 5 Leibniz rule})} \big(\xoverline{(\gamma-1)\Cv  }\lb\rho  \vecv\rb-\xoverline{\press_\infty}\lb\zeta \vecv\rb\big)\cdot\nor.
\end{equation}

Using the definition of $\zeta$ in (\cref{eq:def_zeta}), the entropy of the mixture (\cref{eqn: chap 5 mixture physical entropy}) may be reformulated as
\begin{equation}\label{eqn: chap 5 physical entropy of conservative}
    s = \Cv  \ln\left(\frac{\press+\press_\infty}{\rho^\gamma}\right) = -\Cv  \ln\zeta-(\gamma-1)\Cv  \ln\rho +\Cv  \ln\big((\gamma-1)\Cv\big),
\end{equation}

\noindent then we have
\begin{align*}
    \textbf{h}^\top_{ec}(\vecu^-,\vecu^+,\nor)\lb\V(\vecu)\rb &= \lb\gamma\Cv  -s-\zeta\frac{\vecv\cdot\vecv}{2}\rb \hat{\rho}\xoverline{\vecv}\cdot\nor +\left(\xoverline{\vecv}\cdot\xoverline{\vecv}\hat{\rho}+\frac{(\gamma-1)\Cv  \xoverline{\rho}}{\xoverline{\zeta}}-\press_\infty\right)\lb\zeta \vecv\rb\cdot\nor\\
    &\hspace{4cm}-\lb\zeta\rb\left(\left(\frac{\Cv  }{\hat{\zeta}}+\frac{\vecv^-\cdot\vecv^+}{2}\right)\hat{\rho}+\frac{(\gamma-1)\Cv  }{\xoverline{\zeta}}\xoverline{\rho}\right)\xoverline{\vecv}\cdot\nor\\
    &\underset{(\cref{eqn: chap 5 Leibniz rule})}{\equaltext{(\cref{eqn: chap 5 physical entropy of conservative})}}\left(\frac{\lb\zeta\rb}{\hat{\zeta}}+(\gamma-1)\frac{\lb\rho\rb}{\hat{\rho}}\right) \Cv  \hat{\rho}\xoverline{\vecv}\cdot\nor +\lb\zeta \vecv\rb\left(\xoverline{\vecv}\cdot\xoverline{\vecv}\hat{\rho}+\frac{(\gamma-1)\Cv  \xoverline{\rho}}{\xoverline{\zeta}}-\press_\infty\right)\cdot\nor\\
    &\hspace{4cm}-\lb\zeta\rb\left(\left(\frac{\Cv  }{\hat{\zeta}}+\frac{\vecv^-\cdot \vecv^+}{2}\right)\hat{\rho}+\frac{(\gamma-1)\Cv  }{\xoverline{\zeta}}\xoverline{\rho}\right)\xoverline{\vecv}\cdot\nor\\
    &=\big((\gamma-1)\Cv  \lb\rho \vecv\rb-\press_\infty\lb\zeta \vecv\rb \big)\cdot\nor,
\end{align*}

\noindent which cancels out with (\cref{eqn: chap 5 lb psi rb}), so the proof is complete.
\end{proof}

\section{An HLLC Riemann solver for the cell interfaces}\label{sec: chap 5 HLLC solver for interface fluxes}

We now look for two-point numerical fluctuations at interfaces and design an HLLC approximate Riemann solver for the SG-gamma model (\cref{eqn: chap 5 GM}) and analyze its properties. For the sake of generality, we assume that the entropy $\eta(\vecu)$ in (\cref{eqn: chap 5 entropy criteria}) is convex, which excludes material interfaces. 

\subsection{One-dimensional Riemann problem}

We are interested in approximating solutions to the following Riemann problem in a given unit direction ${\bf n}$ in $\mathbb{R}^d$:
\begin{subequations}\label{eq:RP_1D_gamma_model}
\begin{equation}\label{eq:RP_1D_gamma_model_a}
    \partial_t\vecu +\partial_x\textbf{f}(\vecu)\cdot\nor + \textbf{c}_{\nor}(\vecu)\partial_x\vecu = 0,
\end{equation}

\noindent where $x={\bf x}\cdot\nor$ and $\textbf{c}_{\nor}(\vecu)=\textbf{c}(\vecu)\nor=\sum^d_{i=1}n_i\textbf{c}_i(\vecu)$, together with initial data
\begin{equation}\label{eq:RP_1D_gamma_model_b}
    \vecu_0(x)=
    \left\{\begin{array}{ll}
    \vecu_L, & x <0,\\
    \vecu_R, & x >0.
    \end{array}\right.
\end{equation}
\end{subequations}

By $\W(\tfrac{x}{t},\vecu_L,\vecu_R,\nor)$ we denote the exact entropy weak solution to (\cref{eq:RP_1D_gamma_model}) for $t>0$. Following \cite{harten1983upstream,gallice_ARS_NC_03}, we integrate (\cref{eq:RP_1D_gamma_model_a}) over the control volume $[-\tfrac{h}{2},\tfrac{h}{2}]\times[0,\Delta t]$ with $h>0$ and $\Delta t>0$ the space and time steps, respectively. Using (\cref{eq:RP_1D_gamma_model_b}), we obtain
\begin{equation}\label{eqn: chap 5 exact integral GM a}
 \int^{\frac{h}{2}}_{-\frac{h}{2}} \W\left(\frac{x}{\Delta t};\vecu_L,\vecu_R,\nor\right)\dd x-\frac{h}{2}(\vecu_L+\vecu_R) + \Delta t\left(\textbf{f}(\vecu_R)-\textbf{f}(\vecu_L)\right)\cdot\nor + \int^{\Delta t}_0\int^{\frac{h}{2}}_{-\frac{h}{2}}\textbf{c}_{\nor}(\vecu)\partial_x\vecu\dd x\dd t = 0.
\end{equation}

Note that both $\Gamma$ and $\Pi$ are continuous across shocks and discontinuous across the intermediate contact wave. Let us introduce the two last components $\eG = (0,0,0,1,0)^\top$ and $\eP = (0,0,0,0,1)^\top$ of the canonical basis in $\mathbb{R}^{n_{eq}}$, and define $u:=\vecv\cdot\nor$. By (\cref{eqn: chap 5 exact integral GM a}) and (\cref{Eqn: chap 5 shorthand GM vectors}) we have
\begin{align*}
    \int^{\Delta t}_0\int^{\frac{h}{2}}_{-\frac{h}{2}}\textbf{c}_{\nor}(\vecu)\partial_x\vecu\dd x\dd t
    &= \frac{h}{2}(\vecu_L+\vecu_R) - \Delta t\left(\textbf{f}(\vecu_R)-\textbf{f}(\vecu_L)\right)\cdot\nor - \int^{\frac{h}{2}}_{-\frac{h}{2}} \W\left(\frac{x}{\Delta t};\vecu_L,\vecu_R,\nor\right)\dd x \nonumber\\
		&= \frac{h}{2}(\Gamma_L+\Gamma_R)\eG + \frac{h}{2}(\Pi_L+\Pi_R)\eP - \left(\frac{h}{2}-u^\star\Delta t\right)(\Gamma_R\eG+\Pi_R\eP) - \left(u^\star\Delta t+\frac{h}{2}\right)(\Gamma_L\eG+\Pi_L\eP) \nonumber\\
    &= \Delta t u^\star\big((\Gamma_R-\Gamma_L)\eG+(\Pi_R-\Pi_L)\eP\big), 
\end{align*}

\noindent where $u_L$, $u^\star$, and $u_R$ are the normal velocity components in the left state, the star region, and the right states, respectively. The integral form for (\cref{eqn: chap 5 GM}) thus reads
\begin{equation}\label{eqn: chap 5 exact integral GM b}
    \frac{1}{\Delta t}\int^{\frac{h}{2}}_{-\frac{h}{2}} \W\left(\frac{x}{\Delta t};\vecu_L,\vecu_R,\nor\right)\dd x -\frac{h}{2\Delta t}(\vecu_L+\vecu_R)+\big(\textbf{f}(\vecu_R)-\textbf{f}(\vecu_L)\big)\cdot\nor + u^\star\big((\Gamma_R-\Gamma_L)\eG + (\Pi_R-\Pi_L)\eP\big) = 0.
\end{equation}

Likewise integrating (\cref{eqn: chap 5 entropy criteria}) in the direction $\nor$ over the control volume $[-\tfrac{h}{2},\tfrac{h}{2}]\times[0,\Delta t]$ gives
\begin{equation}\label{eqn: chap 5 integral form for entropy ineq}
    \frac{1}{\Delta t}\int^{\frac{h}{2}}_{-\frac{h}{2}}\eta\left(\W\left(\frac{x}{\Delta t};\vecu_L,\vecu_R,\nor\right)\right)\dd x\dd t-\frac{h}{2\Delta t}\big(\eta(\vecu_L)+\eta(\vecu_R)\big) + \q(\vecu_R)\cdot\nor - \q(\vecu_L)\cdot\nor \leqslant 0.
\end{equation}

\subsection{Three-point schemes and the Godunov method}

\begin{figure}[ht]
\begin{center}
\begin{tikzpicture}[scale=0.8]
\draw (3.0,0.8) node[above] {$h$};
\draw [>=stealth,<->] (2.0,0.8) -- (4.0,0.8) ;
\draw [>=stealth,->] (-1.0,0) -- (7.0,0) ;
\draw (7.0,0) node[below right] {$x$};
\draw (0.0,-0.1) -- (0.0,0.1);
\draw (2.0,-0.1) -- (2.0,0.1);
\draw (4,-0.1) -- (4,0.1);
\draw (6.0,-0.1) -- (6.0,0.1);
\draw (3.0,0)   node[circle,fill,inner sep=1pt](a){};
\draw (2.0,0)   node[above] {$x_{j-\frac{1}{2}}$};
\draw (3.0,0)   node[above] {$x_{j}$};
\draw (4.0,0)   node[above] {$x_{j+\frac{1}{2}}$};
\draw (1.0,-0.25) node[below] {${j-1}$};
\draw (3.0,-0.25) node[below] {${j}$};
\draw (5.0,-0.25) node[below] {${j+1}$};
\end{tikzpicture}
\caption{Notations for the mesh used for the three-point scheme (\cref{eq:3pt-scheme-NC}).}
\label{fig:stencil_1D}
\end{center}
\end{figure}
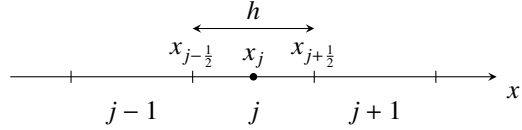

It will be convenient for the analysis of the HLLC solver to consider one-dimensional three-point numerical schemes in fluctuation form \cite{pares2006numerical}
\begin{equation}\label{eq:3pt-scheme-NC}
  {\bf U}_j^{n+1} - {\bf U}_j^{n} + \frac{\Delta t}{h}\big({\bf D}^-({\bf U}_{j}^{n},{\bf U}_{j+1}^{n},{\bf n})+{\bf D}^+({\bf U}_{j-1}^{n},{\bf U}_{j}^{n},{\bf n})\big) = 0,
\end{equation}

\noindent where ${\bf D}^\pm(\cdot,\cdot,\cdot)$ are assumed to be consistent (\cref{eq:consistency_interf_fluct}). Here, ${\bf U}_{j}^{n}$ approximates the cell-averaged solution in the $j$th cell of size $h$ at time $t^{(n)}=n\Delta t$. In the Godunov method, the fluctuations are defined by solving exact Riemann problems (\cref{eq:RP_1D_gamma_model}) centered at every interface $j+\tfrac{1}{2}$ of coordinate $x_{j+\frac{1}{2}}$, between cells $j$ and $j+1$ (see \cref{fig:stencil_1D}), with ${\bf U}_{j}^{n}$ and ${\bf U}_{j+1}^{n}$ as the left and right initial data, respectively. The solution at time $t^{n+1}$ in the $j$th cell is defined as the cell-average of the exact solution at $t^{(n+1)}=t^{(n)}+\Delta t$:
\begin{align}
\U^{n+1}_j &=\frac{1}{h}\left(\int^{x_j}_{x_{j-\frac{1}{2}}} \W\left(\frac{x}{\Delta t};\U^{n}_{j-1},\U^{n}_j,\nor\right)\dd x + \int^{x_{j+\frac{1}{2}}}_{x_j} \W\left(\frac{x}{\Delta t};\U^{n}_{j},\U^{n}_{j+1},\nor\right)\dd x\right),  \label{eqn: chap 5 cell-averaged self-similar solution}\\
&= \U^{n}_j - \frac{\Delta t}{h}\left(\frac{h}{2\Delta t}\U^{n}_j - \frac{1}{\Delta t}\int^{x_j}_{x_{j-\frac{1}{2}}} \W\left(\frac{x}{\Delta t};\U^{n}_{j-1},\U^{n}_j,\nor\right)\dd x + \frac{h}{2\Delta t}\U^{n}_j - \frac{1}{\Delta t}\int^{x_{j+\frac{1}{2}}}_{x_j} \W\left(\frac{x}{\Delta t};\U^{n}_{j},\U^{n}_{j+1},\nor\right)\dd x \right), \nonumber
\end{align}

\noindent with $x_j=x_{j-\frac{1}{2}}+\tfrac{h}{2}=x_{j+\frac{1}{2}}-\tfrac{h}{2}$ (see \cref{fig:stencil_1D}), and takes the form (\cref{eq:3pt-scheme-NC}) with the fluctuations defined as
\begin{subequations}\label{eqn: chap 5 ES fluctuation fluxes}
\begin{align}
\D^-(\vecu_L,\vecu_R,\nor)&= \frac{h}{2\Delta t}\vecu_L - \frac{1}{\Delta t}\int^0_{-\frac{h}{2}}\W\left(\frac{x}{\Delta t};\vecu_L,\vecu_R,\nor\right)\dd x \label{eqn: chap 5 D- fluctuation form} \\
 &= \textbf{f}\big(\W(0;\vecu_L,\vecu_R,\nor)\big)\cdot\nor - \textbf{f}(\vecu_L)\cdot\nor + \min(u^\star,0)\big((\Gamma_R-\Gamma_L)\eG + (\Pi_R-\Pi_L)\eP\big), \nonumber\\
\D^+(\vecu_L,\vecu_R,\nor):&= \frac{h}{2\Delta t}\vecu_R - \frac{1}{\Delta t}\int^{\frac{h}{2}}_0\W\left(\frac{x}{\Delta t};\vecu_L,\vecu_R,\nor\right)\dd x \label{eqn: chap 5 D+ fluctuation form}\\
 &= \textbf{f}(\vecu_R)\cdot\nor - \textbf{f}\big(\W(0;\vecu_L,\vecu_R,\nor)\big)\cdot\nor + \max(u^\star,0)\big((\Gamma_R-\Gamma_L)\eG + (\Pi_R-\Pi_L)\eP\big). \nonumber
\end{align}
\end{subequations}

Note that $\D^\pm(\cdot,\cdot,\cdot)$ in (\cref{eqn: chap 5 ES fluctuation fluxes}) satisfy the following path-conservation property \cite{pares2006numerical}
\begin{equation}\label{eq:path_conserv_prop}
    \D^-(\vecu_L,\vecu_R,\nor)+\D^+(\vecu_L,\vecu_R,\nor) = \textbf{f}(\vecu_R)\cdot\nor - \textbf{f}(\vecu_L)\cdot\nor + u^\star(\Gamma_R-\Gamma_L)\eG + u^\star(\Gamma_R-\Gamma_L)\eP,
\end{equation}
for a path $\phib:[0,1]\times\Omega_{\GM}\times\Omega_{\GM}\rightarrow\Omega_{\GM}$ such that
\begin{equation*}
 u^\star(\Gamma_R-\Gamma_L) = \int_0^1u\big(\phib(s;\vecu_L,\vecu_R)\big)\partial_s\phib(s;\vecu_L,\vecu_R)\cdot\eG ds, \quad u^\star(\Pi_R-\Pi_L) = \int_0^1u\big(\phib(s;\vecu_L,\vecu_R)\big)\partial_s\phib(s;\vecu_L,\vecu_R)\cdot\eP ds.
\end{equation*}

The interface Riemann problems are assumed to be noninteracting through the imposition of a half CFL condition: 
\begin{equation}\label{eqn: chap 5 HLLC CFL condition}
 \frac{\Delta t}{h}\max_{j\in\mathbb{Z}} |\lambda|_{max}(\U^n_j,\U^n_{j+1},{\bf n}) \leqslant \frac{1}{2},
\end{equation}

\noindent where $|\lambda|_{max}(\vecu_L,\vecu_R,\nor)$ is an upper bound of the absolute value of the signal speeds in the Riemann problem (\cref{eq:RP_1D_gamma_model}).

Finally, invoking a Jensen's inequality in (\cref{eqn: chap 5 cell-averaged self-similar solution}) for any convex entropy function in (\cref{eqn: chap 5 entropy pair}), we obtain the following discrete entropy inequality \cite{harten1983upstream} consistent with (\cref{eqn: chap 5 entropy criteria}):
\begin{equation}\label{eqn: chap 5 Jensen entropy inequality}
	\eta(\U^{n+1}_j) \leqslant \eta(\U_j)-\frac{\Delta t}{h}\left(Q(\U^n_{j},\U^n_{j+1},\nor)-Q(\U^n_{j-1},\U^n_{j},\nor)\right),
\end{equation}

\noindent with the consistent entropy flux $Q(\vecu_L,\vecu_R,{\bf n})=\q\left(\W(0;\vecu_L,\vecu_R,\nor)\right)\cdot\nor$.


%
\subsection{HLLC Riemann solver}\label{ssec: chap 5 HLLC Riemann solver}

\begin{figure}[h]
    \centering
\tikzset{every picture/.style={line width=0.75pt}} 

\begin{tikzpicture}[x=0.75pt,y=0.75pt,yscale=-1,xscale=1]

\draw    (132,281) -- (569,281) ;
\draw [shift={(572,281)}, rotate = 180] [fill={rgb, 255:red, 0; green, 0; blue, 0 }  ][line width=0.08]  [draw opacity=0] (10.72,-5.15) -- (0,0) -- (10.72,5.15) -- (7.12,0) -- cycle    ;
\draw    (351,84) -- (351,281) ;
\draw [shift={(351,81)}, rotate = 90] [fill={rgb, 255:red, 0; green, 0; blue, 0 }  ][line width=0.08]  [draw opacity=0] (10.72,-5.15) -- (0,0) -- (10.72,5.15) -- (7.12,0) -- cycle    ;
\draw    (151,109) -- (551,109) ;
\draw    (151,109) -- (151,281) ;
\draw    (551,109) -- (551,281) ;
\draw  [dash pattern={on 4.5pt off 4.5pt}]  (400,109) -- (352,281) ;
\draw    (199,109) -- (351,281) ;
\draw    (500,109) -- (351,281) ;
\draw    (299.58,276.71) .. controls (306.2,242.28) and (327.39,219.94) .. (366.58,219.97) ;
\draw [shift={(369,220)}, rotate = 181.36] [fill={rgb, 255:red, 0; green, 0; blue, 0 }  ][line width=0.08]  [draw opacity=0] (10.72,-5.15) -- (0,0) -- (10.72,5.15) -- (7.12,0) -- cycle    ;
\draw [shift={(299,280)}, rotate = 279.21] [fill={rgb, 255:red, 0; green, 0; blue, 0 }  ][line width=0.08]  [draw opacity=0] (10.72,-5.15) -- (0,0) -- (10.72,5.15) -- (7.12,0) -- cycle    ;
\draw  [dash pattern={on 0.84pt off 2.51pt}]  (199,109) -- (199,280) ;
\draw  [dash pattern={on 0.84pt off 2.51pt}]  (500,109) -- (500,280) ;
\draw    (372.11,220.5) .. controls (399.32,225.7) and (416.75,250.21) .. (417.94,277.07) ;
\draw [shift={(418,280)}, rotate = 270] [fill={rgb, 255:red, 0; green, 0; blue, 0 }  ][line width=0.08]  [draw opacity=0] (10.72,-5.15) -- (0,0) -- (10.72,5.15) -- (7.12,0) -- cycle    ;
\draw [shift={(369,220)}, rotate = 7.59] [fill={rgb, 255:red, 0; green, 0; blue, 0 }  ][line width=0.08]  [draw opacity=0] (10.72,-5.15) -- (0,0) -- (10.72,5.15) -- (7.12,0) -- cycle    ;

\draw (133,284.4) node [anchor=north west][inner sep=0.75pt]    {$-\frac{h}{2}$};
\draw (545,285.4) node [anchor=north west][inner sep=0.75pt]    {$\frac{h}{2}$};
\draw (201,89.4) node [anchor=north west][inner sep=0.75pt]    {$s_{L}$};
\draw (481,89.4) node [anchor=north west][inner sep=0.75pt]    {$s_{R}$};
\draw (401,89.4) node [anchor=north west][inner sep=0.75pt]    {$s^{*}$};
\draw (219,180.4) node [anchor=north west][inner sep=0.75pt]    {$\mathbf{u}_{L}$};
\draw (462,180.4) node [anchor=north west][inner sep=0.75pt]    {$\mathbf{u}_{R}$};
\draw (346,287.4) node [anchor=north west][inner sep=0.75pt]    {$0$};
\draw (328,89.4) node [anchor=north west][inner sep=0.75pt]    {$\Delta t$};
\draw (343,59.4) node [anchor=north west][inner sep=0.75pt]    {$t$};
\draw (576,274.4) node [anchor=north west][inner sep=0.75pt]    {$x$};
\draw (255,231.4) node [anchor=north west][inner sep=0.75pt]    {$\Gamma _{L} ,\Pi _{L}$};
\draw (420,231.4) node [anchor=north west][inner sep=0.75pt]    {$\Gamma _{R} ,\Pi _{R}$};
\draw (183,284.4) node [anchor=north west][inner sep=0.75pt]    {$s_{L} \Delta t$};
\draw (484,284.4) node [anchor=north west][inner sep=0.75pt]    {$s_{R} \Delta t$};
\draw (298,139.4) node [anchor=north west][inner sep=0.75pt]    {$\mathbf{u}_{L}^{*}$};
\draw (407,139.4) node [anchor=north west][inner sep=0.75pt]    {$\mathbf{u}_{R}^{*}$};
\end{tikzpicture}
\caption{Wave pattern of the HLLC solver (\cref{eqn: chap 5 approx. soln. HLLC}).}
 \label{fig: chap 5 HLL-HLLC mesh}
\end{figure}
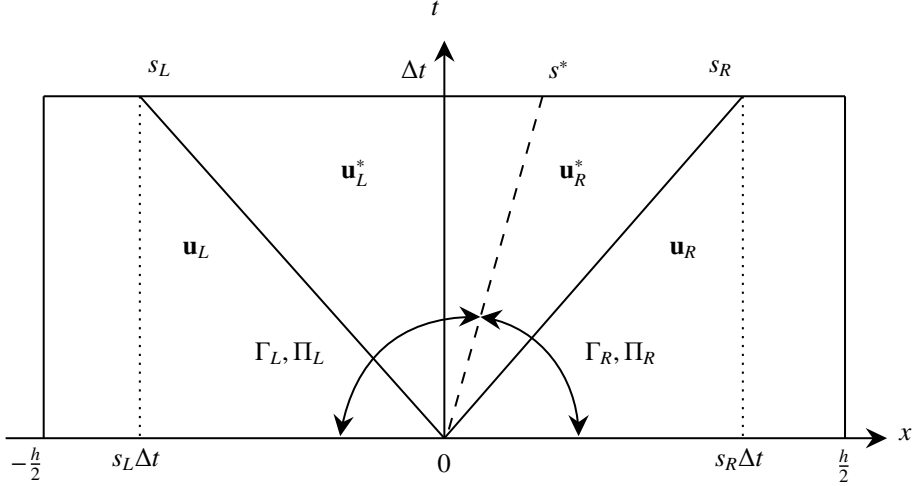

The HLLC solver \cite{toro1994restoration,batten1997choice} is a simple solver \cite{castro2017well} with four uniform states separated by simple discontinuities, see \cref{fig: chap 5 HLL-HLLC mesh}:
\begin{equation}\label{eqn: chap 5 approx. soln. HLLC}
    \W^{\HLLC}\left(\frac{x}{ t};\vecu_L,\vecu_R,\nor\right)=
    \left\{\begin{array}{ll}
    \vecu_L, & \frac{x}{t}<s_L,\\
    \vecu^\star_L, & s_L<\frac{x}{t}<s^\star,\\
    \vecu^\star_R, & s^\star<\frac{x}{t}<s_R,\\
    \vecu_R, & s_R<\frac{x}{t},
    \end{array}\right.
\end{equation}

\noindent where the wave $s^\star$ approximates the speed of the intermediate contact wave. In contrast to \cite{harten1983upstream}, we here require an approximate consistency of the HLLC solver (\cref{eqn: chap 5 approx. soln. HLLC}) with the integral form
(\cref{eqn: chap 5 exact integral GM b}). Integrating (\cref{eqn: chap 5 approx. soln. HLLC}) over $[-\tfrac{h}{2},\tfrac{h}{2}]$ at time $\Delta t$ gives
\begin{equation}\label{eqn: chap 5 split four-way HLLC solution}
\int^{\frac{h}{2}}_{-\frac{h}{2}}\W^{\HLLC}\left(\frac{x}{\Delta t};\vecu_L,\vecu_R,\nor\right)\dd x=\left(\frac{h}{2}-s_R\Delta t\right)\vecu_R + \Delta t(s_R - s^\star)\vecu^\star_R + \Delta t(s^\star-s_L)\vecu^\star_L +\left(s_L\Delta t+\frac{h}{2}\right)\vecu_L,
\end{equation}

\noindent and using (\cref{eqn: chap 5 exact integral GM b}) with $\W^{\HLLC}$ and $s^\star$ in place of $\W$ and $u^\star$, we obtain
\begin{equation}\label{eq:global_cons_HLLC}
 s_R(\vecu_R^\star-\vecu_R) + s^\star(\vecu^\star_L-\vecu^\star_R) + s_L(\vecu_L-\vecu_L^\star) + {\bf f}(\vecu_R)\cdot\nor-{\bf f}(\vecu_L)\cdot\nor + s^\star\big((\Gamma_R-\Gamma_L)\eG + (\Pi_R-\Pi_L)\eP\big) = 0.
\end{equation}

The component on mass conservation in the above relation gives
\begin{equation*}
 \rho_R(s_R-u_R) + \rho_L(u_L-s_L) = \rho_R^\star(s_R-s^\star) + \rho_L^\star(s^\star-s_L),
\end{equation*}

\noindent which is satisfied by further requiring the half-consistency conditions \cite{berthon_etal_ES_HLLC_12} which give
\begin{equation}\label{eqn: mass flux HLLC}
    Q_L=\rho_L(u_L-s_L)=\rho_L^\star(s^\star-s_L)>0, \quad Q_R=\rho_R(s_R-u_R)=\rho_R^\star(s_R-s^\star)>0,
\end{equation}

\noindent and will be referred to as the mass fluxes \cite{toro2013riemann}. Note that (\cref{eqn: mass flux HLLC}) will be shown to be important for the proof of entropy stability in \cref{ssec: chap 5 HLLC Discrete entropy inequality} and are also usually invoked through the satisfaction of some jump relations across the $s_L$ and $s_R$ waves \cite{toro1994restoration,batten1997choice} (see below).

From the approximate global consistency relation (\cref{eq:global_cons_HLLC}) we also deduce
\begin{align*}
 Q_R(u_R^\star-u_R) + Q_L(u_L^\star-u_L) &= \press_L - \press_R, \\
 Q_R(\vecv_R^\perp-\vecv_R^{\star\perp}) + Q_L(\vecv_L^\perp-\vecv_L^{\star\perp}) &= 0, \\
 Q_R(E_R^\star-E_R) + Q_L(E_L^\star-E_L) &= \press_Lu_L - \press_Ru_R, \\
 s_R(\Gamma_R-\Gamma_R^\star) + s_L(\Gamma_L^\star-\Gamma_L) &= s^\star(\Gamma_R-\Gamma_R^\star) + s^\star(\Gamma_L^\star-\Gamma_L), \\
 s_R(\Pi_R-\Pi_R^\star) + s_L(\Pi_L^\star-\Pi_L) &= s^\star(\Pi_R-\Pi_R^\star) + s^\star(\Pi_L^\star-\Pi_L), \\
\end{align*}

\noindent where $\vecv^\perp=\vecv-u\nor$ denotes the velocity component perpendicular to $\nor$. The second and two latter conditions impose
\begin{equation}\label{eqn: chap 5 gamma star eq gamma}
\vecv_L^{\star\perp}=\vecv_L^\perp, \quad \vecv_R^{\star\perp}=\vecv_R^\perp, \quad
\Gamma^\star_L=\Gamma_L, \quad \Gamma^\star_R=\Gamma_R,\quad
\Pi^\star_L=\Pi_L, \quad \Pi^\star_R=\Pi_R,
\end{equation}

\noindent meaning that $\vecv^\perp$, $\Gamma$, and $\Pi$ are continuous across shocks and may be discontinuous across $s^\star$ in agreement with the fact that $\vecv^\perp$, $\Gamma$, and $\Pi$ are associated to LD fields. The remaining unknowns can be computed by imposing the Rankine-Hugoniot relations across $s_L$ and $s_R$

\begin{equation}\label{eqn: chap 5 HLLC RH relation}
 \textbf{f}^\star_L-\textbf{f}(\vecu_L)\cdot\nor=s_L(\vecu^\star_L-\vecu_L),\quad 
 \textbf{f}^\star_R-\textbf{f}(\vecu_R)\cdot\nor=s_R(\vecu^\star_R-\vecu_R),
\end{equation}

\noindent where ${\bf f}_{X=L,R}^\star$ means that we are considering $\press^\star_X$ as an unknown instead of evaluating the pressure via the EOS (\cref{eqn: chap 5 mixture EOS}) and $\vecu_X^\star$. This leads to the following definition of the velocity in the star regions:

\begin{equation}\label{eqn: chap 5 velocity in star region}
    u^\star_L=u_L-\frac{\press^\star_L-\press_L}{Q_L}, \quad u^\star_R=u_R+\frac{\press^\star_R-\press_R}{Q_R},
\end{equation}

\noindent and to the following jump relation across the $s^\star$ wave:

 \begin{equation*}
 {\bf f}_R^\star-{\bf f}_L^\star = s^\star\big({\bf u}_R^\star-{\bf u}_L^\star-(\Gamma_R-\Gamma_L)\eG-(\Pi_R-\Pi_L)\eP\big).
\end{equation*}

Further imposing continuity of the velocity and the pressure across the intermediate wave, $u^\star_L=u^\star_R=s^\star$ and $p^\star_L=p^\star_R=p^\star$, gives the expression for the pressure and velocity in the star region:

\begin{equation}\label{eqn: chap 5 pstar HLLC}
    \press^\star=\frac{Q_L\press_R+Q_R\press_L+Q_RQ_L(u_L-u_R)}{Q_L+Q_R}, \quad s^\star=\frac{Q_Lu_L +Q_Ru_R + \press_L-\press_R}{Q_L+Q_R}.
\end{equation}

The other states are then directly obtained:

\begin{subequations}\label{eqn: chap 5 HLLC star variables}
\begin{align}
&\rho^\star_L=\frac{u_L-s_L}{s^\star-s_L}\rho_L, \;
e^\star_L=e_L+(s^\star-u_L)\left(\frac{s^\star-u_L}{2}-\frac{\press_L}{Q_L}\right), \;
E^\star_L=e^\star_L+\frac{1}{2}\vecv_L^\star\cdot\vecv_L^\star=E_L+(s^\star-u_L)\left(s^\star-\frac{\press_L}{Q_L}\right),\label{eqn: chap 5 rho starL}\\
&\rho^\star_R=\frac{s_R-u_R}{s_R-s^\star}\rho_R, \; 
e^\star_R=e_R+(s^\star-u_R)\left(\frac{s^\star-u_R}{2}+\frac{\press_R}{Q_R}\right),\;
E^\star_R=e^\star_R+\frac{1}{2}\vecv_R^\star\cdot\vecv_R^\star=E_R+(s^\star-u_R)\left(s^\star+\frac{\press_R}{Q_R}\right). \label{eqn: chap 5 rho starR}
\end{align}
\end{subequations}

As for the Godunov method in (\cref{eqn: chap 5 ES fluctuation fluxes}), we now define the fluctuations through

\begin{subequations}\label{eq:HLLC-fluct-fluxes}
\begin{align}
 {\bf D}^-({\bf u}_L,{\bf u}_R,{\bf n}) &= \frac{h}{2\Delta t}{\bf u}_L - \frac{1}{\Delta t}\int_{-\frac{h}{2}}^0\W^{\HLLC}\Big(\frac{x}{\Delta t};{\bf u}_L,{\bf u}_R,{\bf n}\Big)dx, \\
 {\bf D}^+({\bf u}_L,{\bf u}_R,{\bf n}) &= \frac{h}{2\Delta t}{\bf u}_R - \frac{1}{\Delta t}\int_0^{\frac{h}{2}}\W^{\HLLC}\Big(\frac{x}{\Delta t};{\bf u}_L,{\bf u}_R,{\bf n}\Big)dx,
\end{align}
\end{subequations}

\noindent and plugging (\cref{eqn: chap 5 approx. soln. HLLC}) into (\cref{eq:HLLC-fluct-fluxes}) gives
\begin{subequations}\label{eqn: chap 5 HLLC numerical flux}
\begin{align}
&\D^-(\vecu_L,\vecu_R,\nor)=
\left\{\begin{array}{ll}
0, & 0<s_L,\\
s_L(\vecu^\star_L-\vecu_L), & s_L<0<s^\star,\\
\textbf{f}_R^\star - \textbf{f}(\vecu_L)\cdot\nor + s^\star(\Gamma_R-\Gamma_L)\eG + s^\star(\Pi_R-\Pi_L)\eP, & s^\star<0<s_R,\\
\textbf{f}(\vecu_R)\cdot\nor - \textbf{f}(\vecu_L)\cdot\nor + s^\star(\Gamma_R-\Gamma_L)\eG + s^\star(\Pi_R-\Pi_L)\eP, & s_R<0,
\end{array}\right.\\
&\D^+(\vecu_L,\vecu_R,\nor)=
\left\{\begin{array}{ll}
\textbf{f}(\vecu_R)\cdot\nor - \textbf{f}(\vecu_L)\cdot\nor + s^\star(\Gamma_R-\Gamma_L)\eG + s^\star(\Pi_R-\Pi_L)\eP, & 0<s_L,\\
\textbf{f}(\vecu_R)\cdot\nor - \textbf{f}_R^\star + s^\star(\Gamma_R-\Gamma_L)\eG + s^\star(\Pi_R-\Pi_L)\eP, & s_L<0<s^\star,\\
s_R(\vecu_R-\vecu_R^\star), & s^\star<0<s_R,\\
0, & s_R<0.
\end{array}\right.
\end{align}
\end{subequations}

Note that by construction the numerical fluxes (\cref{eqn: chap 5 HLLC numerical flux}) satisfy the relation
\begin{equation*}
\D^-(\vecu_L,\vecu_R,\nor) + \D^+(\vecu_L,\vecu_R,\nor) = \textbf{f}(\vecu_R)\cdot\nor - \textbf{f}(\vecu_L)\cdot\nor + s^\star(\Gamma_R-\Gamma_L)\eG + s^*(\Pi_R-\Pi_L)\eP
\end{equation*}

\noindent similar to the path-conservation property (\cref{eq:path_conserv_prop}) \cite{pares2006numerical}. Note that for (\cref{eq:path_conserv_prop}) to hold, one needs $s^\star=u^\star$ which would require a root-finding algorithm to evaluate $u^\star$. We here follow another strategy where we approximate $u^\star$ by $s^\star$ in (\cref{eqn: chap 5 pstar HLLC}), this is justified by the fact that the nonconservative product is here associated to a LD field.

Finally, we define the updated cell-averaged solution as
\begin{equation}\label{eqn: chap 5 HLLC U at n+1}
 {\bf U}_j^{(n+1)} = \frac{1}{h} \int_{x_{j-\frac{1}{2}}}^{x_{j}} \W^{\HLLC}\Big(\frac{x}{\Delta t};{\bf U}_{j-1}^{n},{\bf U}_{j}^{n},\nor\Big) dx + \frac{1}{h} \int_{x_{j}}^{x_{j+\frac{1}{2}}} \W^{\HLLC}\Big(\frac{x}{\Delta t};{\bf U}_{j}^{n},{\bf U}_{j+1}^{n},\nor\Big) dx,
\end{equation}

\noindent so the HLLC solver may be recast as a three-point scheme (\cref{eq:3pt-scheme-NC}) with (\cref{eqn: chap 5 HLLC numerical flux}).

\subsection{Properties of the HLLC solver}\label{sec: chap 5 properties of the numerical scheme}

In this section, we analyse the properties of the numerical scheme (\cref{eq:3pt-scheme-NC}) using fluxes (\cref{eqn: chap 5 HLLC numerical flux}), where the time step $\Delta t >0$ is assumed to satisfy the CFL condition
\begin{equation}\label{eqn: chap 5 HLLC CFL condition}
 \frac{\Delta t}{h}\max_{j\in\mathbb{Z}} \left(\left|s_R(\U^n_j,\U^n_{j+1},{\bf n})\right|,\left|s_L(\U^n_j,\U^n_{j+1},{\bf n})\right|\right)\leqslant \frac{1}{2},
\end{equation}

\noindent where the wave speeds $s_L$ and $s_R$ will be defined in \cref{ssec: chap 5 wave speed estimates}.

\subsubsection{Discrete entropy inequality}\label{ssec: chap 5 HLLC Discrete entropy inequality}

We are here interested in the nonlinear stability of the scheme (\cref{eq:3pt-scheme-NC}) and follow \cite{berthon_etal_ES_HLLC_12} to use the local entropy minimum principles and first prove in \cref{th:entropy_ineq_HLLC} the entropy inequality in integral form (\cref{eqn: chap 5 Jensen entropy inequality}) for the HLLC solution (\cref{eqn: chap 5 approx. soln. HLLC}). We then prove in \cref{lemma: existence of wave speeds and invariance relations} that the local entropy minimum principles hold for the intermediate states.


\begin{mythm}\label{th:entropy_ineq_HLLC}
Suppose that condition (\cref{eqn: chap 5 HLLC CFL condition}) on the time step holds and that the intermediate states in the HLLC solver (\cref{eqn: chap 5 approx. soln. HLLC}), satisfy $\vecu^\star_L,\vecu^\star_R\in\Omega_{\GM}$ together with the following local minimum entropy principles
\begin{equation}\label{eqn: chap 5 intermediate entropy condition}
    s(\vecu^\star_L)\geqslant s(\vecu_L), \quad s(\vecu^\star_R)\geqslant s(\vecu_R),
\end{equation}
for the specific entropy (\cref{eqn: chap 5 mixture physical entropy}). Then, the three-point scheme (\cref{eq:3pt-scheme-NC}) satisfies an entropy inequality (\cref{eqn: chap 5 Jensen entropy inequality}) with the consistent numerical flux
\begin{equation}\label{eqn: chap 5 entropy flux for entropy ineq}
Q(\U^n_{j-1},\U^n_{j},\nor) = \q(\U^n_{j})\cdot\nor + \frac{1}{\Delta t}\int^{x_j}_{x_{j-\frac{1}{2}}}\eta\left(\W^{\HLLC}\left(\frac{x}{\Delta t};\U^n_{j-1},\U^n_{j},\nor\right)\right)\dd x -\frac{h}{2\Delta t}\eta\left(\U^{n}_{j}\right).
\end{equation}
\end{mythm}

\begin{proof}
We first prove the entropy inequality in integral form (\cref{eqn: chap 5 integral form for entropy ineq}) using (\cref{eqn: chap 5 split four-way HLLC solution}), therefore, we have
\begin{align*}
\int^{\frac{h}{2}}_{-\frac{h}{2}}\eta\left(\W^{\HLLC}\left(\frac{x}{\Delta t};\vecu_L,\vecu_R,\nor\right)\right)\dd x
=& \left(s_L\Delta t+\frac{h}{2}\right)\eta(\vecu_L)+(s^\star-s_L)\Delta t\eta(\vecu^\star_L)+(s_R-s^\star)\Delta t\eta(\vecu^\star_R)+\left(\frac{h}{2}-s_R\Delta t\right)\eta(\vecu_R),\\
\overset{(\cref{eqn: chap 5 entropy pair})}{=}& \frac{h}{2}\big(\eta(\vecu_L)+\eta(\vecu_R)\big)-\Delta t\big(s_L\rho_L s(\vecu_L) +(s^\star-s_L)\rho^\star_Ls(\vecu^\star_L)\big) \\ &- \Delta t\big((s_R-s^\star)\rho^\star_Rs(\vecu^\star_R)-s_R\rho_Rs(\vecu_R)\big),\\
\overset{(\cref{eqn: chap 5 intermediate entropy condition})}{\leqslant}& \frac{h}{2}\big(\eta(\vecu_L)+\eta(\vecu_R)\big)-\Delta t s(\vecu_L)\left(s_L\rho_L+(s^\star-s_L)\rho^\star_L\right) -\Delta t\big((s_R-s^\star)\rho^\star_R-s_R\rho_R\big)s(\vecu_R),\\
=& \frac{h}{2}\big(\eta(\vecu_L)+\eta(\vecu_R)\big) - \Delta t\big(\rho_R\vecv_Rs(\vecu_R)-\rho_L\vecv_Ls(\vecu_L)\big)\cdot\nor,
\end{align*}

\noindent where we have used the half-consistency conditions (\cref{eqn: mass flux HLLC}). As a consequence, setting $\vecu_L=\U^n_j$ and $\vecu_R=\U^n_{j+1}$, the numerical flux (\cref{eqn: chap 5 entropy flux for entropy ineq}) satisfies 
\begin{equation}\label{eqn: chap 5 HLLC discrete entropy flux inequality}
    Q(\U^n_j,\U^n_{j+1},\nor)\leqslant \q(\U^n_{j})\cdot\nor-\frac{1}{h}\int^{x_{j+\frac{1}{2}}}_{x_j}\eta\left(\W^{\HLLC}\left(\frac{x}{\Delta t};\U^n_{j},\U^n_{j+1},\nor\right)\right)\dd x + \frac{h}{2\Delta t}\eta(\U^n_j),
\end{equation}

\noindent and using (\cref{eqn: chap 5 HLLC U at n+1}) we have the following relation through Jensen's inequality for the convex entropy function (\cref{eqn: chap 5 entropy pair})
\begin{align*}
    \eta(\U^{n+1}_j) &\leqslant \frac{1}{h}\int^{x_j}_{x_{j-\frac{1}{2}}} \eta\left(\W^{\HLLC} \left(\frac{x}{\Delta t};\U^n_{j-1},\U^n_j,\nor\right) \right)\dd x + \frac{1}{h}\int^{x_{j+\frac{1}{2}}}_{x_{j}} \eta\left(\W^{\HLLC}\left(\frac{x}{\Delta t};\U^n_{j},\U^n_{j+1},\nor\right)\right)\dd x,\\
    &\underset{(\cref{eqn: chap 5 HLLC discrete entropy flux inequality})}{\overset{(\cref{eqn: chap 5 entropy flux for entropy ineq})}{\leqslant}}\frac{\Delta t}{h} \left(Q(\U^n_{j-1},\U^n_{j},\nor)-\q(\U^n_j)\cdot\nor+\frac{h}{2\Delta t}\eta(\U^n_j))\right) + \frac{\Delta t}{h}\left(-Q(\U^n_j,\U^n_{j+1},\nor)+\q(\U^n_j)\cdot\nor+\frac{h}{2\Delta t}\eta(\U^n_j)\right),\\
    &=\eta(\U^n_j) - \frac{\Delta t}{h}\left(Q(\U^n_j,\U^n_{j+1},\nor) - Q(\U^n_{j-1},\U^n_j,\nor)\right).
\end{align*}
\end{proof}

\begin{lemma}\label{lemma: existence of wave speeds and invariance relations}
There exist wave speed estimates $s_L$ and $s_R$ large enough that: (i) bound the minimum and maximum wave speeds in the exact entropy weak solution of the Riemann problem (\cref{eq:RP_1D_gamma_model}), (ii) satisfy the interlacing condition $s_L<s^\star<s_R$, (iii) ensure that the local minimum entropy principles (\cref{eqn: chap 5 intermediate entropy condition}) hold.
\end{lemma}

\begin{proof}
These results are consequences of \cite[Prop.~3.2]{berthon_etal_ES_HLLC_12} and we now show that the required assumptions on the half domains separated by $\tfrac{x}{t}=s^\star$ hold, see \cite[after Prop.~3.2]{berthon_etal_ES_HLLC_12}. This is justified here because the EOS (\cref{eqn: chap 5 mixture EOS}) does not change accross the extreme waves $s_L$ and $s_R$, hence $\tfrac{x}{t}=s^\star$ separates domains with one unique equivalent pure phase and its associated entropy (\cref{eqn: chap 5 mixture physical entropy}) which is strictly convex and satisfies the Gibbs principle (\cref{eq:Gibbs_pcple_math_entropy}). The required conditions across the extreme waves are: first, the half-consistency relations hold through (\cref{eqn: mass flux HLLC}); then, the following quantities must be invariant across the $s_L$ and $s_R$ waves \cite[\S~4.2]{berthon_etal_ES_HLLC_12}:
\begin{equation*}
  \press^\star_X+\frac{Q^2_X}{\rho^\star_X} = \press_X+\frac{Q^2_X}{\rho_X},\quad
  e^\star_X-\frac{(\press^\star_X)^2}{2Q^2_X} = e_X-\frac{\press^2_X}{2Q^2_X}, \quad X=L,R. 
\end{equation*}

The first relation is a direct consequence of (\cref{eqn: chap 5 velocity in star region}), with $u_L^\star=u_R^\star=s^\star$:
\begin{equation*}
 \press_L^\star = \press_L + Q_L(u_L-s^\star), \quad \press_R^\star = \press_R + Q_R(s^\star-u_R),
\end{equation*}

\noindent and (\cref{eqn: mass flux HLLC}). For the second relation, we inject the above relations in the expression of $e_{X=L,R}^\star$ in (\cref{eqn: chap 5 HLLC star variables}) to get
\begin{equation*}
 e_X^\star - \press_X =  \frac{\press_X-\press_X^\star}{Q_X}\left(\frac{\press_X-\press_X^\star}{2Q_X} - \frac{\press_X-}{Q_X} \right)
 = -\frac{(\press_X^\star)^2-\press_X^2}{2Q_X},\quad X=L,R,
\end{equation*}

\noindent which concludes the proof.
\end{proof}

We finally link the discrete entropy inequality (\cref{eqn: chap 5 integral form for entropy ineq}) to the entropy stable character of the numerical fluctuations \cite{castro2013entropy}, so that the HLLC fluxes can be used at the interfaces in the \red{DG scheme} to prove the semi-discrete entropy inequality established in \cref{thm: HO accuracy and entropy inequality}.

\begin{corollary}
Let a three-point scheme of the form (\cref{eq:3pt-scheme-NC}) to discretize (\cref{eqn: chap 5 GM}). Then, the entropy inequality (\cref{eqn: chap 5 Jensen entropy inequality}) implies the entropy stability of the numerical fluxes in the sense of (\cref{Eqn: chap 3 entropy stable defn nc}).
\end{corollary}

\begin{proof}
The proof relies on similar arguments as the ones used in \cite[Lemma~2.8]{bouchut2004nonlinear} in the conservative setting. Let $\U^n_{j-1}=\U^n_j=\vecu_L$ and $\U^n_{j+1}=\vecu_R$, then from  (\cref{eq:3pt-scheme-NC}), we obtain $\U^{n+1}_j=\vecu_L-\tfrac{\Delta t}{h}\textbf{D}^-(\vecu_L,\vecu_R,\nor)$ and (\cref{eqn: chap 5 Jensen entropy inequality}) gives
\begin{equation*}
    \eta\left(\vecu_L-\frac{\Delta t}{h}\textbf{D}^-(\vecu_L,\vecu_R,\nor)\right) \leqslant \eta(\vecu_L)-\frac{\Delta t}{h}\left(Q(\vecu_L,\vecu_R,\nor)-\q(\vecu_L)\cdot\nor\right).
\end{equation*}
Likewise, using $\U^n_{j-1}=\vecu_L$ and $\U^n_j=\U^n_{j+1}=\vecu_R$, we get $\U^{n+1}_j=\vecu_R-\tfrac{\Delta t}{h}\textbf{D}^+(\vecu_L,\vecu_R,\nor)$ and
\begin{equation*}
    \eta\left(\vecu_R-\frac{\Delta t}{h}\textbf{D}^+(\vecu_L,\vecu_R,\nor)\right)\leqslant \eta(\vecu_R)-\frac{\Delta t}{h}\left(\q(\vecu_R)\cdot\nor-Q(\vecu_L,\vecu_R,\nor)\right).
\end{equation*}
Summing both equations and letting $\Delta t \to 0^+$ with fixed $h$, we have up to ${\cal O}(\Delta t)$:
\begin{equation*}
\eta(\vecu_L)-\frac{\Delta t}{h}\V(\vecu_L)\cdot\textbf{D}^-(\vecu_L,\vecu_R,\nor) + \eta(\vecu_R)-\frac{\Delta t}{h}\V(\vecu_R)\cdot\textbf{D}^+(\vecu_L,\vecu_R,\nor)
\leqslant \eta(\vecu_L)+\eta(\vecu_R)-\frac{\Delta t}{h}\left(\q(\vecu_R)-\q(\vecu_L)\right)\cdot\nor,
\end{equation*}
and simplifying terms gives (\cref{Eqn: chap 3 entropy stable defn nc}).
\end{proof}

Finally, as an immediate consequence of the local minimum entropy principles (\cref{eqn: chap 5 intermediate entropy condition}) and the definition of the updated solution (\cref{eqn: chap 5 HLLC U at n+1}), the HLLC solver also satisfies a discrete minimum principle on the specific physical entropy $s({\bf u})$:

\begin{equation}\label{eq:3pt-scheme-entropy-minmax_NC}
 s({\bf U}_j^{n+1})\geqslant\min\big(s({\bf U}_{j-1}^{n}),s({\bf U}_j^{n}),s({\bf U}_{j+1}^{n})\big).
\end{equation}

\subsubsection{Preservation of material interfaces and pure phases}\label{ssec: chap 5 preservation of uniform states}

We, first, prove that the three-point scheme (\cref{eq:3pt-scheme-NC}), with numerical fluxes (\cref{eqn: chap 5 HLLC numerical flux}), preserves material interfaces \cite{abgrall1996prevent} by following the same analysis as in \cref{ssec: chap 5 contact preserving numerical fluxes}. Let us assume that the left and right states satisfy $\vecv_L=\vecv_R=\vecv$ and $\press_L=\press_R=\press$, then (\cref{eqn: chap 5 pstar HLLC}) gives $\press^\star=\press$ and $s^\star=u$, respectively. As a result, the discrete requirements $\dd(\rho\vecv)=\vecv\dd\rho$ and $\dd\rho E = (\tfrac{1}{2}\vecv\cdot\vecv)\dd\rho + \press\dd\Gamma + \dd\Pi$ applied to the three-point scheme (\cref{eq:3pt-scheme-NC}) impose 
\begin{equation}
    D^\pm_{\rho u} = u D^\pm_\rho, \quad D^\pm_{\rho E} = \frac{\vecv\cdot\vecv}{2}D^\pm_\rho + \press D^\pm_\Gamma + D^\pm_\Pi,
\end{equation}

\noindent and are obviously satisfied. Now assume that $\Gamma_L=\Gamma_R$ and $\Pi_L=\Pi_R$, then $D^\pm_\Gamma = 0$ and $D^\pm_\Pi = 0$ and the fluxes reduce to the conservative HLLC solver for the Euler equations so pure phases are preserved.

\subsubsection{Positivity of the solution}\label{ssec: chap 5 positivity of solution}

Since the updated solution (\cref{eqn: chap 5 HLLC U at n+1}) is the cell-average of the superposition of approximate Riemann solutions from the HLLC solver (\cref{eqn: chap 5 approx. soln. HLLC}), and assuming that the left and right states are positive, it is sufficient to prove positivity of the intermediate states in the Riemann solution \cite{einfeldt1991godunov}. Note that the intermediate states should be evaluated from the intermediate fluxes \cite{batten1997choice}, ${\bf f}_{X=L,R}^\star$ in (\cref{eqn: chap 5 HLLC RH relation}) since we are imposing the pressure via (\cref{eqn: chap 5 pstar HLLC}), and not evaluating it from the EOS (\cref{eqn: chap 5 mixture EOS}), see \cref{ssec: chap 5 HLLC Riemann solver}. 

The proof for positivity of density in the star region can be directly stated following its definition in (\cref{eqn: chap 5 HLLC star variables}) since $s_L<u_L$, $u_R<s_R$ from the wave estimates in \cref{ssec: chap 5 wave speed estimates}, and the interlacing property $s_L<s^\star<s_R$ in \cref{lemma: existence of wave speeds and invariance relations}. According to (\cref{eq:set_of_states_GM}) hyperbolicity of (\cref{eqn: chap 5 GM}) and positivity of the solution in the star region require satisfying $\rho^\star_Xe^\star_X>\press_{\infty_X}$, $X=L,R$, which gives for the left intermediate state:
\begin{align*}
\rho^\star_L\left(E^\star_L-\frac{(s^\star)^2}{2}\right)>\press_{\infty_L}
\overset{(\cref{eqn: chap 5 HLLC star variables})}{\underset{(\cref{eqn: mass flux HLLC})}{\Leftrightarrow}} &\rho^\star_L\left(E_L+s^\star(s^\star-u_L)-(s^\star-u_L)\frac{\press_L}{\rho_L(s_L-u_L)}-\frac{(s^\star)^2}{2}\right)>\press_{\infty_L}\\
\Leftrightarrow &\rho^\star_L\left(e_L+\frac{(u_L-s^\star)^2}{2}-(s^\star-u_L)\frac{\press_L}{\rho_L(u_L-s_L)}\right)>\press_{\infty_L}\\
\overset{(\cref{eqn: chap 5 mixture EOS})}{\underset{(\cref{eqn: mass flux HLLC})}{\Leftrightarrow}} &\left(\frac{u_L-s_L}{s^\star-s_L}\right)\left(\frac{\press_L+\gamma_L\press_{\infty_L}}{\gamma_L-1}\right)+\left(\frac{u_L-s_L}{s^\star-s_L}\right)\frac{\rho_L(u_L-s^\star)^2}{2}-\left(\frac{s^\star-u_L}{s^\star-s_L}\right)\press_L-\press_{\infty_L}>0\\
\overset{s^\star>s_L}{\Leftrightarrow} &(u_L-s_L)\left(\frac{\press_L+\gamma_L\press_{\infty_L}}{\gamma_L-1}\right) + (u_L-s_L)\frac{\rho_L}{2}\sigma^2 + \sigma\press_L-(s^\star-s_L)\press_{\infty_L}>0 \\
\Leftrightarrow & (u_L-s_L)\frac{\rho_L}{2}\sigma^2+(\press_L+\press_{\infty_L})\sigma+(u_L-s_L)\left(\frac{\press_L+\press_{\infty_L}}{\gamma_L-1}\right)>0,
\end{align*}

\noindent where $\sigma = u_L-s^\star$. The sign for this inequality holds for all $\sigma \in\mathbb{R}$ if the discriminant $\mathcal{D}$ of the above quadratic equation is negative:
\begin{equation}
    \mathcal{D} = (\press_L+\press_{\infty_L})^2-2\rho_L(u_L-s_L)^2\left(\frac{\press_L+\press_{\infty_L}}{\gamma_L-1}\right) 
    \label{eqn: mathcal D}
\end{equation}

Using $\gamma_L(\press_L+\press_{\infty_L})=\rho_Lc_L^2$, $\mathcal{D}<0$ implies $s_L<u_L-\sqrt{(\gamma_L-1)/2\gamma_L}c_L$ which is satisfied by the wave speed estimates in \cref{ssec: chap 5 wave speed estimates} since $\tfrac{\gamma_L-1}{2\gamma_L}<1$. A similar result holds for the right intermediate state. Note that \red{(\ref{eqn: mathcal D})} is the same as in \cite{batten1997choice} with $\press_X+\press_{\infty_X}$ instead of $\press_X$.

Finally, let \red{us} prove positivity of $\Gamma$ and $\Pi$ through a discrete maximum principle. The discrete equation for $X\in\{\Gamma,\Pi\}$ in (\ref{eq:3pt-scheme-NC}) with fluctuations (\cref{eqn: chap 5 HLLC numerical flux}) reads
\begin{align}
    X^{n+1}_j &= X^n_j - \frac{\Delta t}{h}\left(\min(s^{*}_{j+\frac{1}{2}},0)\left(X^n_{j+1}-X^n_j\right) + \max(s^{*}_{j-\frac{1}{2}},0)\left(X^n_{j}-X^n_{j-1}\right)\right) \nonumber\\
    &=\left(1-\frac{\Delta t}{h}\left(\max(s^{*}_{j-\frac{1}{2}},0)-\min(s^{*}_{j+\frac{1}{2}},0)\right)\right)X^n_{j} - \frac{\Delta t}{h}\min(s^{*}_{j+\frac{1}{2}},0)X^n_{j+1} + \frac{\Delta t}{h}\max(s^{*}_{j-\frac{1}{2}},0)X^n_{j-1} \label{eq:3pt-scheme-max_prcples}
\end{align}

\noindent which shows that $X_j^{n+1}$ is a convex combination of the $X_{i=j\pm1,j}^n$ under the CFL condition (\cref{eqn: chap 5 HLLC CFL condition}).

\subsubsection{Wave speed estimates}\label{ssec: chap 5 wave speed estimates}

The fan of waves of the HLLC solver must contain the fan of waves of the exact Riemann problem (\cref{eq:RP_1D_gamma_model}). This is in particular required to ensure the local entropy minimum principles (\cref{eqn: chap 5 intermediate entropy condition}), see \cite[Prop.~3.2]{berthon_etal_ES_HLLC_12}. Direct wave speed estimates have been proposed that do not \red{require one to} solve the exact Riemann problem \cite{TORO_bound_wave_speed_2020,bouchut2004nonlinear}. Note that the usual estimate $S_R=-S_L=\max(|u_L|+c_L,|u_R|+c_R)$ may be wrong due to the \red{Lax entropy condition \cite{lax1973hyperbolic}} across a shock. On the other hand, the time steps can be affected by overestimating wave speed estimates through (\cref{eqn: chap 5 HLLC CFL condition}). We propose the following wave speeds estimates
\begin{equation*}
    s_L=u_L-\tilde{c}_L, \quad s_R=u_R+\tilde{c}_R,
\end{equation*}

\noindent where $u_X=\vecv_X\cdot\nor$, $X=L,R$, and

\begin{equation*}
    \text{if}\; \press_R\geqslant\press_L: \;
    \left\{\begin{array}{l}
    \tilde{c}_L = c_L + \frac{\gamma+1}{2}\max\left(\frac{\press_R-\press_L}{\rho_Rc_R}+u_L-u_R,0\right),\\
    \tilde{c}_R = c_R +\frac{\gamma+1}{2}\max\left(\frac{\press_L-\press_R}{\rho_L\tilde{c}_L}+u_L-u_R,0\right),
    \end{array}\right. \quad \text{else}:\;
    \left\{\begin{array}{l}
    \tilde{c}_R = c_R + \frac{\gamma+1}{2}\max\left(\frac{\press_L-\press_R}{\rho_Lc_L}+u_L-u_R,0\right),\\
    \tilde{c}_L = c_L +\frac{\gamma+1}{2}\max\left(\frac{\press_R-\press_L}{\rho_R\tilde{c}_R}+u_L-u_R,0\right),
    \end{array}\right.
\end{equation*}

\noindent and $\gamma = \max(\gamma_L,\gamma_R)$, $c_X=\sqrt{\gamma(\press_X+\press_{\infty_X})/\rho_X}$ for $X=L,R$. These estimates will bound the wave speeds in the exact Riemann solution in the case of pure phases \cite{bouchut2004nonlinear}. Moreover, the above definition of $\gamma$ \red{enables one to} bound the signal speeds in the case of polytropic gases, $\press_{\infty_L}=\press_{\infty_R}=0$ \cite{renac_etal_ES_noneq_21,renac2019entropy}. \red{Although}, it is difficult to guarantee such properties in the general case when $\Gamma_L\neq\Gamma_R$ and $\Pi_L\neq\Pi_r$,  these estimates proved to be robust in the present numerical experiments.

%
%
\section{Properties of the \red{proposed DGSEM}}\label{sec: chap 5 Properties of the DGSEM scheme}

We recall here the main properties of the \red{DGSEM} proposed in this work for the discretization of the SG-gamma model (\cref{Eqn: chap 5 Cauchy prob}) with a stiffened gas EOS (\cref{eqn: chap 5 mixture EOS}).

\subsection{Semi-discrete scheme}

The semi-discrete scheme (\cref{eqn: chap 5 modified DG semi-discrete}) with the EC fluxes (\cref{eqn: chap 5 defn of EC fluxes}) in the volume integral and the HLLC flux (\cref{eqn: chap 5 HLLC numerical flux}) at interfaces satisfies a semi-discrete entropy inequality, see \cref{thm: HO accuracy and entropy inequality}. In contrast, the scheme with the CP fluxes (\cref{eqn: chap 5 material-interface preserving flux b}) in the volume integral, along with the HLLC solver at the interfaces, preserves uniform pressure and velocity profiles across material interfaces. Both are high-order accurate and preserve uniform states.

\subsection{Fully discrete scheme}

We again consider the HLLC solver (\cref{eqn: chap 5 HLLC numerical flux}) at interfaces and both CP and EC fluctuations in the volume integrals. We restrict ourselves to the use of a one-step first-order explicit time discretization. High-order time integration will be achieved by using a strong-stability preserving explicit Runge-Kutta method from \cite{shu1988efficient} that is a convex combination of forward Euler steps and thus keeps the properties of the first-order in time scheme under some condition on the time step. The fully discrete \red{DGSEM} reads

\begin{equation}\label{eq:fully-discr_DGSEM}
 \omega_i\omega_j J_\kappa^{ij} \frac{{\bf U}_\kappa^{ij,n+1}-{\bf U}_\kappa^{ij,n}}{\Delta t^{(n)}} + {\bf R}_\kappa^{ij}({\bf u}_h^{(n)}) = 0 \quad \forall \kappa\in \Omega_h, \; 0\leqslant i,j\leqslant p, \; n\geqslant 0,
\end{equation}

\noindent where $\Delta t^{(n)}=t^{(n+1)}-t^{(n)}>0$ is the time step, ${\bf U}_\kappa^{ij,n}={\bf U}_\kappa^{ij}(t^{(n)})$, ${\bf u}_h^{(n)}={\bf u}_h(\cdot,t^{(n)})$, and the vector of space residuals ${\bf R}_\kappa^{ij}(\cdot)$ is defined by (\cref{eqn: chap 5 modified DG semi-discrete}). \Cref{th:properties_discr_DGSEM} below summarizes the properties of the above scheme with fluctuations (\cref{eqn: chap 5 EC fluxes abstract form}) and (\cref{eq:def_D+_fromD-}). These properties are independent of the volume fluctuations and therefore hold for both EC and CP fluctuations for the SG-gamma model (\cref{Eqn: chap 5 shorthand GM vectors}), though the minimum entropy principle holds when excluding material interfaces. Likewise, any other interface fluctuation that satisfy properties (i) to (iv) below can be used in place of the HLLC solver. The derivation of the CFL condition will rely on the work in \cite[Lemma~3.4]{carlier_renac_IDP_22}\red{,} that proves that there exist pseudo-equilibrium states $\vecu_\kappa^{\star,n}$ in $\Omega_{\GM}$ and finite wave speed estimates $\lambda_\kappa^{\star,n}>0$ such that

\begin{equation}\label{eq:pseudo_equ_state}
 \sum_{e\in\partial\kappa}\sum_{k=0}^p\omega_kJ_e^k{\bf h}^{\Rus}\big({\bf u}_\kappa^{\star,n},{\bf u}_h^-({\bf x}_e^{k},t^{(n)}),{\bf n}_e^k\big) =  0, \quad \lambda_\kappa^{\star,n} \geqslant \max_{0\leqslant k\leqslant p,\; e\in\partial\kappa}\big(|\lambda|_{max}\Big({\bf u}_\kappa^{\star,n},{\bf u}_h^-({\bf x}_e^{k},t^{(n)}),\nor_e^k\big)\Big),
\end{equation}

\noindent where ${\bf h}^{\Rus}({\bf u}^-,{\bf u}^+,\nor)=\tfrac{1}{2}\big({\bf f}({\bf u}^-)+{\bf f}({\bf u}^+)\big)\cdot\nor-\tfrac{\lambda_\kappa^{\star,n}}{2}({\bf u}^+-{\bf u}^-)$ is the Rusanov flux and $\lambda_\kappa^{\star,n}$ bounds the fans of waves in the exact Riemann problems (\cref{eq:RP_1D_gamma_model}) with $\nor=\nor_e^k$ and left and right states ${\bf u}_\kappa^{\star,n}$ and ${\bf u}_h^-({\bf x}_e^{k},t^{(n)})$, see (\cref{eqn: chap 5 time restriction}).

\begin{mythm}\label{th:properties_discr_DGSEM}
Let us consider the \red{numerical scheme} (\cref{eq:fully-discr_DGSEM}) with consistent fluctuations (\cref{eqn: chap 5 EC fluxes abstract form}) of the form (\cref{eqn: high-order approximation proof}) in the volume integrals and consistent interface fluctuations (\cref{eq:def_D+_fromD-}) such that the associated three-point scheme (\cref{eq:3pt-scheme-NC}): (i) is robust: ${\bf U}_{j\in\mathbb{Z}}^{n\geqslant 0}\in\Omega_{\GM}$; (ii) satisfies a discrete entropy inequality (\cref{eqn: chap 5 Jensen entropy inequality}) with consistent numerical flux; (iii) satisfies a discrete minimum principle (\cref{eq:3pt-scheme-entropy-minmax_NC}) on the specific physical entropy $s(\vecu)$; (iv) satisfies discrete maximum principles (\cref{eq:3pt-scheme-max_prcples}) on $\Gamma$ and $\Pi$. Then, the updated cell-averaged solution $\langle{\bf u}_h^{(n+1)}\rangle_\kappa$ of the \red{numerical scheme} is a convex combination of DOFs at time $t^{(n)}$ and updates of three-point schemes under the following conditions on the time step:

\begin{subequations}\label{eqn: chap 5 time restriction}
\begin{align}
 \Delta t^{(n)}\max_{e\in{\cal E}_h}\max_{0\leqslant k\leqslant p}\frac{\omega_kJ_e^k}{\min(\tilde\omega_{\kappa^\pm}(\x_e^k)J_{\kappa^\pm}(\x_e^k))}\max \big(\big|s^n_{L_{j+\frac{1}{2}}}\big|,\big|s^n_{R_{j-\frac{1}{2}}}\big|,\lambda_{\kappa^-}^{\star,n},\lambda_{\kappa^+}^{\star,n}\big) &\leqslant \frac{1}{2},\label{eqn: chap 5 time restriction A} \\
 \Delta t^{(n)} \max_{\kappa\in\Omega_h}\max_{0\leqslant i,j\leqslant p} \sum_{k=0}^p \frac{\omega_k\Big( \omega_jD_{ki}\vecv_\kappa^{kj,n}\cdot\nor_{(i,k)j}+\omega_iD_{kj}\vecv_\kappa^{ik,n}\cdot\nor_{i(j,k)} - \sum\limits_{e\in\partial\kappa}\phi_\kappa^{ij}(\x_e^k)J_e^k\min\big(s^\star(\x_e^k,t^{(n)}),0\big) \Big)}{\omega_i\omega_jJ_\kappa^{ij}} &\leqslant 1, \label{eqn: chap 5 time restriction B}
\end{align}
\end{subequations}

\noindent where $\lambda_{\kappa}^{\star,n}$ is defined in (\cref{eq:pseudo_equ_state}), $\tilde\omega_{\kappa}(\x_e^k)=\tfrac{\omega_i\omega_j}{d}$ if $\x_e^k=\x_\kappa^{ij}$ is a vertex of $\kappa$, $\tilde\omega_{\kappa}(\x_e^k)=\tfrac{\omega_i\omega_j}{d-1}$ if $\x_e^k=\x_\kappa^{ij}$ is on an edge, and  $\tilde\omega_{\kappa}(\x_e^k)=\omega_i\omega_j$ if $\x_e^k=\x_\kappa^{ij}$ is
on a face.
Assuming that the DOFs $\U^{ij,n}_\kappa$ are in $\Omega_{\GM}$ for all $\kappa\in\Omega_h$ and $0\leqslant i,j\leqslant p$, the \red{scheme} (\cref{eq:fully-discr_DGSEM}) guarantees positivity of the cell-averaged solution:
\begin{equation*}
    \langle\vecu^{n+1}_h\rangle_\kappa\in\Omega_{\GM} \quad \forall\kappa\in\Omega_h,
\end{equation*}

\noindent together with a minimum principle on the specific entropy when excluding material interfaces 
\begin{equation*}
 s(\langle{\bf u}_h^{(n+1)}\rangle_\kappa) \geqslant \min \big\{s({\bf U}_\kappa^{ij,n}):\; 0\leqslant i,j\leqslant p\big\} \cup \big\{s\big({\bf u}_h^+({\bf x}_e^{k},t^{(n)})\big):\; e\in\partial\kappa, 0\leqslant k\leqslant p \big\} \quad \forall\kappa\in \Omega_h,
\end{equation*}

\noindent and maximum principles on the EOS parameters $Y$ in $\{\Gamma,\Pi\}$:
\begin{equation*}
 \min{\cal S}_\kappa(Y_h^{(n)}) \leqslant \langle Y_h^{(n+1)}\rangle_\kappa \leqslant \max{\cal S}_\kappa(Y_h^{(n)}) \; \forall\kappa\in \Omega_h, \quad  {\cal S}_\kappa(Y_h^{(n)}) := \{Y_\kappa^{ij}: 0\leqslant i,j\leqslant p\} \cup \big\{Y_h^+({\bf x}_e^{k},t): e\in\partial\kappa, 0\leqslant k\leqslant p\big\}.
\end{equation*}
\end{mythm}

\begin{proof}
From (\cref{eq:cell_averaged_DGSEM}) and (\cref{eq:fully-discr_DGSEM}), the cell-averaged discrete scheme for the $d+2$ first components of (\cref{eqn: chap 5 GM}), $Y$ in $\{\rho,\rho\vecv,\rho E\}$, reads
\begin{equation*}
 \langle Y_h^{(n+1)}\rangle_\kappa = \langle Y_h^{(n)}\rangle_\kappa - \frac{\Delta t^{(n)}}{|\kappa|}\sum_{e\in\partial\kappa}\sum_{k=0}^p\omega_kJ_e^k{\bf h}_Y^{\HLLC}\big({\bf u}_h^-({\bf x} _e^{k},t^{(n)}),{\bf u}_h^+({\bf x}_e^{k},t^{(n)}),{\bf n}_e^k\big), 
\end{equation*}

\noindent where ${\bf h}_Y^{\HLLC}(\vecu^-,\vecu^+,\nor)=\D_Y^-(\vecu^-,\vecu^+,\nor)+{\bf f}_Y(\vecu^-)\cdot\nor$ is the numerical flux in conservation form associated to the HLLC fluctuations (\cref{eqn: chap 5 HLLC numerical flux}). We first decompose $\langle Y_h^{(n)}\rangle_\kappa$ in (\cref{eqn: cell-averaged soln}) as the convex combination
\begin{align*}
 \langle Y_h^{(n)}\rangle_\kappa = \sum_{i,j=1}^{p-1}\omega_{i}\omega_j \frac{J_\kappa^{ij}}{|\kappa|} Y_\kappa^{ij,n} + \sum_{e\in\partial\kappa}\sum_{k=0}^p\tilde\omega_\kappa(\x_e^k)\frac{J_\kappa(\x_e^k)}{|\kappa|}Y_h(x_e^k,t^{(n)}),
\end{align*}

\noindent then following \cite[Sec.~4.2]{carlier_renac_IDP_22}, we add $\tfrac{\Delta t^{(n)}}{|\kappa|}$ times the flux balance in (\cref{eq:pseudo_equ_state}) to the cell-averaged scheme to get

\begin{align*}
 \langle Y_h^{(n+1)}\rangle_\kappa =& \sum_{i,j=1}^{p-1}\omega_{i}\omega_j \frac{J_\kappa^{ij}}{|\kappa|} Y_\kappa^{ij,n} + \sum_{e\in\partial\kappa}\sum_{k=0}^p\tilde\omega_\kappa(\x_e^k)\frac{J_\kappa(\x_e^k)}{|\kappa|}\bigg(Y_h^-({\bf x}_e^{k},t^{(n)}) \\
 &- \frac{\Delta t^{(n)}\omega_kJ_e^k}{\tilde\omega_\kappa(\x_e^k)J_\kappa(\x_e^k)} \Big({\bf h}_Y^{\HLLC}\big({\bf u}_h^-({\bf x}_e^{k},t^{(n)}),{\bf u}_h^+({\bf x}_e^{k},t^{(n)}),{\bf n}_e^k\big) - {\bf h}_Y^{\Rus}\big({\bf u}_\kappa^{\star,n},{\bf u}_h^-({\bf x}_e^{k},t^{(n)}),{\bf n}_e^k\big) \Big)\bigg).
\end{align*}

As a consequence, the first $d+2$ components of $\langle \vecu_h^{(n+1)}\rangle_\kappa$ are convex combinations of quantities in $\Omega_{\GM}$ under (\cref{eqn: chap 5 time restriction A}) which confirms positivity of $\langle \rho_h^{(n+1)}\rangle_\kappa$ and $\rho e(\langle\vecu_h^{(n+1)}\rangle_\kappa)$ by concavity of $\rho e(\vecu)=\rho E-\tfrac{\rho\vecv\cdot\rho\vecv}{2\rho}$. Excluding material interfaces, we also get the minimum entropy principle.

Likewise, the two last components $Y$ in $\{\Gamma, \Pi\}$ satisfy
\begin{align*}
 &\langle Y_h^{(n+1)}\rangle_\kappa\\ 
 &= \langle Y_h^{(n)}\rangle_\kappa -\frac{\Delta t^{(n)}}{|\kappa|}\Bigg( \sum^p_{i,j,k=0} \omega_i\omega_j\left(D_{ik}\vecv^{ij,n}_\kappa\cdot\nor_{(i,k)j} Y_\kappa^{kj,n}+ D_{jk}\vecv^{ij,n}_\kappa\cdot\nor_{i(j,k)}Y_\kappa^{ik,n}\right) +\sum_{e\in\partial\kappa}\sum_{k=0}^p\omega_kJ_e^k \D_Y^-\big({\bf u}_h^-({\bf x} _e^{k},t^{(n)}),{\bf u}_h^+({\bf x}_e^{k},t^{(n)}),{\bf n}_e^k\big) \Bigg).
\end{align*}

Inverting indices $i$ and $k$, then $j$ and $k$ in the volume integral, and using the form of the interface fluctuations $D_Y^-({\bf u}^-,{\bf u}^+,{\bf n})=\min(s^\star,0)(Y^+-Y^-)$, and observing that by (\cref{eq:cardinality_relation}) we have $\sum_{0\leq i,j\leq p}\sum_{e\in\partial\kappa}\sum_{k=0}^p\phi_\kappa^{ij}(\x_e^k)f(\x_e^k)=\sum_{e\in\partial\kappa}\sum_{k=0}^pf(\x_e^k)$, we get
\begin{align*}
 \langle Y_h^{(n+1)}\rangle_\kappa &= \sum_{i,j=0}^p  \left(\omega_i\omega_j \frac{J_\kappa^{ij}}{|\kappa|}  - \frac{\Delta t^{(n)}}{|\kappa|}\sum_{k=0}^p\Big(\omega_k\omega_jD_{ki}{\bf v}_\kappa^{kj,n}\cdot\nor_{(i,k)j}+\omega_i\omega_kD_{kj}{\bf v}_\kappa^{ik,n}\cdot\nor_{i(j,k)} -\sum_{e\in\partial\kappa}\omega_kJ_e^k \min\big(s^\star(\x_e^k,t^{(n)}),0\big)\Big)\right) Y_\kappa^{ij,n} \\
 &-\frac{\Delta t^{(n)}}{|\kappa|}\sum_{e\in\partial\kappa}\sum_{k=0}^p \omega_kJ_e^k \min\big(s^\star(\x_e^k,t^{(n)}),0\big)Y_h^+(\x_e^k,t^{(n)}),
\end{align*}

\noindent and invoking the metric identities in (\cref{eq:metId4nor}), $\langle Y_h^{(n+1)}\rangle_\kappa$ is indeed a	convex combination of the DOFs associated to $Y_h^{(n)}$ under (\cref{eqn: chap 5 time restriction}). Positivity and the minimum and maximum principles follow directly.
\end{proof}

Note that (\cref{eqn: chap 5 time restriction}) correspond to sharp estimates of the time step. Assuming hypercube elements of size $h$ and estimating the wave speeds by $\lambda_\infty$, (\cref{eqn: chap 5 time restriction A}) leads to $\Delta t^{(n)}\tfrac{\lambda_\infty}{h}\leqslant\tfrac{1}{2dp(p+1)}$, which results in a time step reduced by a factor $2$ when compared with \cite[Th.~3.1]{zhang2010positivity} which however assumes Cartesian meshes. Then, assuming a uniform velocity field, $\vecv_\kappa^{ij,n}=\vecv$, (\cref{eqn: chap 5 time restriction B}) leads to $\Delta t^{(n)}\tfrac{\|{\bf v}\|_\infty}{h}\leqslant\tfrac{1}{dp(p+1)}$\red{,} thus indicating that this condition will be usually satisfied by (\cref{eqn: chap 5 time restriction A}).

%
%
\section{A posteriori limiters}\label{sec: chap 5 a posterioiri limiters}

Properties of the discrete \red{DGSEM} in \cref{th:properties_discr_DGSEM} hold for the cell-averaged solution and a posteriori limiters are applied at the end of each Runge-Kutta stage to extend these properties to all DOFs within elements. Here we describe these limiters for (\cref{eq:fully-discr_DGSEM}) that are similar to the ones proposed in \cite{zhang2010maximum,zhang2010positivity,rai2021entropy,wang2012robust}. \red{As per their basic principle the DOFs are limited} at time $t^{(n+1)}$ through a linear scaling around the cell-average (\cref{eqn: cell-averaged soln}):
\begin{equation*}
    \Tilde{\U}^{ij,n+1}_\kappa = \theta_\kappa \big(\U^{ij,n+1}_\kappa - \langle \vecu_h^{(n+1)}\rangle_\kappa\big) + \langle \vecu_h^{(n+1)}\rangle_\kappa \quad \forall 0\leqslant i,j\leqslant p, \quad \kappa\in\Omega_h,
\end{equation*}

\noindent where $0\leqslant \theta_\kappa\leqslant 1$ is the limiter coefficient. We here apply successive limiters $(\theta_\kappa^{\rho},\theta_\kappa^{\rho e},\theta_\kappa^{\Gamma},\theta_\kappa^{\Pi})$ on:
\begin{itemize}
    \item mixture density:
    \begin{equation}\label{eqn: chap 5 limiter_density}
    \Tilde{\rho}^{ij,n+1}_\kappa = \theta^{\rho}_\kappa \big(\rho^{ij,n+1}_\kappa - \langle \rho_h^{(n+1)}\rangle_\kappa\big) + \langle \rho_h^{(n+1)}\rangle_\kappa, \quad 
    \theta^{\rho}_\kappa = \min \left(\frac{\langle \rho_{h}^{(n+1)}\rangle_\kappa-\epsilon}{\langle \rho_{h}^{(n+1)}\rangle_\kappa -  \rho^{\min}_{\kappa}}, 1 \right), \quad \rho^{\min}_{\kappa} = \min_{0\leqslant i,j\leqslant p} \rho^{ij,n+1}_\kappa;
    \end{equation}
		
    \item EOS parameters $\Gamma$ and $\Pi$:
    \begin{equation}\label{eqn: chap 5 limiter_gamma}
        \Tilde{Y}^{ij,n+1}_\kappa = \theta^Y_\kappa\left(Y^{ij,n+1}_\kappa - \langle Y_h^{(n+1)}\rangle_\kappa\right) + \langle Y_h^{(n+1)}\rangle_\kappa, \quad
    \theta^{Y}_\kappa = \min \left(\frac{\langle Y_{h}^{(n+1)}\rangle_\kappa-m_{Y}}{\langle Y_{h}^{(n+1)}\rangle_\kappa -  Y^{\min}_{\kappa}}, \frac{M_{Y}-\langle Y_{h}^{(n+1)}\rangle_\kappa}{Y^{\max}_{\kappa}-\langle Y_{h}^{(n+1)}\rangle_\kappa}, 1 \right), \quad Y\in\{\Gamma,\Pi\},
    \end{equation}
    where $Y^{\min}_{\kappa} = \min_{0\leqslant i,j\leqslant p} Y^{ij,n+1}_{\kappa}$, $Y^{\max}_{\kappa} = \max_{0\leqslant i,j\leqslant p} Y^{ij,n+1}_\kappa$, and
    \begin{equation*}
    m_{\Gamma} = \min_{1\leqslant i\leqslant n_s} \frac{1}{\gamma_i-1}, \quad M_{\Gamma} = \max_{1\leqslant i\leqslant n_s}\frac{1}{\gamma_i-1}, \quad
		m_{\Pi} = \min_{1\leqslant i\leqslant n_s} \frac{\gamma_i\press_{\infty_i}}{\gamma_i-1}, \quad M_{\Pi} = \max_{1\leqslant i\leqslant n_s}\frac{\gamma_i\press_{\infty_i}}{\gamma_i-1};
    \end{equation*}
		
		\item mixture total internal energy:
    \begin{equation}\label{eqn: chap 5 limiter_internal_nrj}
		 \Tilde{Y}^{ij,n+1}_\kappa = \theta^{\rho e}_\kappa \big(Y^{ij,n+1}_\kappa - \langle Y_h^{(n+1)}\rangle_\kappa\big) + \langle Y_h^{(n+1)}\rangle_\kappa, \quad Y\in\{\rho,\rho\vecv,\rho E\},
    \end{equation}
    where $\rho e(\vecu)=\rho E-\tfrac{\rho\vecv\cdot\rho\vecv}{2\rho}$ and
    \begin{equation}\label{eqn: chap 5 internal_nrj_theta}
    \theta^{\rho e}_\kappa = \min \left(\min_{0\leqslant i,j\leqslant p}\left(\frac{\rho e(\langle\vecu_h\rangle^{n+1}_\kappa)-\Tilde\press_{\infty_\kappa}^{ij,n+1}-\epsilon}{\rho e(\langle\vecu_h\rangle^{n+1}_\kappa) - \rho e^{ij,n+1}_\kappa}\right), 1 \right). 
    \end{equation}
\end{itemize}

Note that in (\cref{eqn: chap 5 limiter_density}) and (\cref{eqn: chap 5 internal_nrj_theta}), $0< \epsilon \ll 1$ is a small parameter, which we set as $\epsilon = 10^{-8}$ in our numerical tests. The limiters (\cref{eqn: chap 5 limiter_density}) and (\cref{eqn: chap 5 limiter_internal_nrj}) thus guarantee that $\Tilde{\rho}^{0\leqslant ij\leqslant p,n+1}_\kappa >0$ and $\Tilde{\rho e}^{0\leqslant ij\leqslant p,n+1}_\kappa >\Tilde\press_{\infty_\kappa}^{0\leqslant ij\leqslant p,n+1}$. Recalling (\cref{eqn: chap 5 Gamma and Pi}), we have $m_Y\leqslant Y\leqslant M_Y$, with $Y$ in $\{\Gamma,\Pi\}$, formally which may be viewed as relaxed maximum principles. The limiters in (\cref{eqn: chap 5 limiter_gamma}) thus impose similar maximum principles:
\begin{equation*}
    m_\Gamma \leq \Tilde{\Gamma}^{0\leqslant k \leqslant p, n+1}_{j} \leqslant M_\Gamma, \quad m_\Pi \leqslant \Tilde{\Pi}^{0\leqslant k \leqslant p, n+1}_{j} \leqslant M_\Pi.
\end{equation*}

%
%
\section{Numerical experiments}\label{sec: chap 5 numerical tests}

We now perform numerical tests \red{with the DGSEM for discretizing} the SG-gamma model (\cref{Eqn: chap 5 Cauchy prob})-(\cref{Eqn: chap 5 shorthand GM vectors}) where the space discretization is defined in (\cref{eqn: chap 5 modified DG semi-discrete}) and where the three-stage third-order strong stability-preserving Runge-Kutta scheme by Shu and Osher \cite{shu1988efficient} is used for the time discretization. The time step is evaluated from the CFL condition (\cref{eqn: chap 5 time restriction A}) and the limiter, introduced in \cref{sec: chap 5 a posterioiri limiters}, is applied at the end of each stage. We consider numerical tests from \cite{cheng2020discontinuous,coralic2014finite,rai2021entropy,johnsen2006implementation,de2015new}, and \red{we refer to the work in \cite{rai2021modelling} for additional details on the numerical setups}. 

The \red{DGSEM presented here} was implemented in the CFD code \textit{Aghora} developed at ONERA \cite{renac2015aghora}. We recall that here we have proposed two {numerical schemes using DGSEM} for the SG-gamma model that differ in the fluctuations (\cref{eqn: chap 3 tilde D}) in the volume integral: a CP scheme with (\cref{eqn: chap 5 material-interface preserving flux b});  a semi-discrete entropy stable scheme with EC fluxes (\cref{eqn: chap 5 defn of EC fluxes}). We will compare the \red{performance} of both schemes in this section. All numerical tests are performed at fourth-order accuracy, $p=3$, in space unless stated otherwise.


%
\subsection{Advection of a density wave}\label{ssec: density wave}

\begin{figure}[htbp]
 \center
  \subfloat[(a)]{\includegraphics[height=.25\paperwidth]{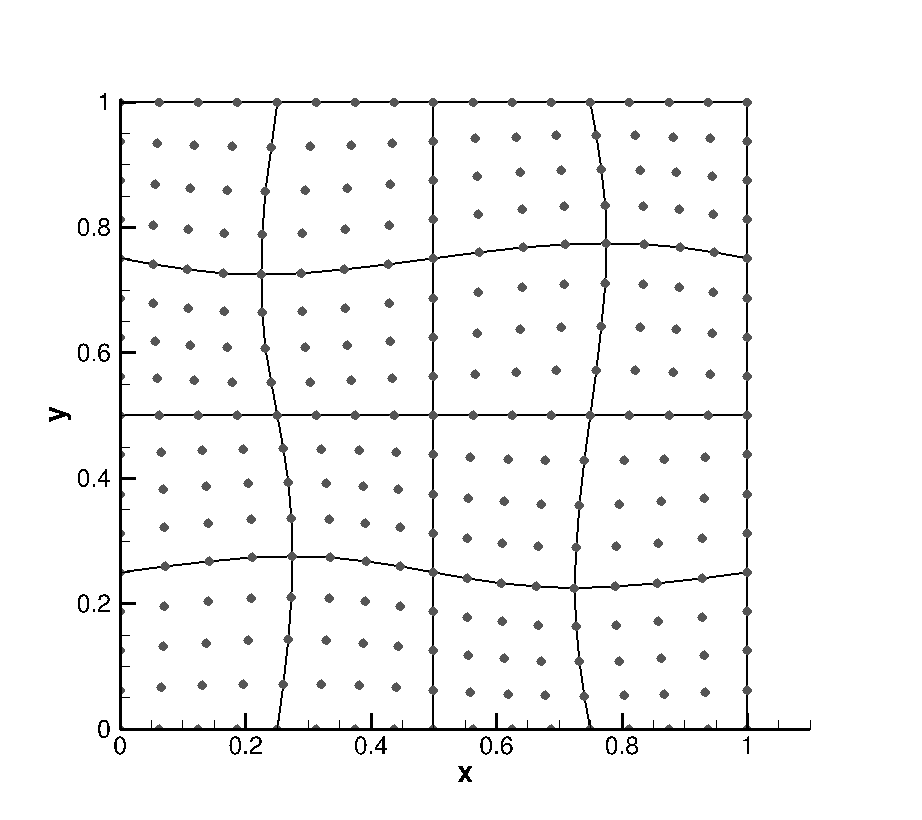}} \quad
  \subfloat[(b)]{\includegraphics[height=.25\paperwidth]{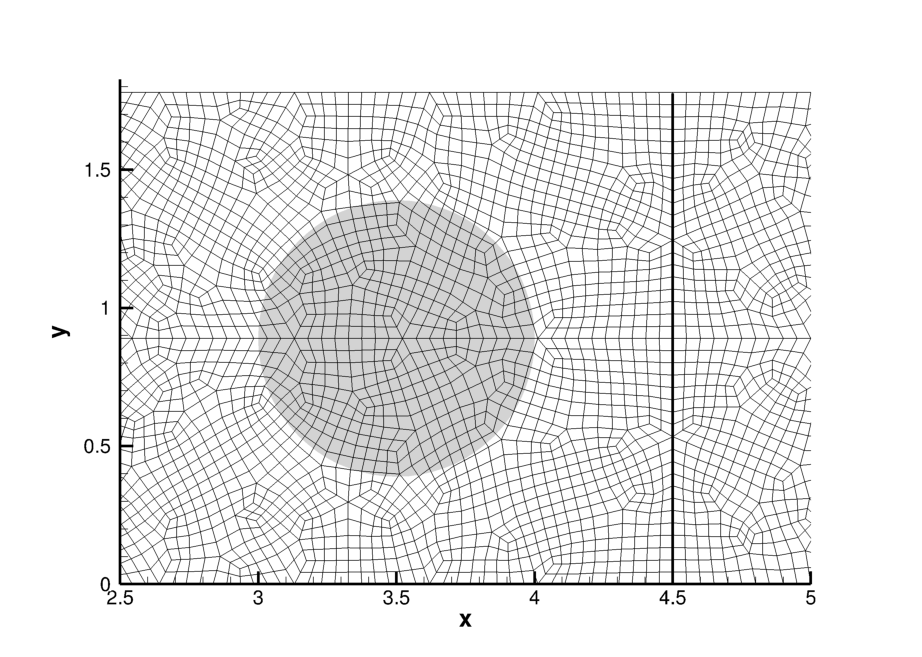}}
  \caption{Examples of meshes used for the numerical tests: (a) square mesh $h=1/4$ with $16$ fourth-order curved quadrangles with representation of the nodes; (b) example of unstructured mesh in the range $2.5\leqslant x\leqslant 5$ with $1580$ elements and initial positions of the shock and bubble (see \cref{ssec: Shock in air interacts with a helium bubble}).}
 \label{fig:test_case_meshes}
\end{figure}

We begin by validating the high-order accuracy of the scheme (\cref{eqn: chap 5 modified DG semi-discrete}) by advecting a density wave in a uniform flow in the domain $\Omega=[0,1]^2$, discretized with unstructured meshes with fourth-order curved elements (see \cref{fig:test_case_meshes}(a)), with periodic conditions and initial condition $\vecu_0(\x)$:
\begin{equation*}
 \alpha_{1_0}(\x)   = \frac{1}{2} + \frac{1}{4} \sin 4\pi (x+y),\quad
 \rho_{0}(\x)   = 1 + \frac{1}{2} \sin 2\pi(x+y),\quad
 u_{0}(\x)      = v_0(\x) = 1, \quad
 \press_{0}(\x) = 1,
\end{equation*}

\noindent along with the EOS parameters $\Cv_1=\Cv_2=1$, $\gamma_1 = 1.4$, $\press_{\infty_1} = 0$, $\gamma_2 = 3$, $\press_{\infty_2} = 2$. In this case, the density and EOS parameters are purely convected in uniform velocity and pressure fields. Norms on the mixture density error, $e_h=\rho_h-\rho$, under $h$- and $p$-refinements are displayed in \cref{tab: high-order accuracy cp_ec fluxes} with either CP fluxes, or EC fluxes in the volume integral. We observe that as the mesh is refined the expected $p+1$ order of convergence is recovered by both schemes, but lower error levels are obtained with the CP flux. 

\begin{table}[htbp]
    \centering
    \begin{tabulary}{1.0\textwidth}{ c|c|l|*{6}{c} }
        & $p$ & $h$ & $\norm{e_h}_{L^1(\Omega_h)}$ & $\mathcal{O}_1$ & $\norm{e_h}_{L^2(\Omega_h)}$ & $\mathcal{O}_2$ & $\norm{e_h}_{L^\infty(\Omega_h)}$ & $\mathcal{O}_\infty$\\ 
				\hline
				\multirow{12}{*}{CP}
        &\multirow{4}{*}{1}
        & 1/8   & 3.76E-01 & --   & 4.16E-01 & --   & 6.25E-01 & --\\
        && 1/16 & 1.82E-01 & 1.05 & 2.04E-01 & 1.03 & 3.69E-01 & 0.76\\
        && 1/32 & 5.47E-02 & 1.73 & 6.11E-02 & 1.74 & 1.08E-01 & 1.77\\
        && 1/64 & 1.43E-02 & 1.93 & 1.59E-02 & 1.94 & 2.56E-02 & 2.08\\ \cline{2-9}
        &\multirow{4}{*}{2}
        &  1/8  & 1.17E-02 & --   & 1.41E-02 & --   & 3.32E-02 & --\\
        && 1/16 & 1.06E-03 & 3.46 & 1.23E-03 & 3.52 & 3.64E-03 & 3.19\\
        && 1/32 & 9.64E-05 & 3.45 & 1.17E-04 & 3.40 & 4.82E-04 & 2.90\\
        && 1/64 & 1.01E-05 & 3.25 & 1.30E-05 & 3.16 & 6.50E-05 & 2.90\\ \cline{2-9}
        &\multirow{4}{*}{3}
        &  1/8  & 3.51E-04 & --   & 4.68E-04 & --   & 3.18E-03 & --\\
        && 1/16 & 1.73E-05 & 4.35 & 2.65E-05 & 4.14 & 1.98E-04 & 4.01\\
        && 1/32 & 1.05E-06 & 4.04 & 1.66E-06 & 4.00 & 1.20E-05 & 4.04\\
        && 1/64 & 6.73E-08 & 3.96 & 1.05E-07 & 3.98 & 7.54E-07 & 3.99\\ \hline
%
				\multirow{12}{*}{EC}
        &\multirow{4}{*}{1}
				& 1/8   & 3.72E-01 & --   & 4.11E-01 & --   & 6.24E-01 & --\\
        && 1/16 & 1.78E-01 & 1.06 & 2.05E-01 & 1.01 & 4.29E-01 & 0.54\\
        && 1/32 & 5.62E-02 & 1.66 & 7.14E-02 & 1.52 & 2.12E-01 & 1.02\\
        && 1/64 & 1.52E-02 & 1.89 & 1.99E-02 & 1.84 & 6.81E-02 & 1.64\\ \cline{2-9}
        &\multirow{4}{*}{2}
        &  1/8  & 2.78E-02 & --   & 3.48E-02 & --   & 9.99E-02 & --\\
        && 1/16 & 3.95E-03 & 2.81 & 5.23E-03 & 2.73 & 1.76E-02 & 2.50\\
        && 1/32 & 3.56E-04 & 3.47 & 4.55E-04 & 3.52 & 1.56E-03 & 3.50\\
        && 1/64 & 2.91E-05 & 3.61 & 3.77E-05 & 3.59 & 1.41E-04 & 3.46\\ \cline{2-9}
        &\multirow{4}{*}{3}
        &  1/8  & 4.58E-03 & --   & 5.78E-03 & --   & 1.62E-02 & --\\
        && 1/16 & 2.35E-04 & 4.28 & 3.10E-04 & 4.22 & 1.10E-03 & 3.88\\
        && 1/32 & 8.53E-06 & 4.78 & 1.21E-05 & 4.68 & 6.62E-05 & 4.06\\
        && 1/64 & 2.91E-07 & 4.87 & 4.63E-07 & 4.71 & 4.01E-06 & 4.05\\
%
    \end{tabulary}
    \caption{Advection of a density wave: results using either CP (top) numerical fluxes (\cref{eqn: chap 5 material-interface preserving flux b}), or EC (bottom) numerical fluxes (\cref{eqn: chap 5 defn of EC fluxes}) in the volume integral. Norms of the error on density under $p$- and $h$-refinements and associated orders of convergence at time $t = 2$.}
    \label{tab: high-order accuracy cp_ec fluxes}
\end{table}

\subsection{Riemann problems}

\begin{figure}[htbp]
 \center
  \subfloat{
  \includegraphics[height=.17\paperwidth,trim=0.2cm 0.2cm 0.2cm 0.2cm,clip=true]{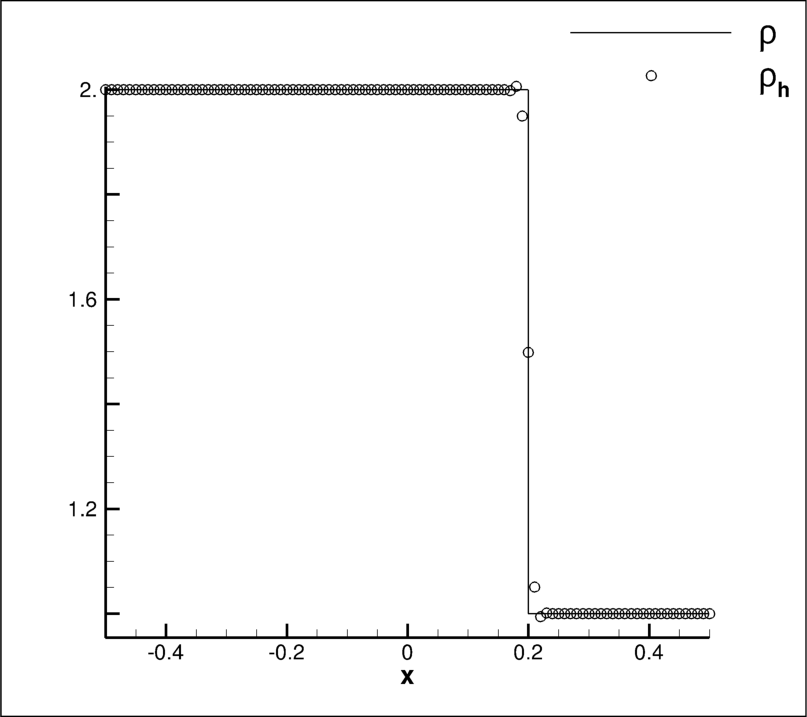}}
  \subfloat{
  \includegraphics[height=.17\paperwidth,trim=0.2cm 0.2cm 0.2cm 0.2cm,clip=true]{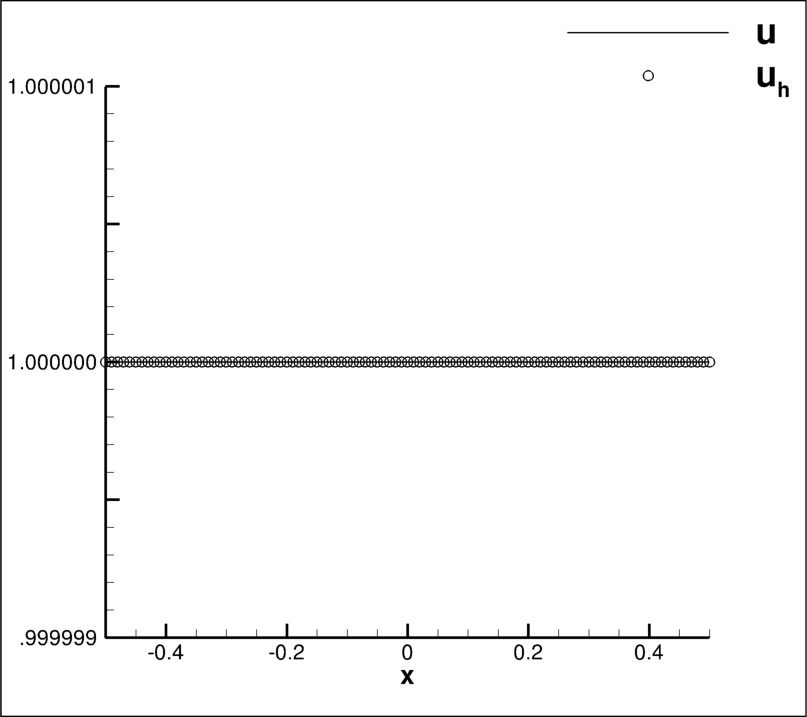}}
  \subfloat{
  \includegraphics[height=.17\paperwidth,trim=0.2cm 0.2cm 0.2cm 0.2cm,clip=true]{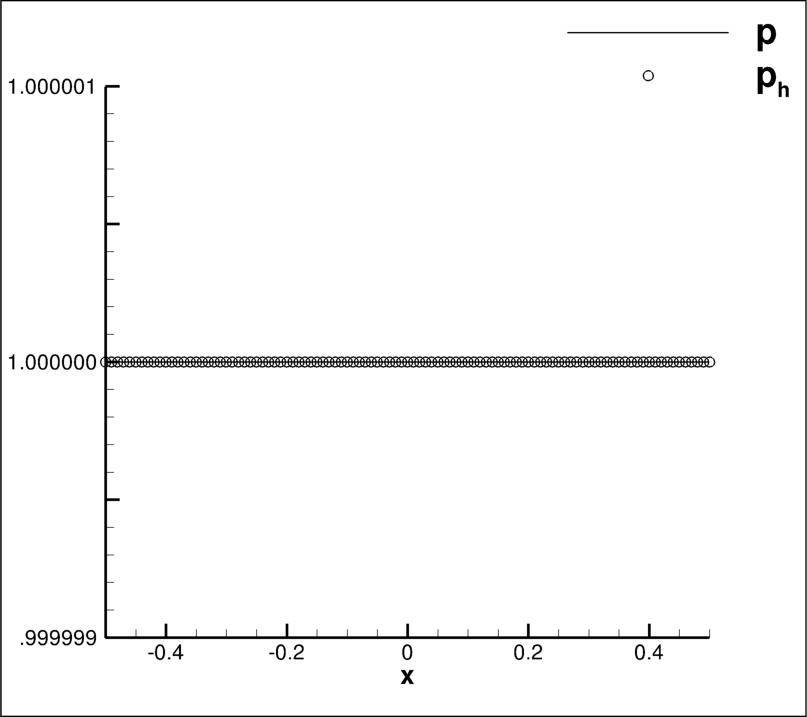}}
  \subfloat{
  \includegraphics[height=.17\paperwidth,trim=0.2cm 0.2cm 0.2cm 0.2cm,clip=true]{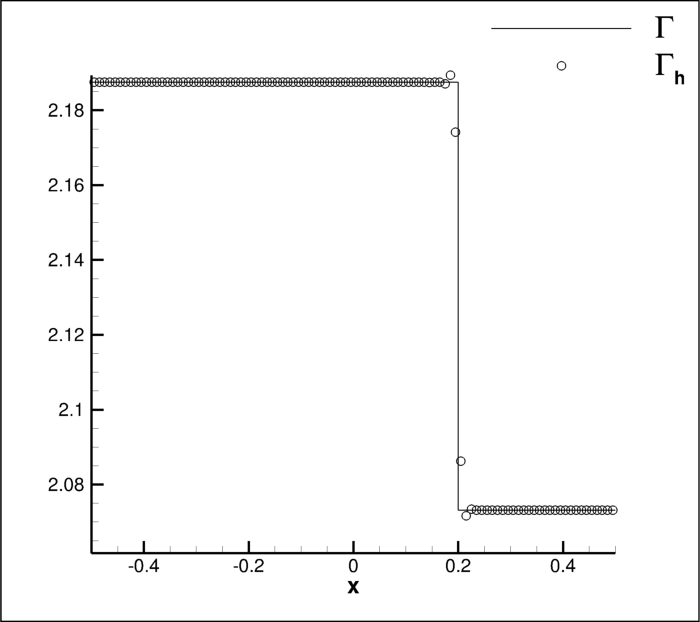}}\\
 \subfloat{
  \includegraphics[height=.17\paperwidth,trim=0.2cm 0.2cm 0.2cm 0.2cm,clip=true]{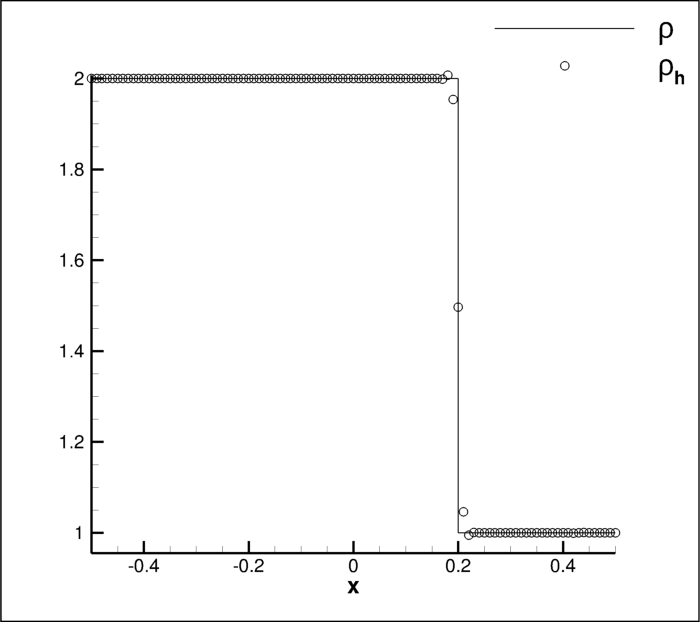}}
  \subfloat{
  \includegraphics[height=.17\paperwidth,trim=0.2cm 0.2cm 0.2cm 0.2cm,clip=true]{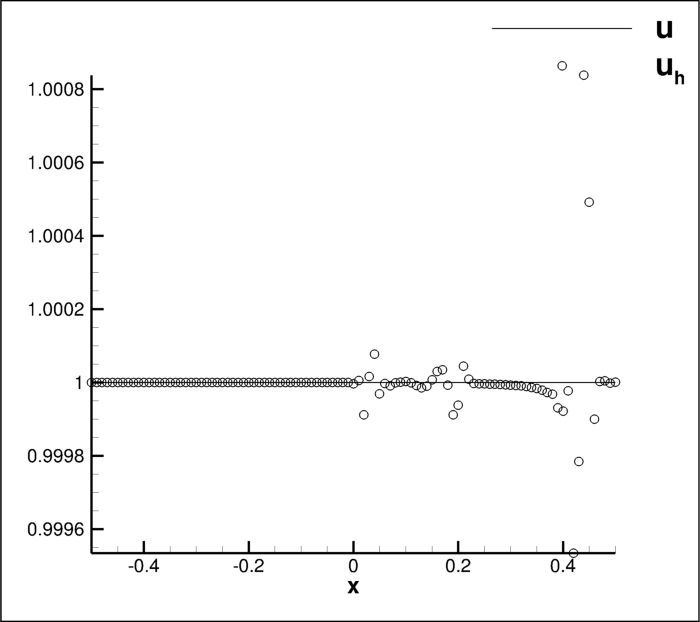}}
  \subfloat{
  \includegraphics[height=.17\paperwidth,trim=0.2cm 0.2cm 0.2cm 0.2cm,clip=true]{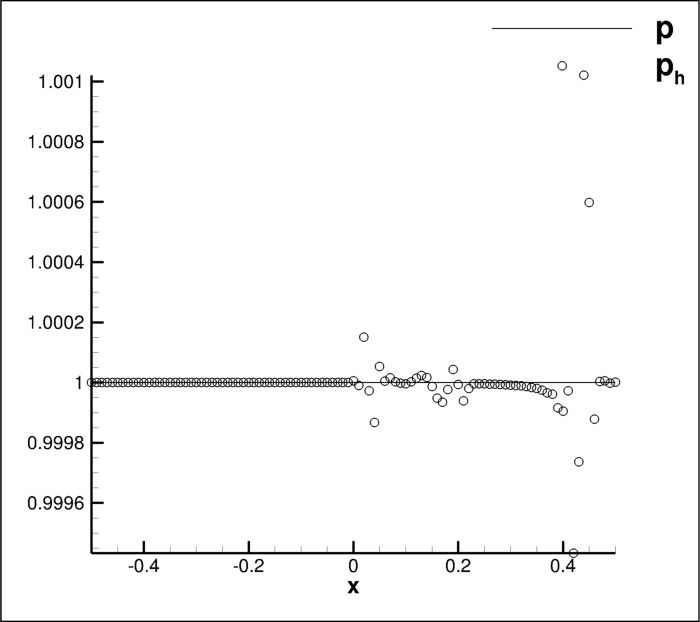}}
  \subfloat{
  \includegraphics[height=.17\paperwidth,trim=0.2cm 0.2cm 0.2cm 0.2cm,clip=true]{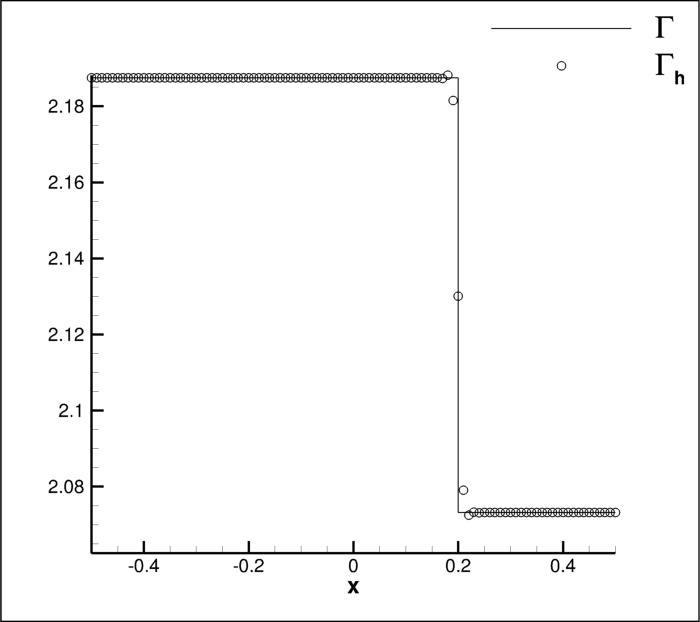}}
  \caption{Advection of an isolated material interface: fourth-order accurate simulations obtained on a mesh with $100$ elements, and using either CP (top), or EC (bottom) fluxes in the volume integral. Approximate results (symbols) are shown at $t = 0.2$ and are compared to the exact solution (lines).}
 \label{result: chap 5 contact}
\end{figure}

We first consider the advection of an isolated material discontinuity to assess the contact preservation property of the present scheme. The computational domain for this test is chosen to be $\Omega_h=[-0.5,0.5]$ with $100$ elements. The Riemann initial data are
\begin{equation*}
    (\alpha_1,\rho,u,\press)=
    \left\{\begin{array}{ll}
    (0.375,2,1,1), & x<0,\\
    (0.146342,1,1,1), & x>0,
    \end{array}\right.
\end{equation*}

\noindent with $\Cv_1 = 1$, $\Cv_2 = 2$, $\gamma_1 = 1.4$, $\gamma_2 =1.5$, $\press_{\infty_1} = \press_{\infty_2} = 0$. Results at time $t = 0.2$ are displayed in \cref{result: chap 5 contact}. As expected, computations with CP fluxes preserve uniform profiles of pressure and velocity across the material interface, while spurious oscillations occur when using EC fluxes. In both cases, the interface is well captured with low amplitude oscillations in the $\rho$ and $\Gamma$ profiles.

%

\begin{figure}[htbp]
 \center
  \subfloat{
  \includegraphics[height=.17\paperwidth,trim=0.2cm 0.2cm 0.2cm 0.2cm,clip=true]{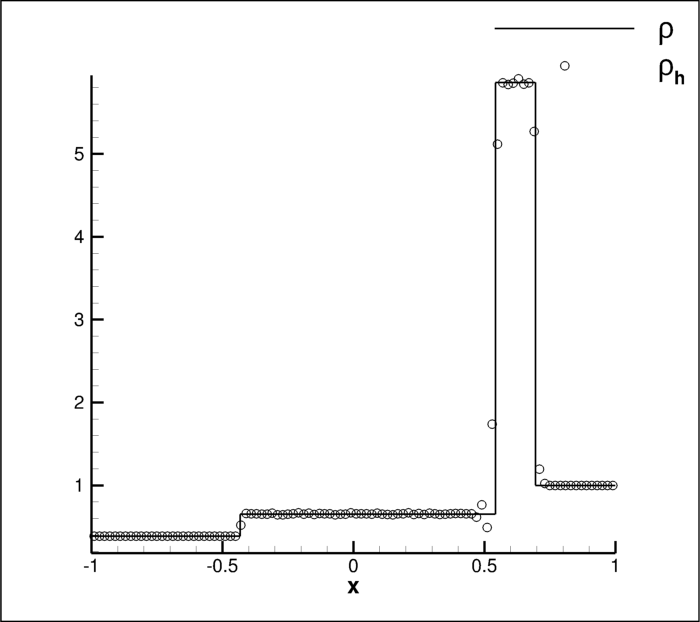}}
  \subfloat{
  \includegraphics[height=.17\paperwidth,trim=0.2cm 0.2cm 0.2cm 0.2cm,clip=true]{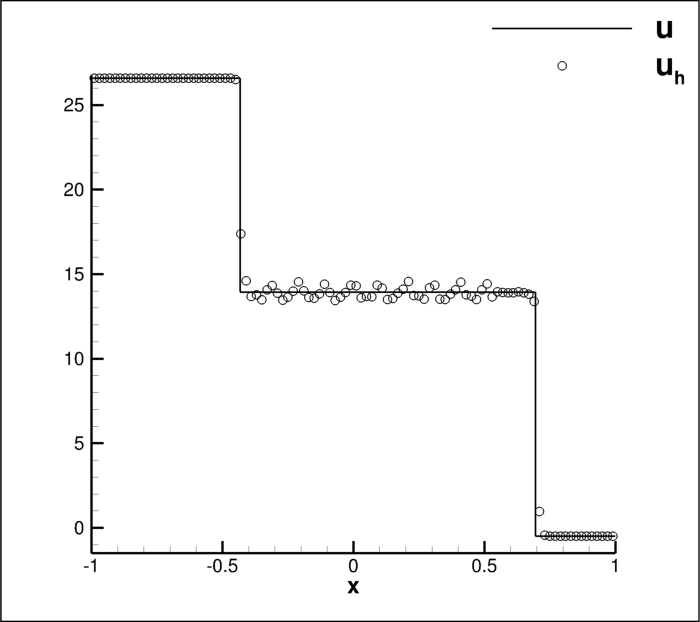}}
  \subfloat{
  \includegraphics[height=.17\paperwidth,trim=0.2cm 0.2cm 0.2cm 0.2cm,clip=true]{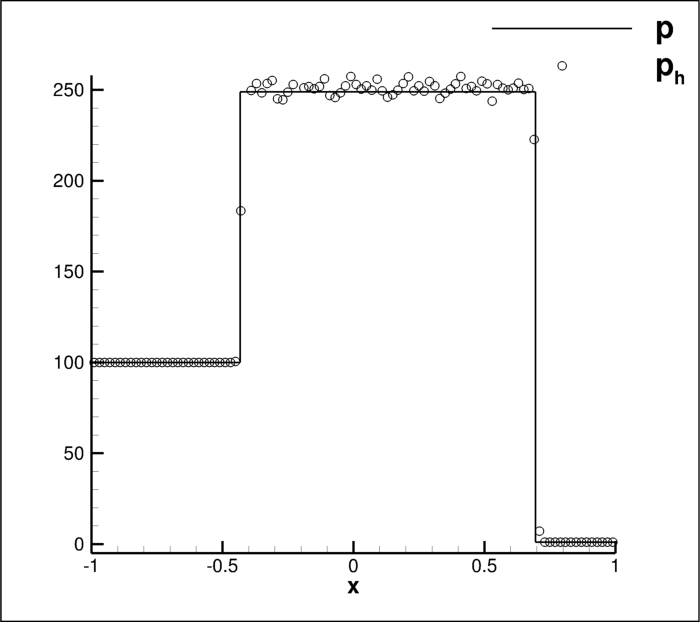}}
  \subfloat{
  \includegraphics[height=.17\paperwidth,trim=0.2cm 0.2cm 0.2cm 0.2cm,clip=true]{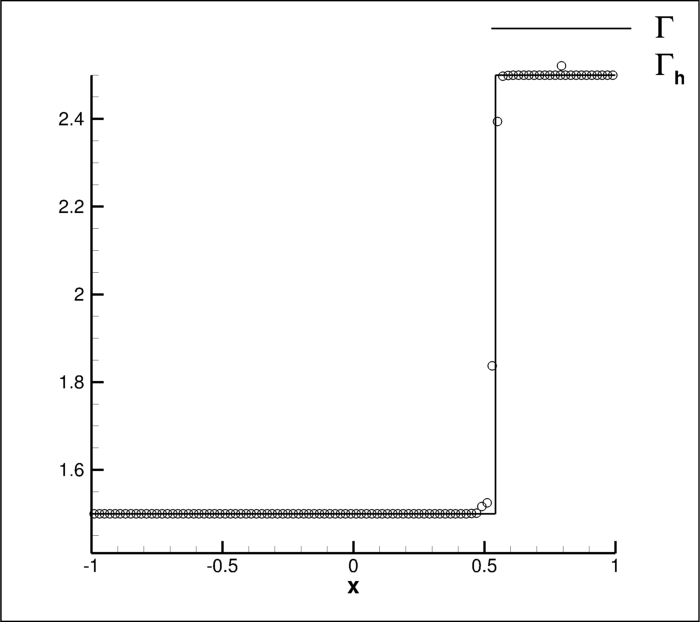}}\\
  \subfloat{
  \includegraphics[height=.17\paperwidth,trim=0.2cm 0.2cm 0.2cm 0.2cm,clip=true]{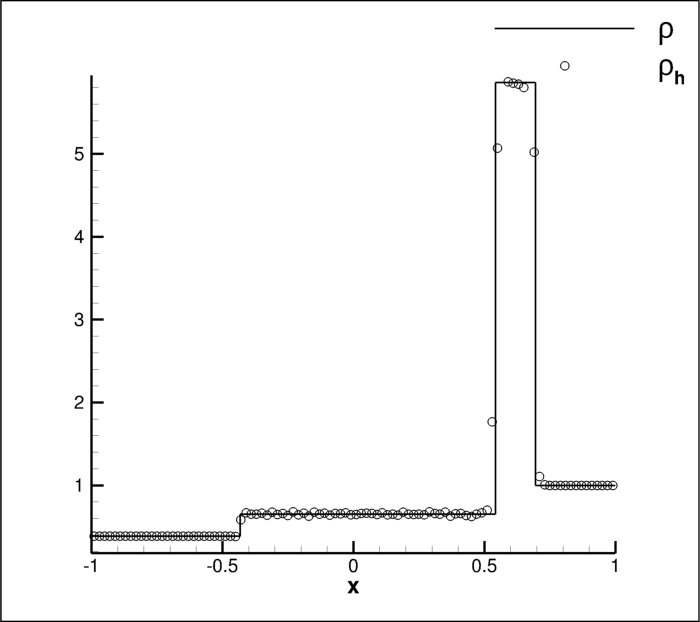}}
  \subfloat{
  \includegraphics[height=.17\paperwidth,trim=0.2cm 0.2cm 0.2cm 0.2cm,clip=true]{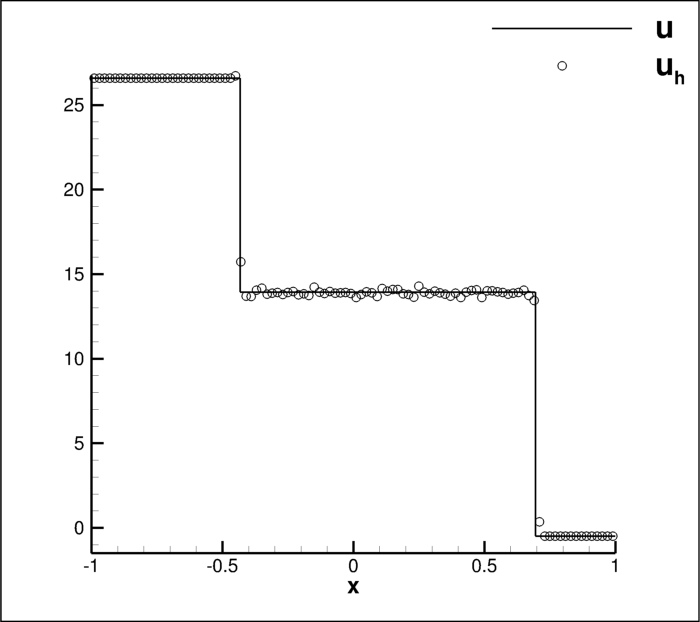}}
  \subfloat{
  \includegraphics[height=.17\paperwidth,trim=0.2cm 0.2cm 0.2cm 0.2cm,clip=true]{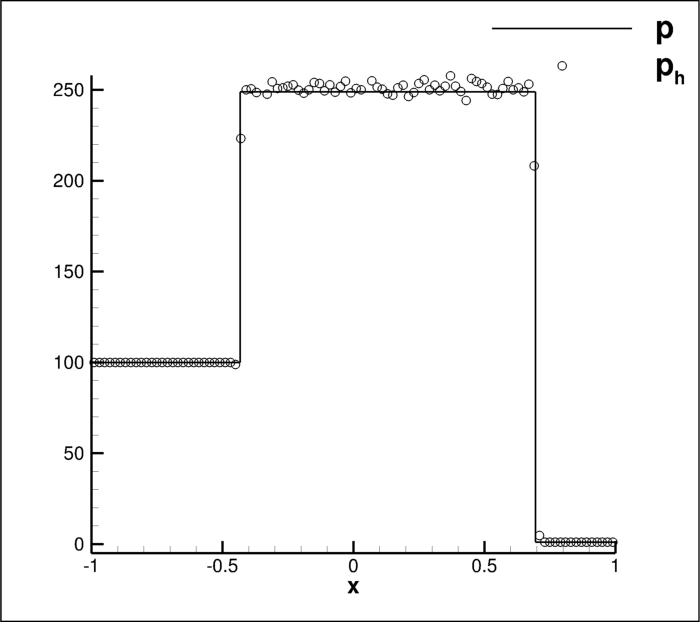}}
  \subfloat{
  \includegraphics[height=.17\paperwidth,trim=0.2cm 0.2cm 0.2cm 0.2cm,clip=true]{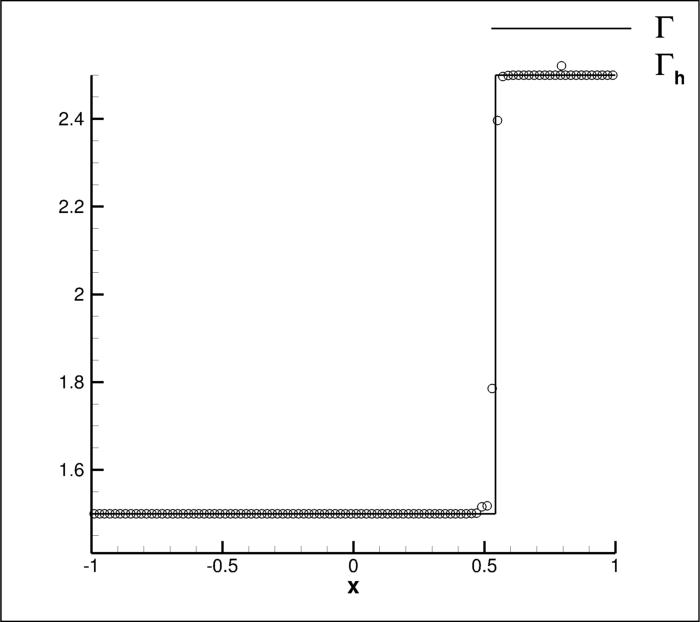}}
  \caption{Shock-material interface interaction: fourth-order accurate, $p=3$, simulations obtained on a mesh with $100$ elements, and using either CP (top), or EC (bottom) fluxes in the volume integral. Approximate results (symbols) are shown at $t = 0.07$ and are compared to the exact solution (lines).}
 \label{result: chap 5 shock-material interaction}
\end{figure}

We now consider a gas-gas shock-interface interaction problem which was originally proposed in \cite{liu2003ghost}. Here a shock wave in helium gas travels at Mach $8.96$ and interacts with an helium-air interface. The computational domain $\Omega_h = [-1,1]$ with $100$ elements and the initial condition is as follows:
\begin{equation*}
    (\alpha_1,\rho,u,\press)=
    \left\{\begin{array}{ll}
    (0,0.386,26.59,100), & x<-0.8,\\
    (0,0.1,-0.5,1), & -0.8<x<-0.2,\\
    (1,1.0,-0.5,1), & x>-0.2,
    \end{array}\right.
\end{equation*}

\noindent with $\Cv_1 = 1$, $\Cv_2 = 2.5$, $\gamma_1 = 1.4$, $\gamma_2 =5/3$, $\press_{\infty_1} = \press_{\infty_2} = 0$. The results are displayed in \cref{result: chap 5 shock-material interaction}, where the solution shows two shocks, one traveling left and the other traveling right, with a right traveling material interface in between. As a result, the shocks and interface are well captured, while spurious oscillations of small amplitude occur in the pressure and velocity fields.

\begin{figure}[htbp]
 \center
  \subfloat[]{
  \includegraphics[height=.17\paperwidth,trim=0.2cm 0.2cm 0.2cm 0.2cm,clip=true]{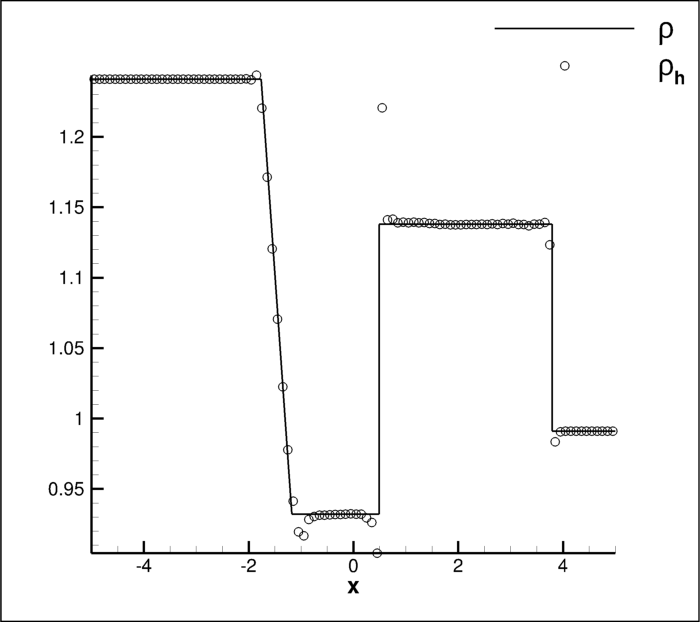}}
  \subfloat{
  \includegraphics[height=.17\paperwidth,trim=0.2cm 0.2cm 0.2cm 0.2cm,clip=true]{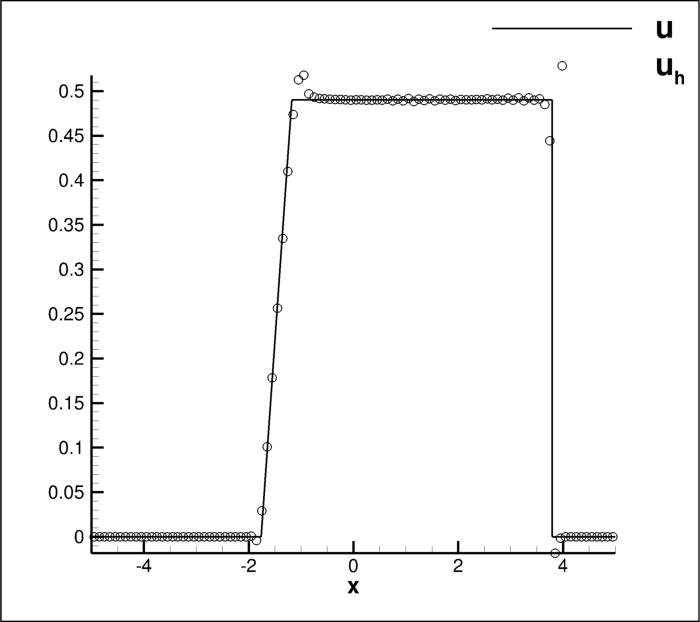}}
  \subfloat{
  \includegraphics[height=.17\paperwidth,trim=0.2cm 0.2cm 0.2cm 0.2cm,clip=true]{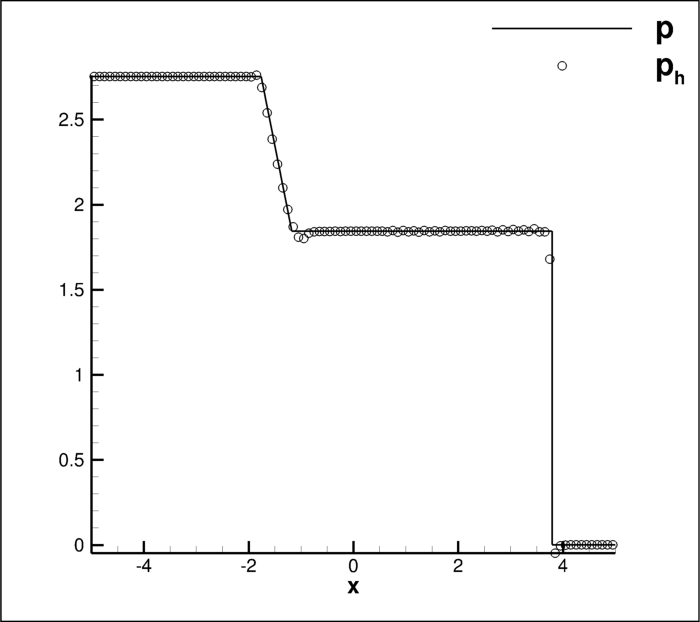}}\\
  \subfloat{
  \includegraphics[height=.17\paperwidth,trim=0.2cm 0.2cm 0.2cm 0.2cm,clip=true]{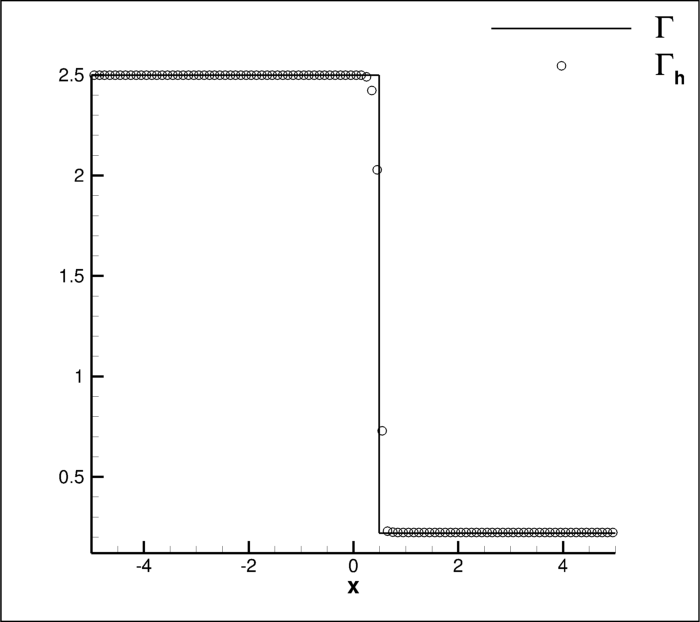}}
  \subfloat{
  \includegraphics[height=.17\paperwidth,trim=0.2cm 0.2cm 0.2cm 0.2cm,clip=true]{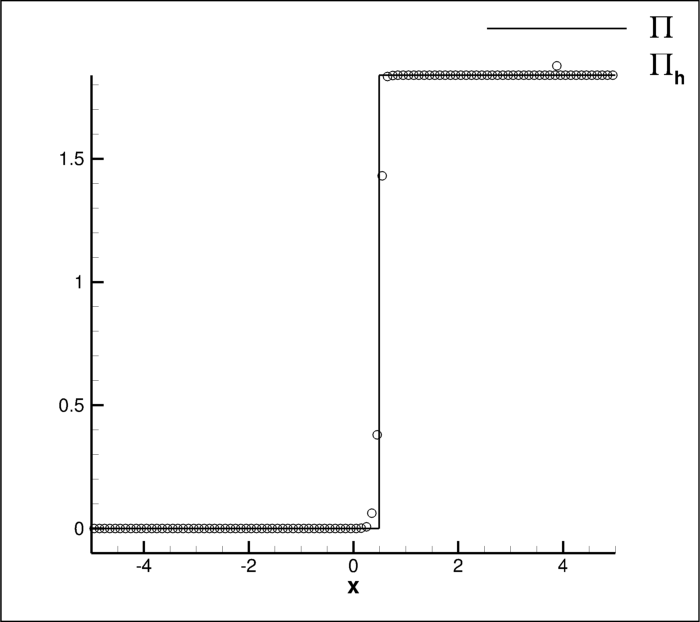}}\\
  \subfloat{
  \includegraphics[height=.17\paperwidth,trim=0.2cm 0.2cm 0.2cm 0.2cm,clip=true]{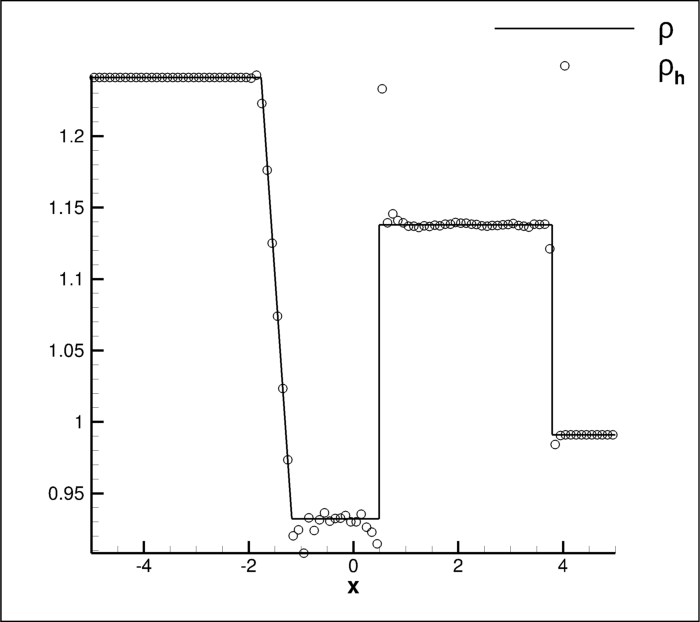}}
  \subfloat{
  \includegraphics[height=.17\paperwidth,trim=0.2cm 0.2cm 0.2cm 0.2cm,clip=true]{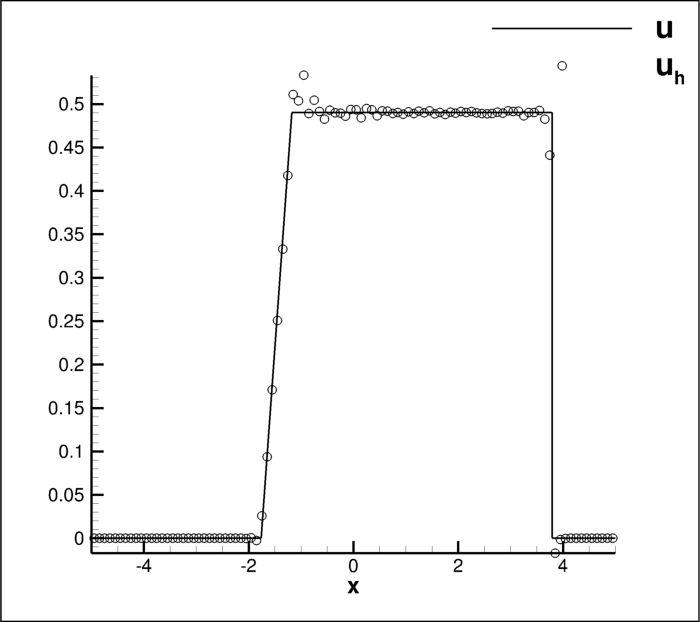}}
  \subfloat{
  \includegraphics[height=.17\paperwidth,trim=0.2cm 0.2cm 0.2cm 0.2cm,clip=true]{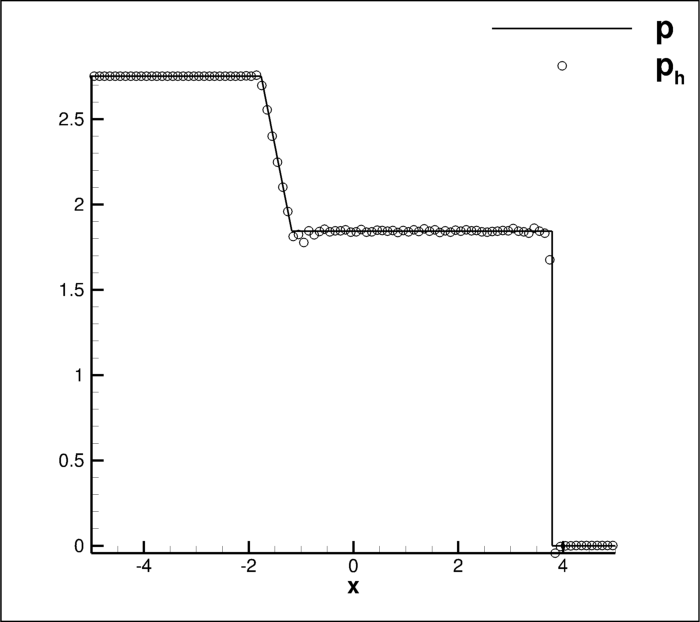}}\\
  \subfloat{
  \includegraphics[height=.17\paperwidth,trim=0.2cm 0.2cm 0.2cm 0.2cm,clip=true]{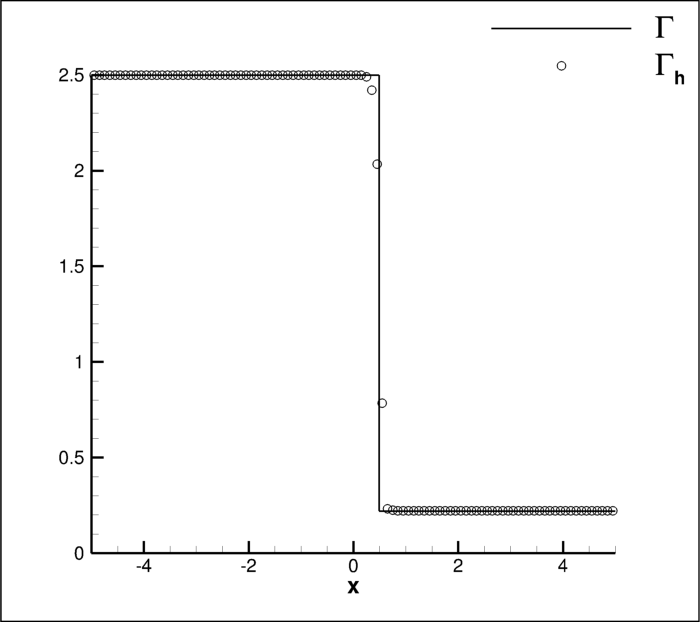}}
  \subfloat{
  \includegraphics[height=.17\paperwidth,trim=0.2cm 0.2cm 0.2cm 0.2cm,clip=true]{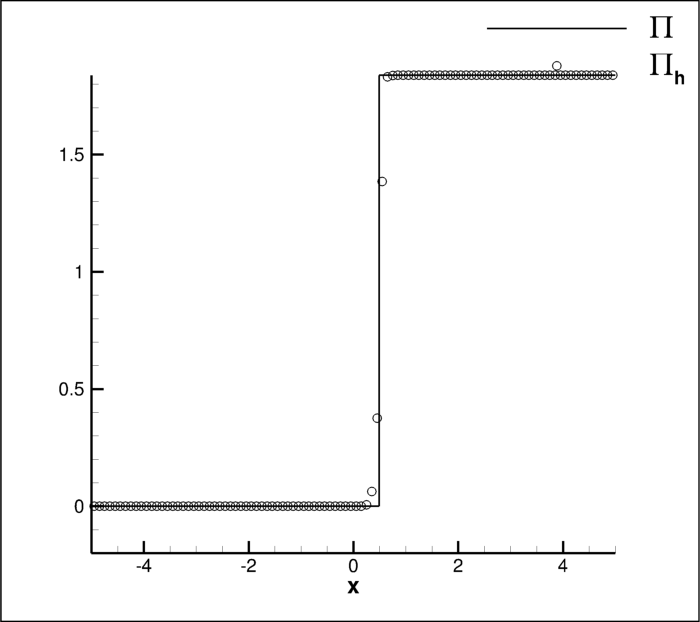}}
  \caption{Gas-water shock-interface interaction problem: fourth-order accurate, $p=3$, simulations obtained on a mesh with $100$ elements, and using either CP (top), or EC (bottom) fluxes in the volume integral. Approximate results (symbols) are shown at $t = 1$ and are compared to the exact solution (lines).}
 \label{result: chap 5 gas-water}
\end{figure}

The last test case concerns a gas-water shock-interface interaction problem and simulates an underwater explosion, where the initial condition consists of a material interface separating highly compressed air to the left and water at atmospheric pressure to the right. The computational domain is $\Omega_h=[-5,5]$ with $100$ elements, and the initial data are given as
\begin{equation*}
    (\alpha_1,\rho,u,\press)=
    \left\{\begin{array}{ll}
    (1,1.241,0,2.753), & x<0,\\
    (0,0.991,0,3.059\times 10^{-4}), & x>0,
    \end{array}\right.
\end{equation*}

\noindent with $\Cv_1 = 1.2$, $\Cv_2 = 0.073037$, $\gamma_1 = 1.4$, $\gamma_2 =5.5$, $\press_{\infty_1} = 0$, and $\press_{\infty_2} = 1.505$. The results, in \cref{result: chap 5 gas-water}, show a right traveling shock, a right traveling contact wave, a right advected material interface and a left rarefaction wave. We observe small oscillations on the velocity and pressure profiles even though the shock is of large amplitude, and the shock and material interface are well captured.

%

%
\subsection{Shock wave-helium bubble interaction}\label{ssec: Shock in air interacts with a helium bubble}

We now consider the interaction of a shock with a helium bubble, which was experimentally investigated in \cite{haas1987interaction} and used to assess numerical schemes for multiphase and multicomponent flows \cite{rai2021entropy,giordano2006richtmyer,hu2006conservative,johnsen2006implementation,kawai2011high,renac2021multicomp,quirk1996dynamics}. The test \red{setup consists of} a stationary helium bubble ($\gamma_1 = 1.648$ and $\Cv_1 = 6.0598$) which is surrounded by air ($\gamma_2 = 1.4$ and $\Cv_2=1.7857$) and interacts with a left moving Mach $1.22$ shock. The computational domain $\Omega_h=[0.0,6.5]\times[0,1.78]$ is discretized with an unstructured mesh with $433016$ elements (see \cref{fig:test_case_meshes}(b)). The helium bubble of unit diameter is centered at $x=3.5$ and $y=0.89$ and the left traveling shock is located at $x=4$. Periodic boundary conditions are imposed on the top and bottom boundaries, while non-reflective conditions are applied at the left and right boundaries. The initial data are made nondimensional with the initial bubble diameter, density, temperature and sound speed of air in the pre-shock region. 

\begin{figure}[htbp]
  \center
  \subfloat[(a) $t=0$]{\includegraphics[height=.17\paperwidth,trim=0.2cm 0.2cm 0.2cm 0.2cm,clip=true]{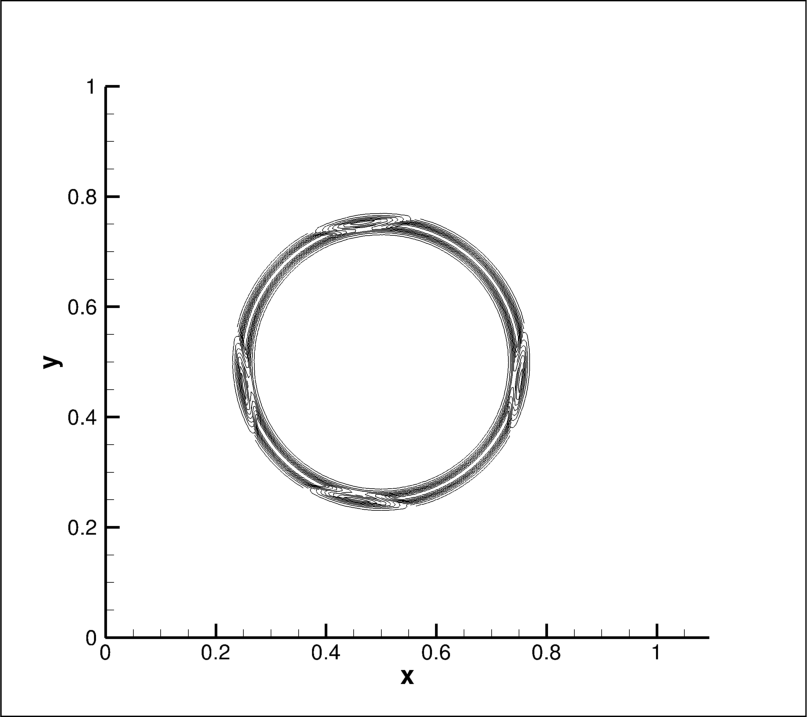}}
  \subfloat[(b) $t=76.19\mu s$, $y=0.5$]{\includegraphics[height=.17\paperwidth,trim=0.2cm 0.2cm 0.2cm 0.2cm,clip=true]{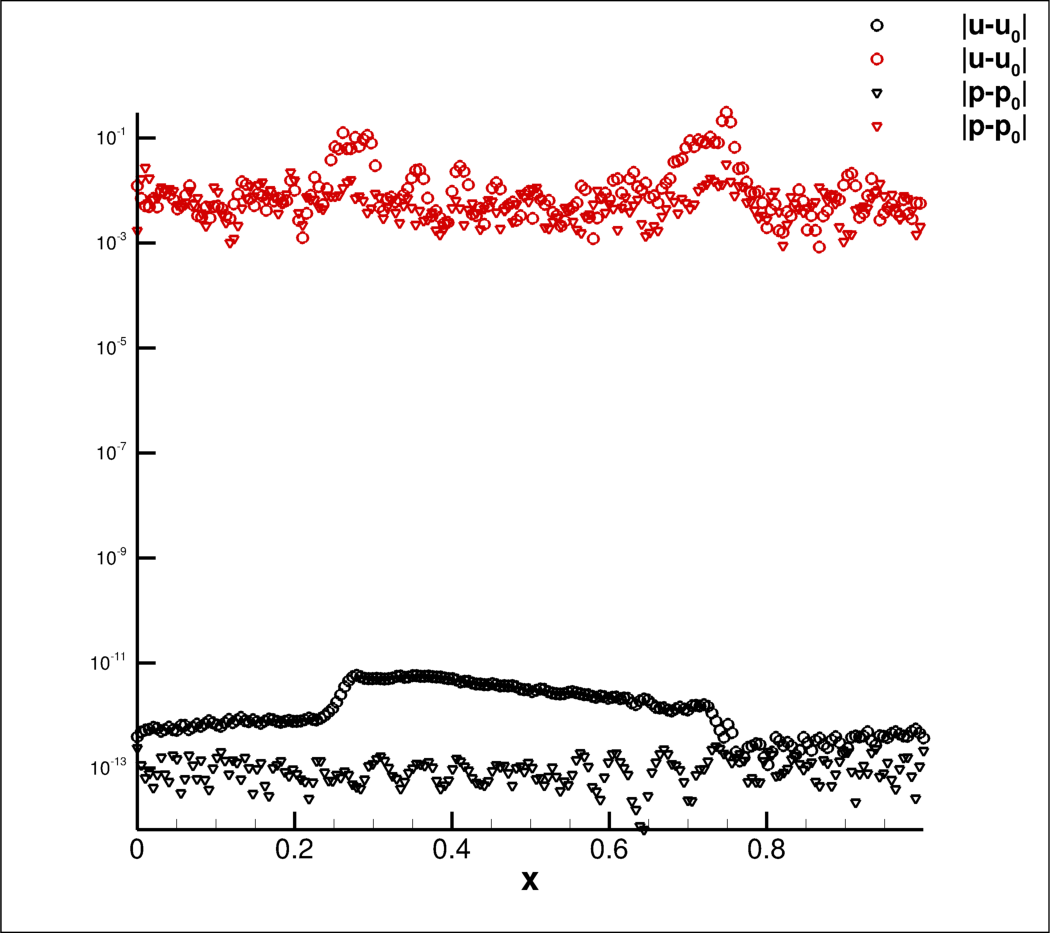}}
  \subfloat[(c) $t=76.19\mu s$, CP]{\includegraphics[height=.17\paperwidth,trim=0.2cm 0.2cm 0.2cm 0.2cm,clip=true]{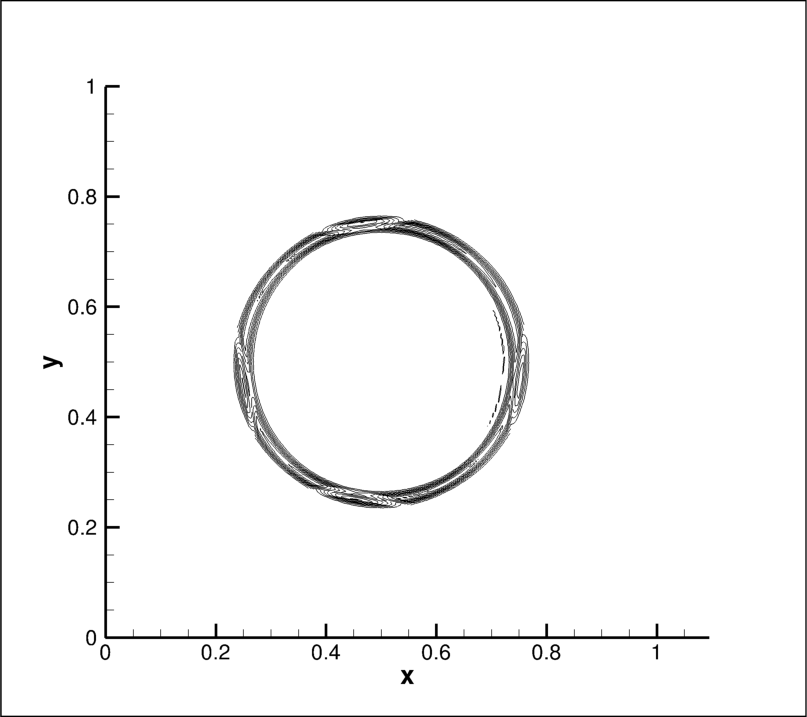}}
  \subfloat[(d) $t=76.19\mu s$, EC]{\includegraphics[height=.17\paperwidth,trim=0.2cm 0.2cm 0.2cm 0.2cm,clip=true]{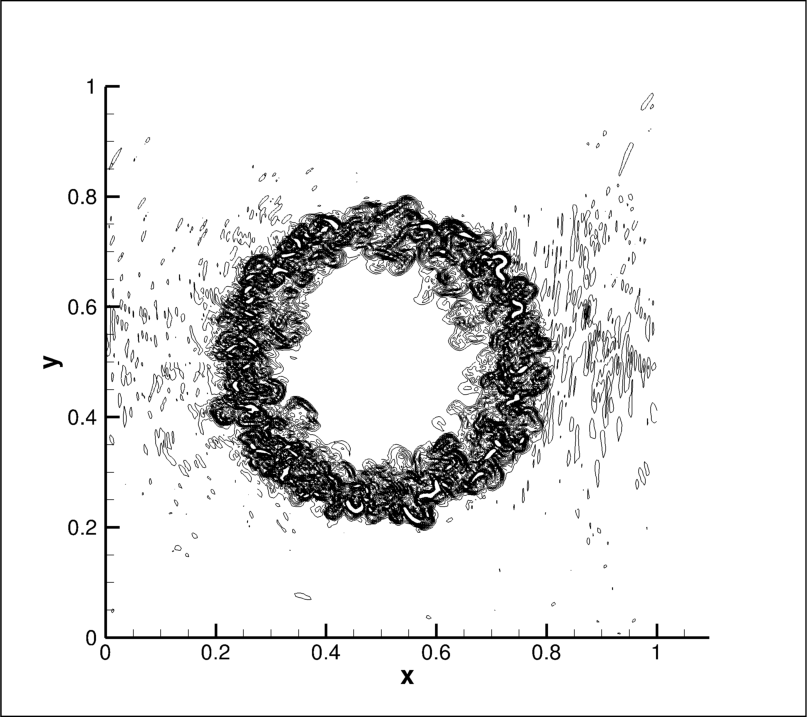}}
  \caption{Advection of a helium bubble in air: (a) initial condition, and results obtained with fourth-order accuracy $(p=3)$ in space with either CP fluxes, or EC fluxes: (b) absolute error levels (logarithmic scale) on the pressure and velocity distributions along $y=0.5$ obtained with either CP (black symbols), or EC (red symbols) fluctuations; (c,d) numerical Schlieren $\phi=\exp(|\nabla\rho|/|\nabla\rho|_{\max})$.}
  \label{result: bubble transport}
\end{figure}

We first test the ability of both schemes to preserve material interfaces and consider the advection of the bubble only on a mesh with $64\times64$ fourth-order curved elements (see \cref{fig:test_case_meshes}(a)). We thus remove the shock wave, impose a uniform velocity field $\vecv_0(\x)=(1,0)^\top$ and reduce the size of the bubble which is now initially centered at $(0.5,0.5)$. \Cref{result: bubble transport} shows the bubble at time $t=76.19\mu s$ corresponding to a transport of the bubble over a unit distance. The scheme with the CP fluxes captures the interface sharply and preserves the uniform velocity and pressure profiles across the interface, while spurious oscillations are observed with the EC flux.

\begin{figure}[htbp]
 \center
 \subfloat{\begin{picture}(0,0) \put(-10,50){\rotatebox{90}{$t=32\mu s$}} \end{picture}}
 \subfloat{\includegraphics[height=0.18\paperwidth,trim=0.2cm 0.2cm 0.2cm 0.2cm,clip=true]{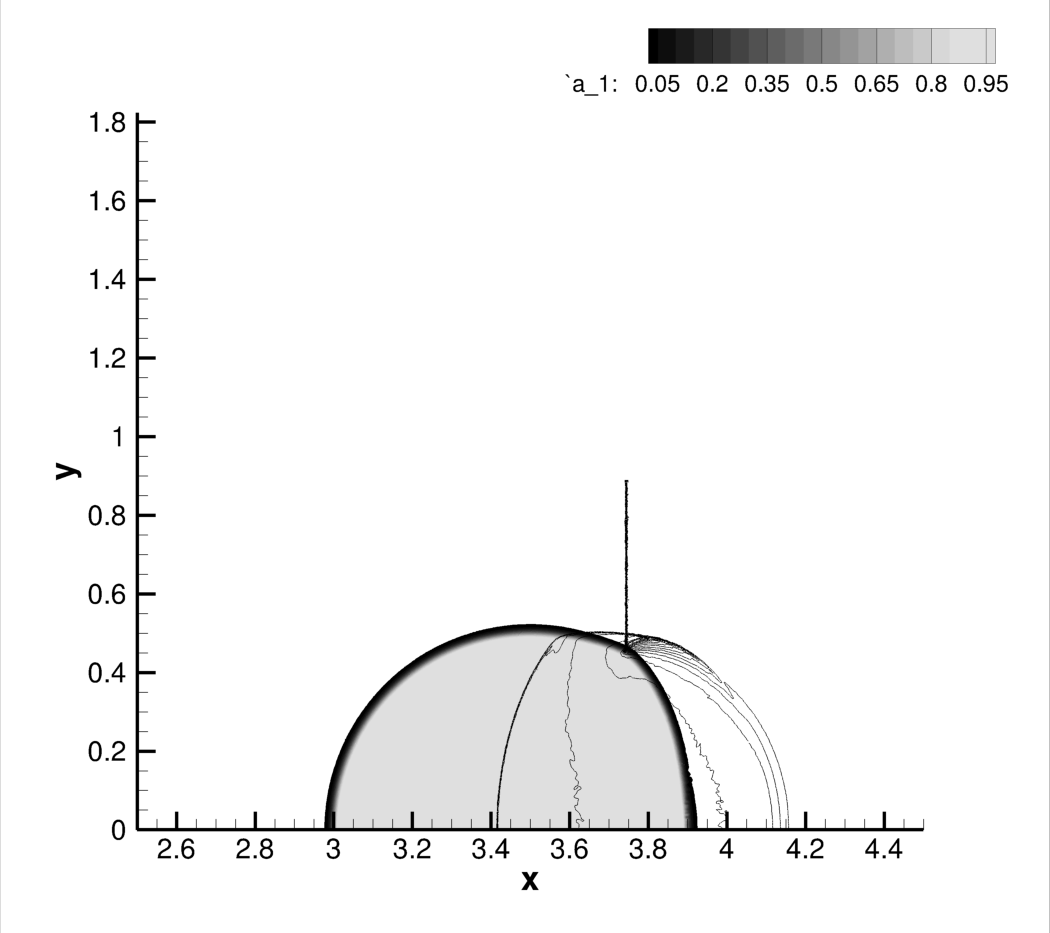}}\hspace*{-0.2cm}
  \subfloat{\includegraphics[height=0.18\paperwidth,trim=0.2cm 0.2cm 0.2cm 0.2cm,clip=true]{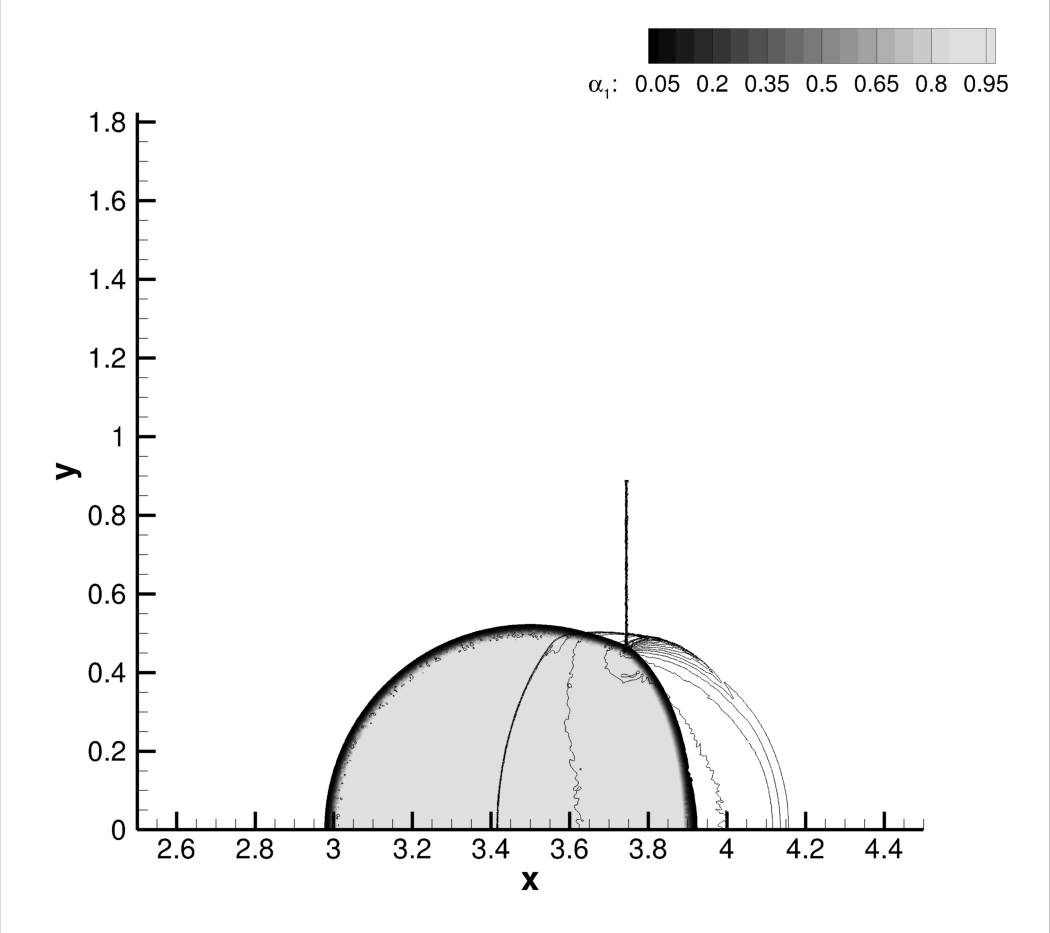}}\hspace*{-0.2cm}
  \subfloat{\includegraphics[height=0.18\paperwidth,trim=0.2cm 0.2cm 0.2cm 0.2cm,clip=true]{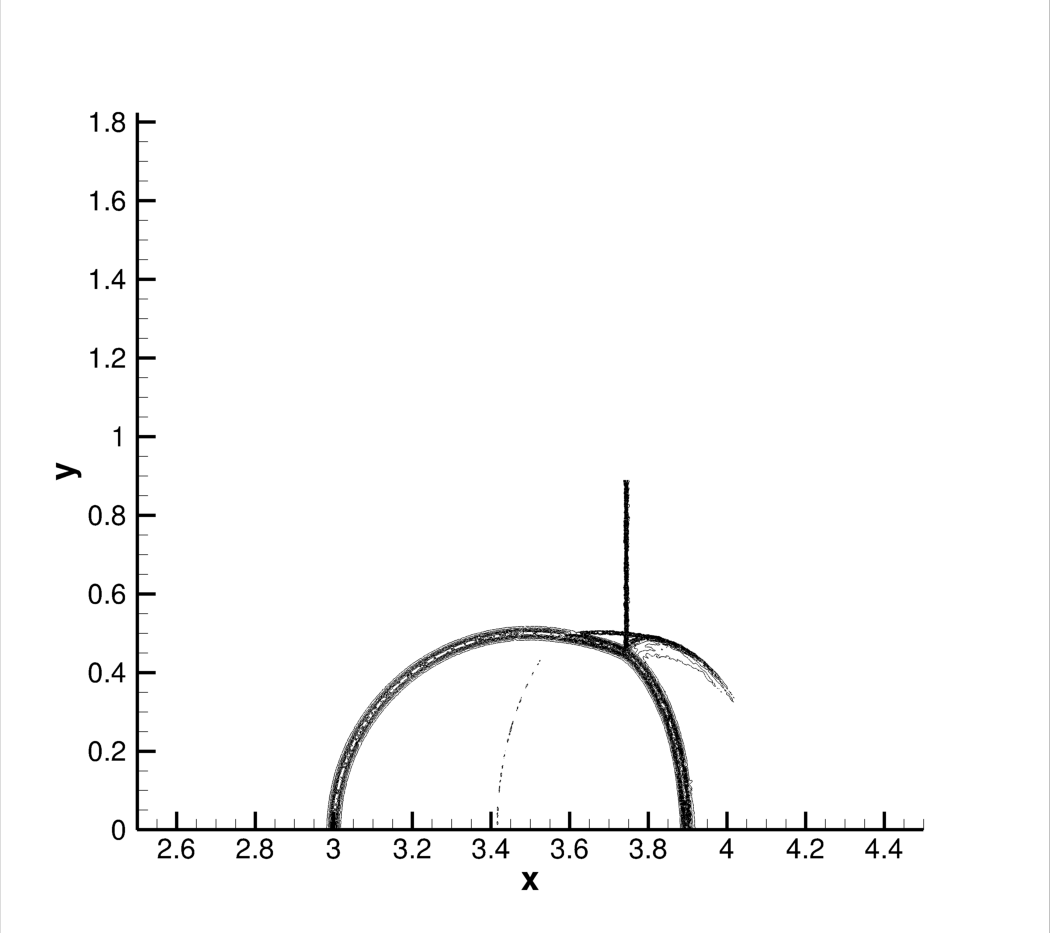}}\hspace*{-0.2cm}
  \subfloat{\includegraphics[height=0.18\paperwidth,trim=0.2cm 0.2cm 0.2cm 0.2cm,clip=true]{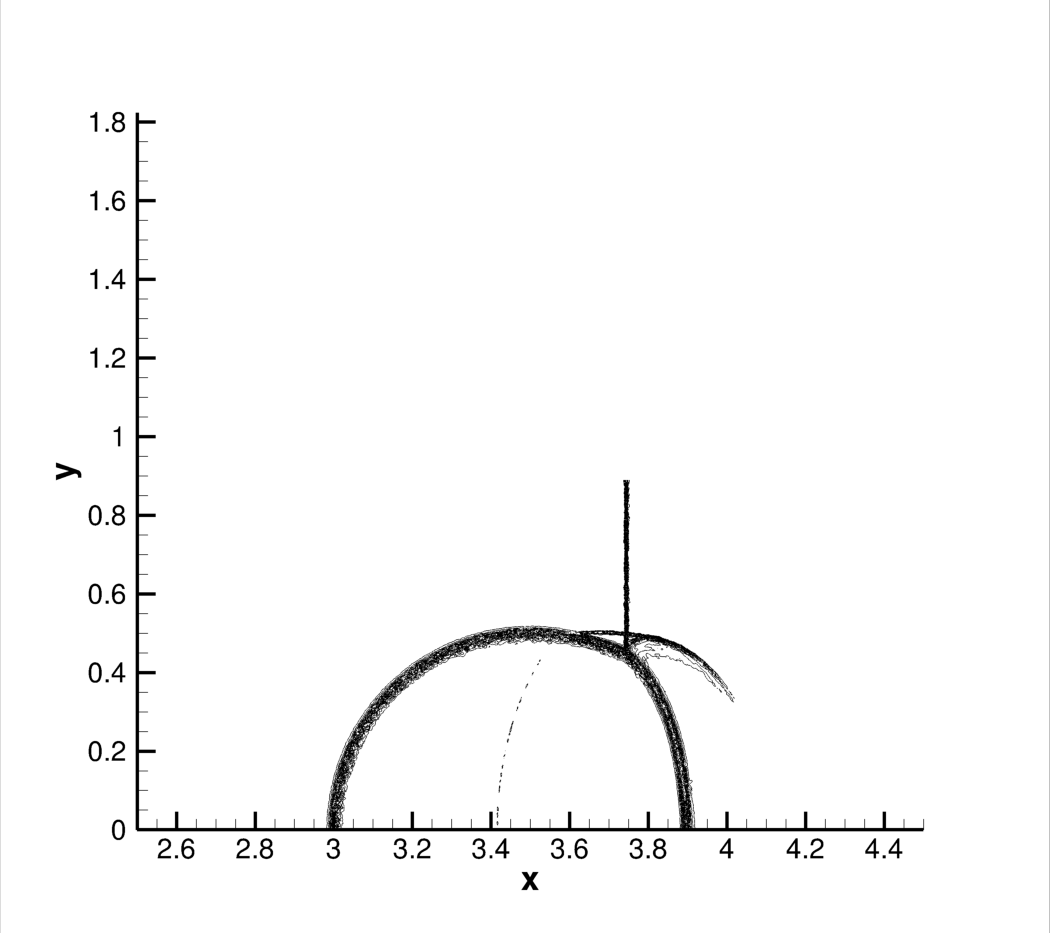}}\\
 \subfloat{\begin{picture}(0,0) \put(-10,50){\rotatebox{90}{$t=102\mu s$}} \end{picture}}
 \subfloat{\includegraphics[height=0.18\paperwidth,trim=0.2cm 0.2cm 0.2cm 0.2cm,clip=true]{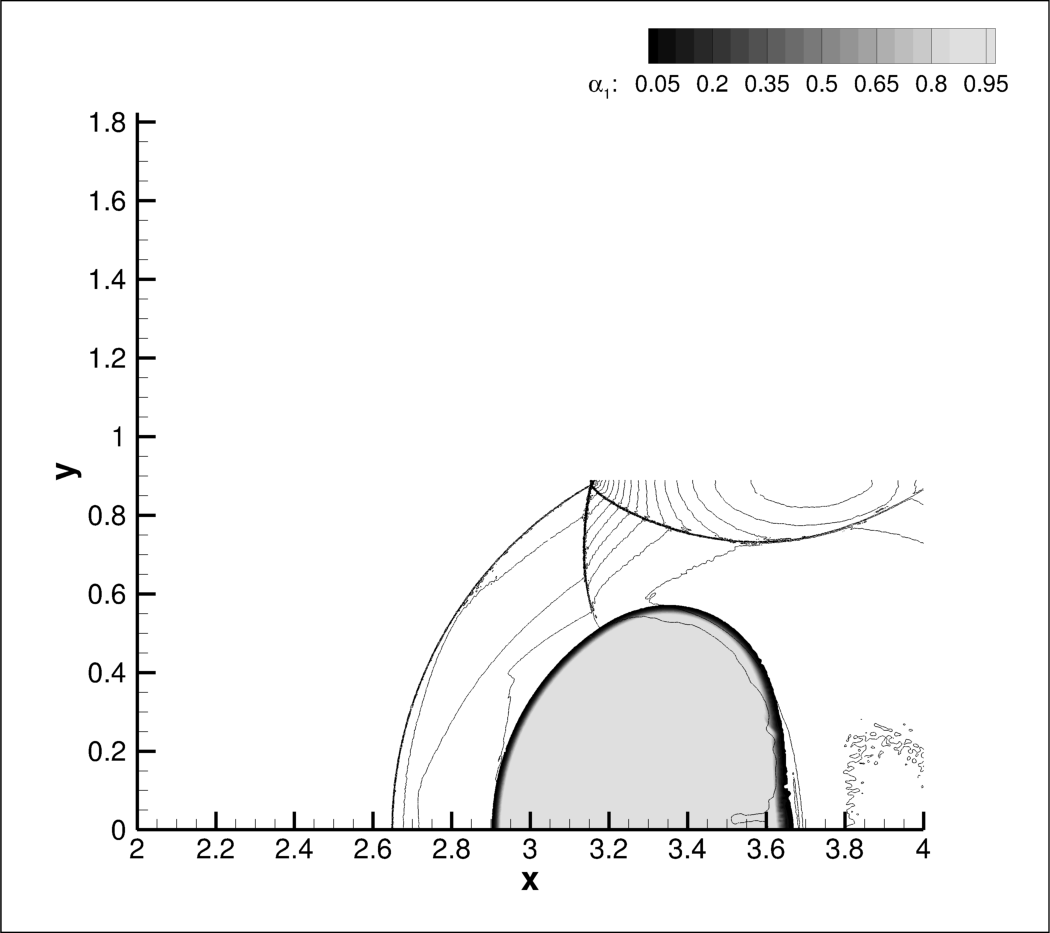}}\hspace*{-0.2cm}
  \subfloat{\includegraphics[height=0.18\paperwidth,trim=0.2cm 0.2cm 0.2cm 0.2cm,clip=true]{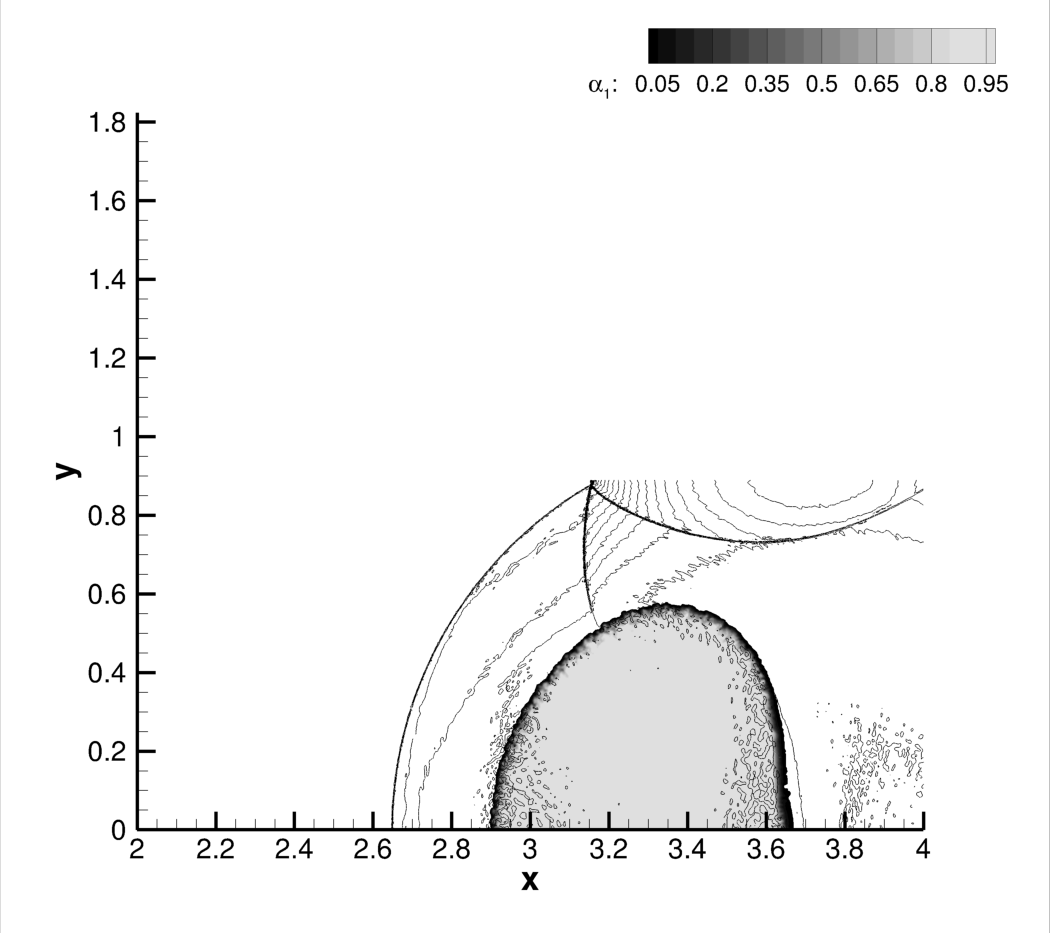}}\hspace*{-0.2cm}
  \subfloat{\includegraphics[height=0.18\paperwidth,trim=0.2cm 0.2cm 0.2cm 0.2cm,clip=true]{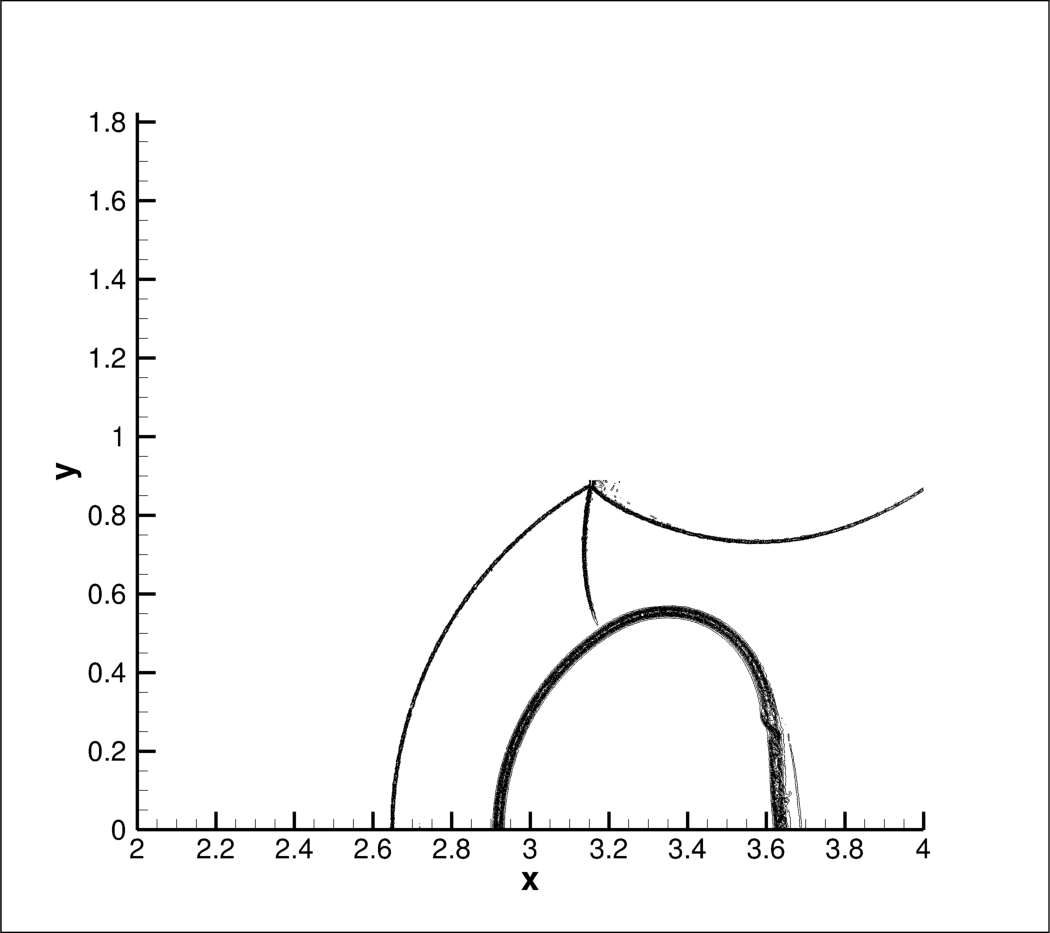}}\hspace*{-0.2cm}
  \subfloat{\includegraphics[height=0.18\paperwidth,trim=0.2cm 0.2cm 0.2cm 0.2cm,clip=true]{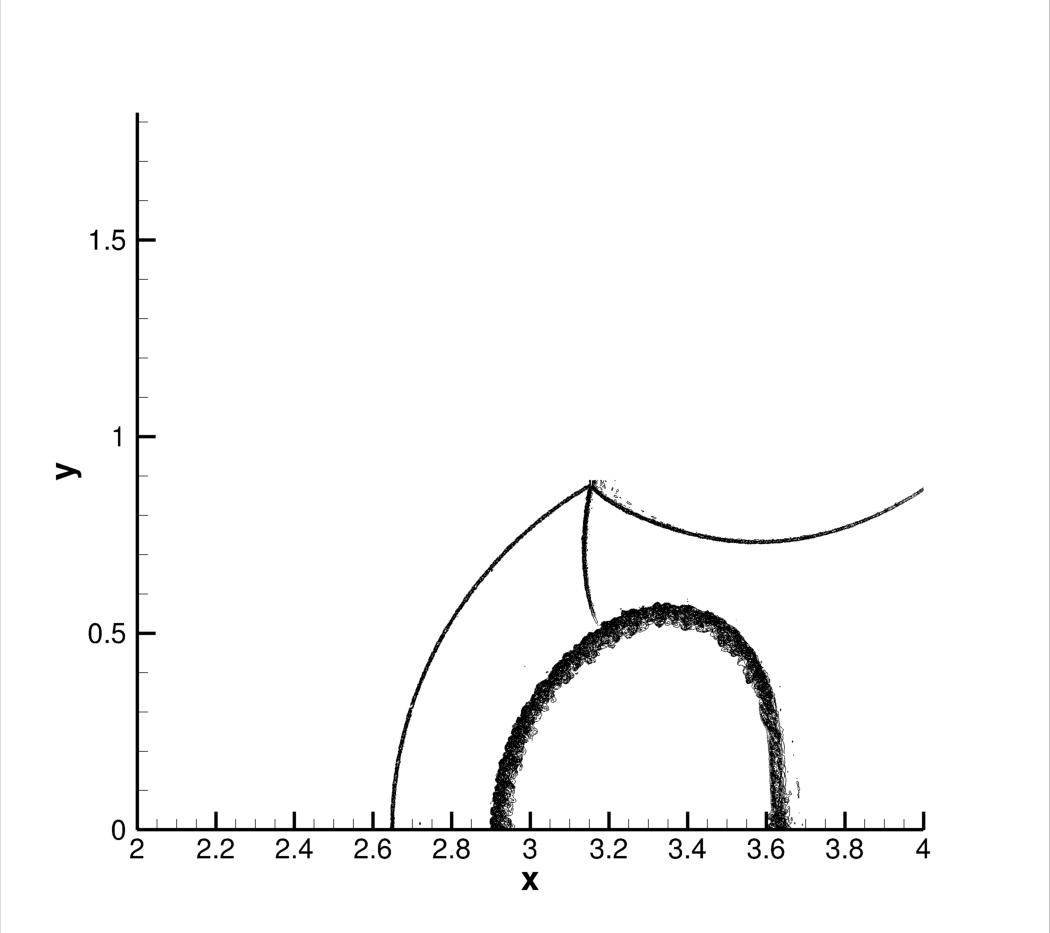}}\\
 \subfloat[]{\begin{picture}(0,0) \put(-10,50){\rotatebox{90}{$t=427\mu s$}} \end{picture}}
 \subfloat{\includegraphics[height=0.18\paperwidth,trim=0.2cm 0.2cm 0.2cm 0.2cm,clip=true]{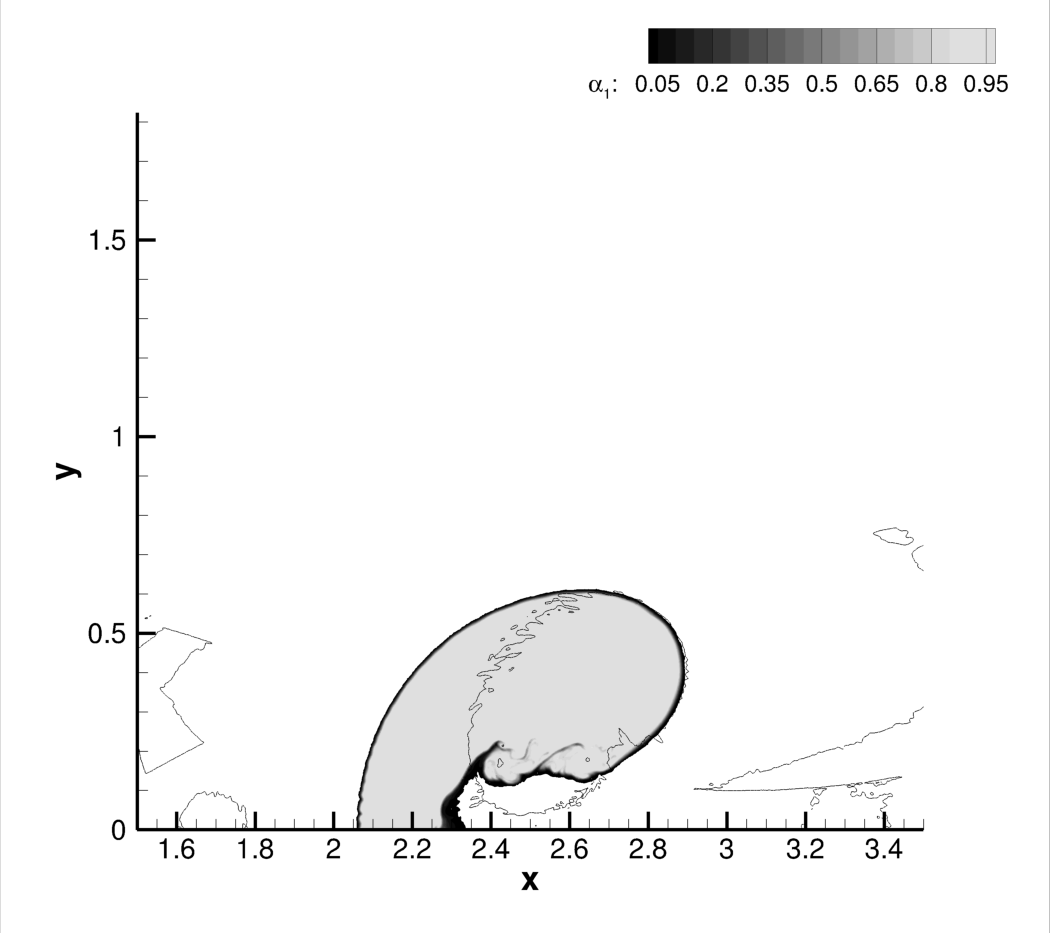}}\hspace*{-0.2cm}
  \subfloat{\includegraphics[height=0.18\paperwidth,trim=0.2cm 0.2cm 0.2cm 0.2cm,clip=true]{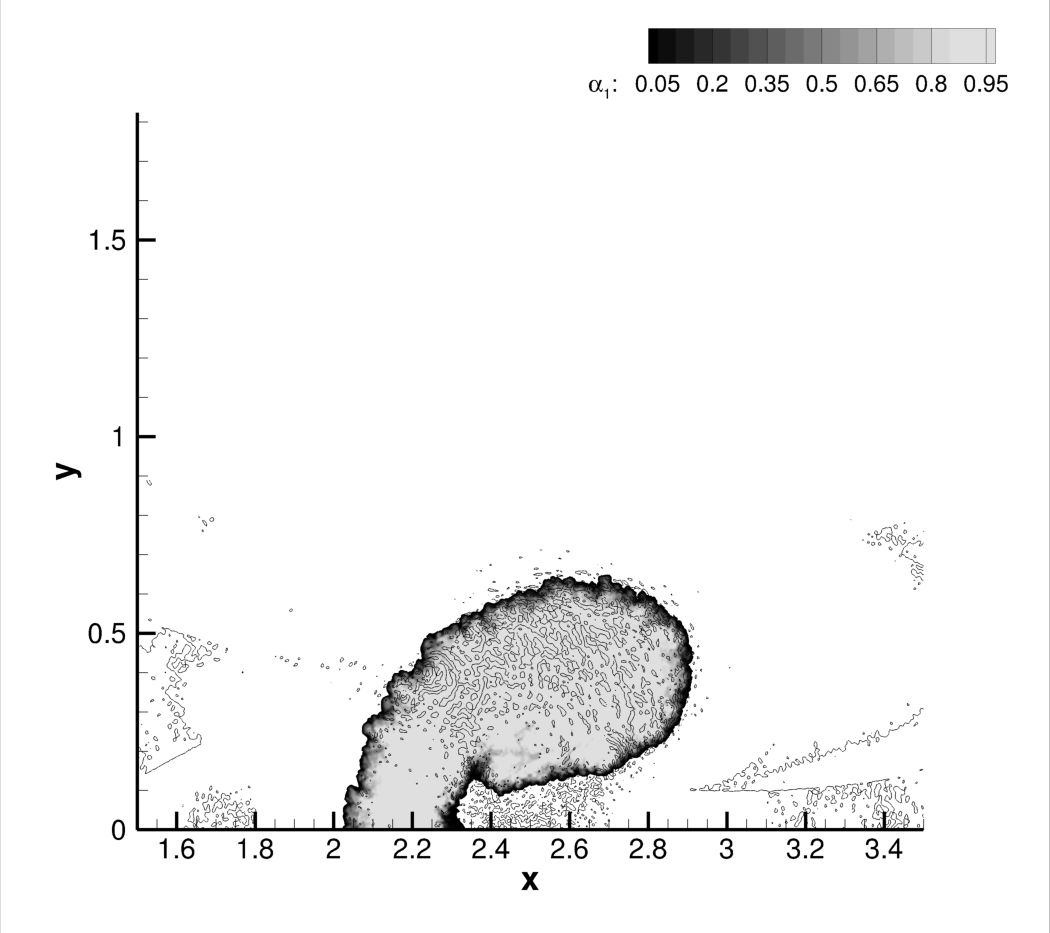}}\hspace*{-0.2cm}
	\subfloat{\includegraphics[height=0.18\paperwidth,trim=0.2cm 0.2cm 0.2cm 0.2cm,clip=true]{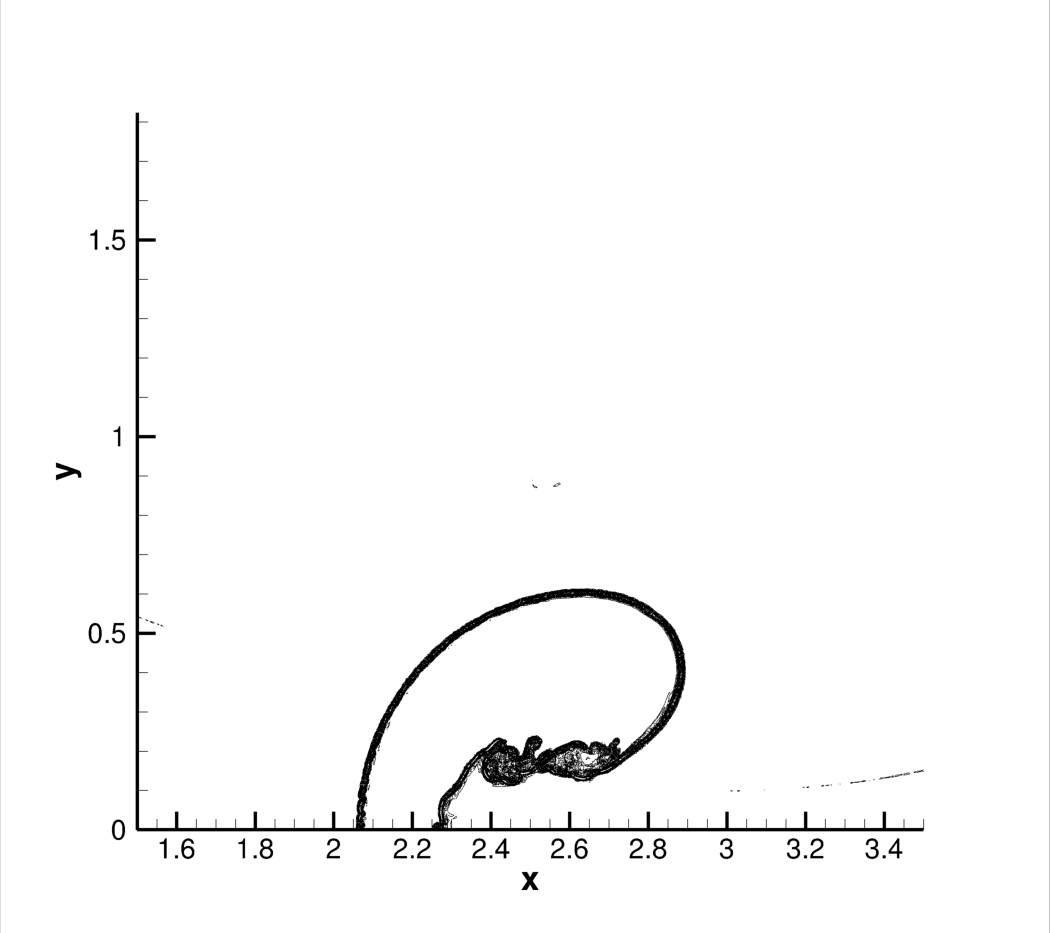}}\hspace*{-0.2cm}
  \subfloat{\includegraphics[height=0.18\paperwidth,trim=0.2cm 0.2cm 0.2cm 0.2cm,clip=true]{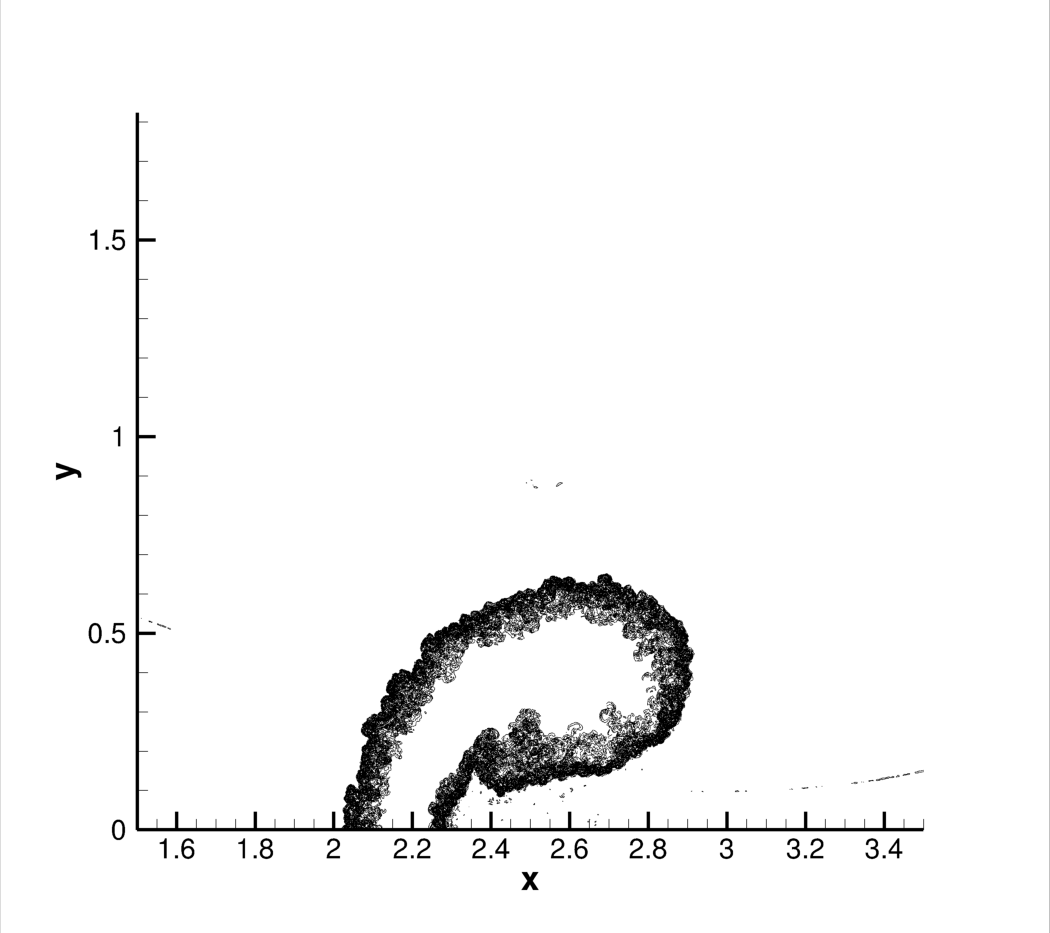}}\\
 \subfloat[]{\begin{picture}(0,0) \put(-10,50){\rotatebox{90}{$t=674\mu s$}} \end{picture}}
 \subfloat[CP fluxes]{\includegraphics[height=0.18\paperwidth,trim=0.2cm 0.2cm 0.2cm 0.2cm,clip=true]{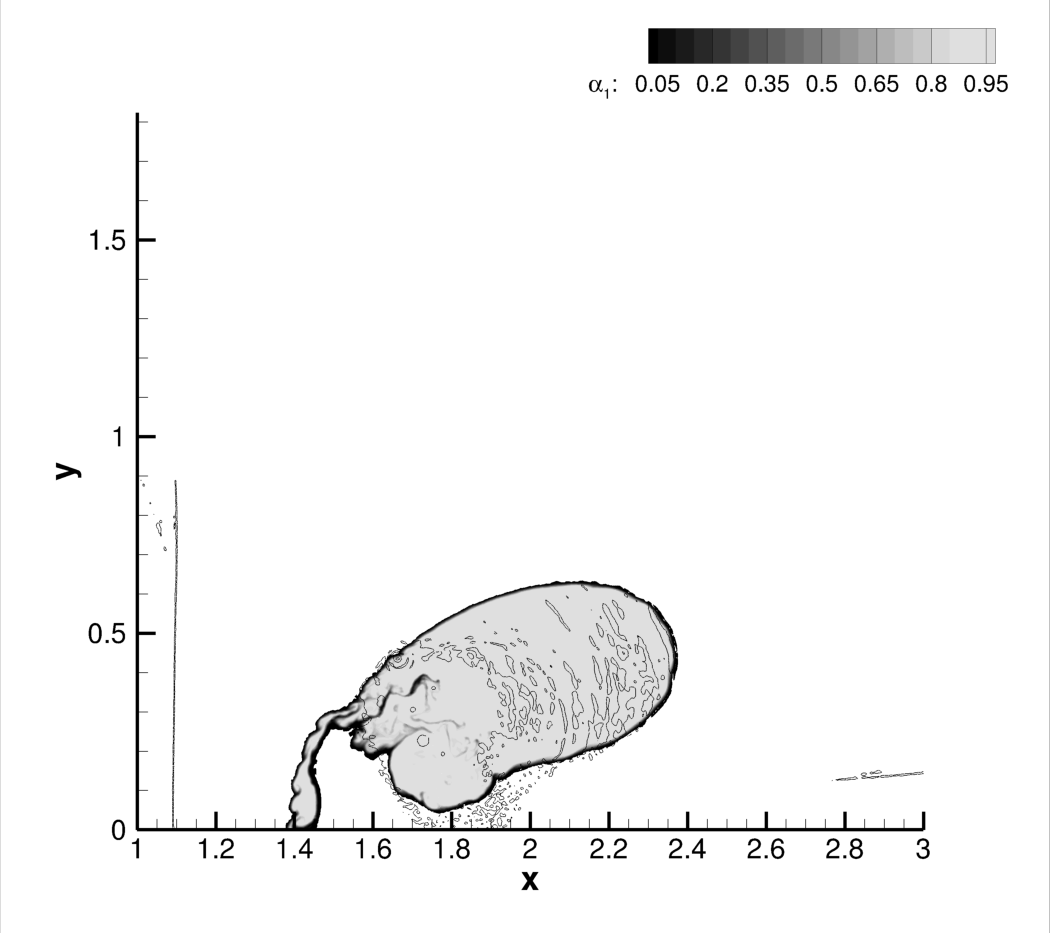}}\hspace*{-0.2cm}
  \subfloat[EC fluxes]{\includegraphics[height=0.18\paperwidth,trim=0.2cm 0.2cm 0.2cm 0.2cm,clip=true]{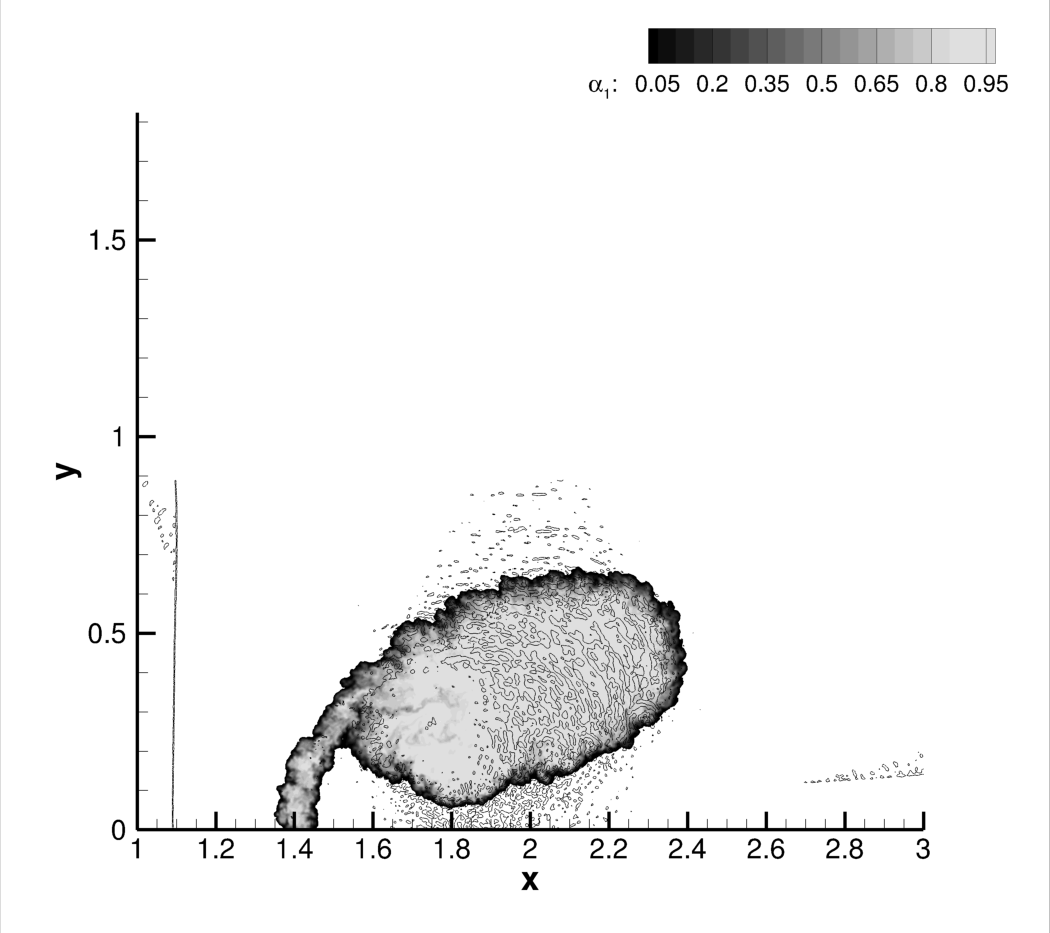}}\hspace*{-0.2cm}
	\subfloat[CP fluxes]{\includegraphics[height=0.18\paperwidth,trim=0.2cm 0.2cm 0.2cm 0.2cm,clip=true]{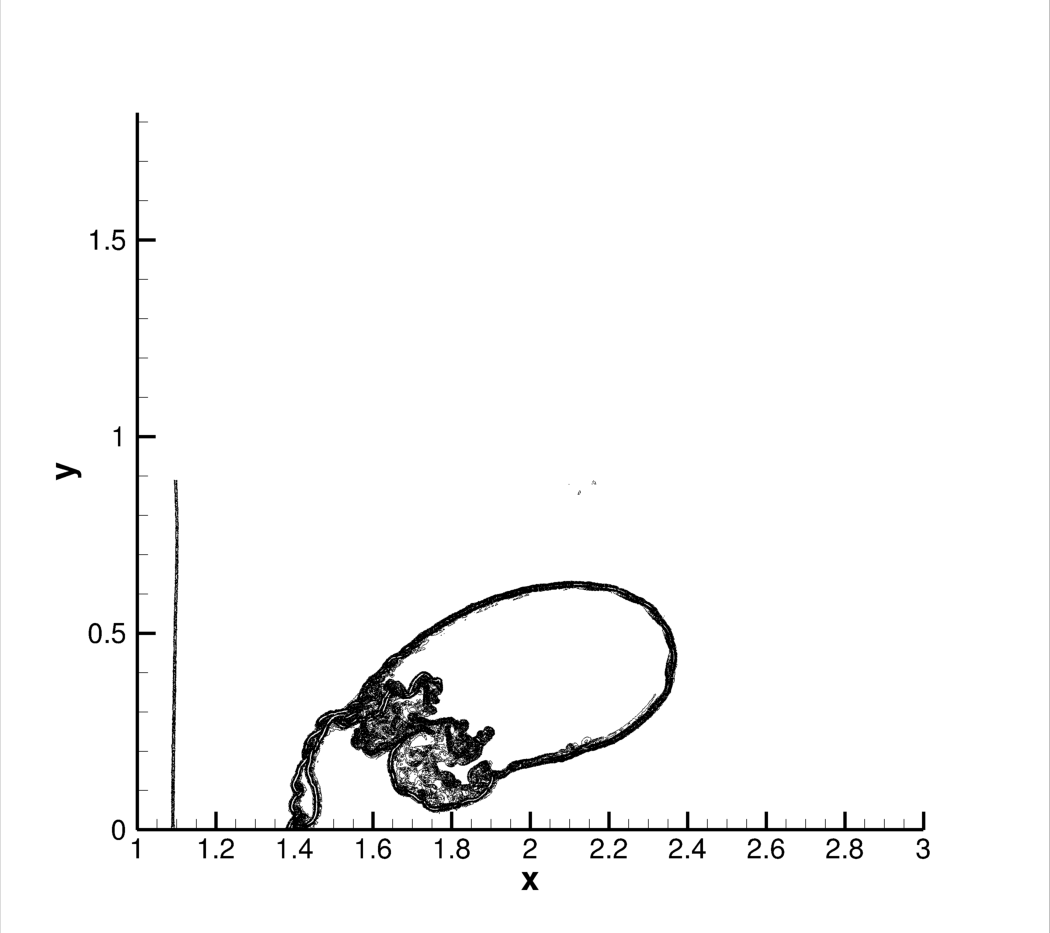}}\hspace*{-0.2cm}
  \subfloat[EC fluxes]{\includegraphics[height=0.18\paperwidth,trim=0.2cm 0.2cm 0.2cm 0.2cm,clip=true]{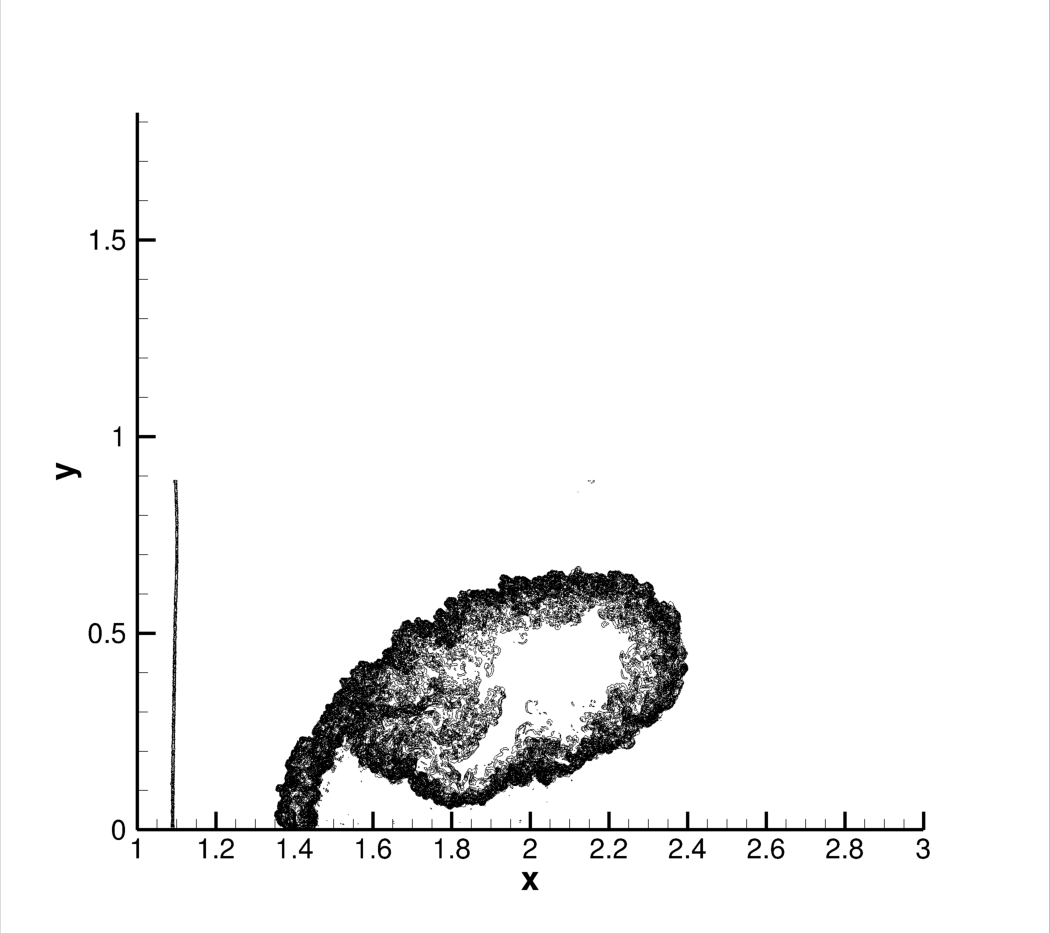}}
  \caption{Interaction of a $M=1.22$ shock in air with a helium bubble: fourth order accurate in space, $p=3$, numerical simulations obtained at different times using either CP, or EC fluxes in the volume integral, on an unstructured mesh with $433016$ elements. For each snapshot the left image shows the helium void fraction contours (color levels) and the pressure contours (lines), while the right image shows the numerical Schlieren of the density $\phi=\exp{|\nabla\rho|/|\nabla\rho_{\max}|}$.}
  \label{result: chap 5 shock-bubble interaction}
\end{figure}
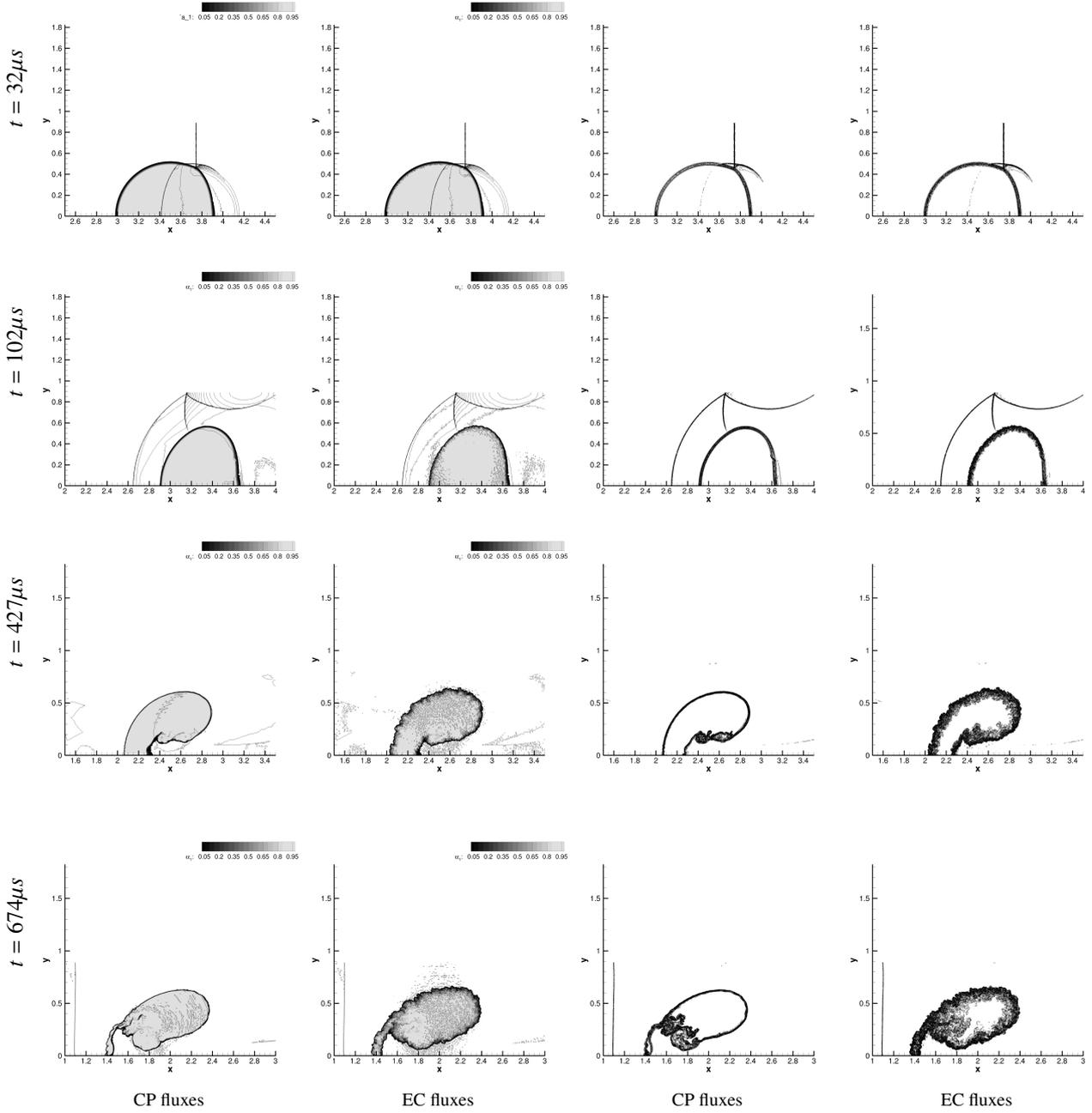

We now consider the shock-bubble interaction. \Cref{result: chap 5 shock-bubble interaction} shows the deformation of the bubble at several physical times as the left traveling shock passes through it and represents contours of void fraction of the helium bubble, $\alpha_1$, and mixture pressure, $\press$, together with numerical Schlieren. We observe that the scheme with CP fluxes allows a better and sharper resolution of the bubble interface for all physical times and is able to accurately capture the shock and bubble dynamics. The bubble interface develops vortices after interacting with shock due to a Kelvin-Helmholtz instability. Once again, results with the EC fluxes show some spurious oscillations at the material interface before and after the interaction in contrast to CP fluxes, while both CP and EC schemes show good resolutions of the shock \red{(shock is captured within, almost, one single cell)}. 


%
\subsection{Strong shock wave-hydrogen bubble interaction}\label{ssec: Strong shock in air interacts with a Hydrogen bubble}

We finally consider an interaction problem of a strong $M=2$ shock in air with a hydrogen bubble that has been numerically investigated in \cite{billet2008impact,sjogreen2003grid}. Compared to \cref{ssec: Shock in air interacts with a helium bubble}, these conditions result in faster shock and bubble dynamics. The computational domain for this test is $\Omega_h=[0,22.5]\times [0,7.5]$ and is discretized with an unstructured mesh with $154622$ elements. The hydrogen bubble ($\gamma_1=1.41$ and $\Cv_1=7.424$) is initially centered at $x=4$ and $y=0$, and the right traveling shock located at $x=7$ in air ($\gamma_2=1.353$ and $\Cv_2=0.523$). The initial data is made nondimensional with the pre-shock density, velocity and temperature of the air and a length scale of $1$mm. We impose symmetry conditions at the top and bottom boundaries, along with supersonic inflow condition at the left boundary and nonreflecting conditions at the right boundary.

\begin{figure}[htbp]
 \center
 \subfloat[]{\begin{picture}(0,0) \put(-10,50){\rotatebox{90}{$t=1.5\mu s$}} \end{picture}}
 \subfloat{\includegraphics[height=0.18\paperwidth,trim=0.2cm 0.2cm 0.2cm 0.2cm,clip=true]{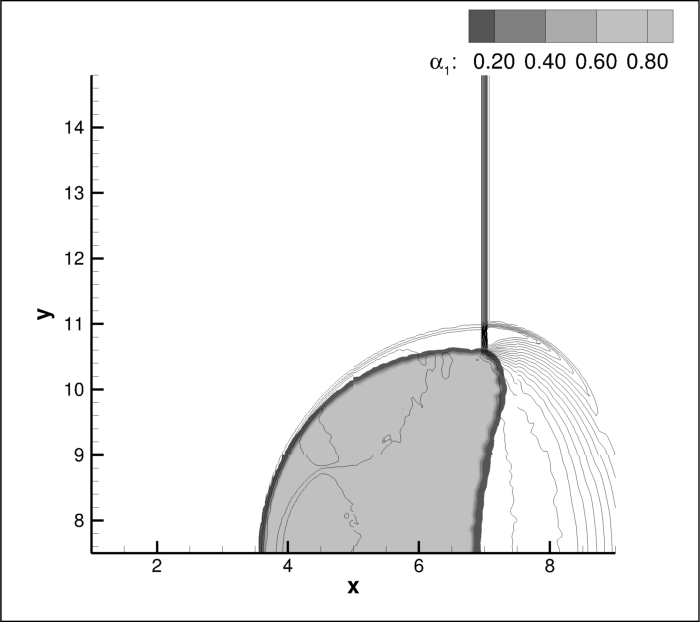}}\hspace*{-0.2cm}
  \subfloat{\includegraphics[height=0.18\paperwidth,trim=0.2cm 0.2cm 0.2cm 0.2cm,clip=true]{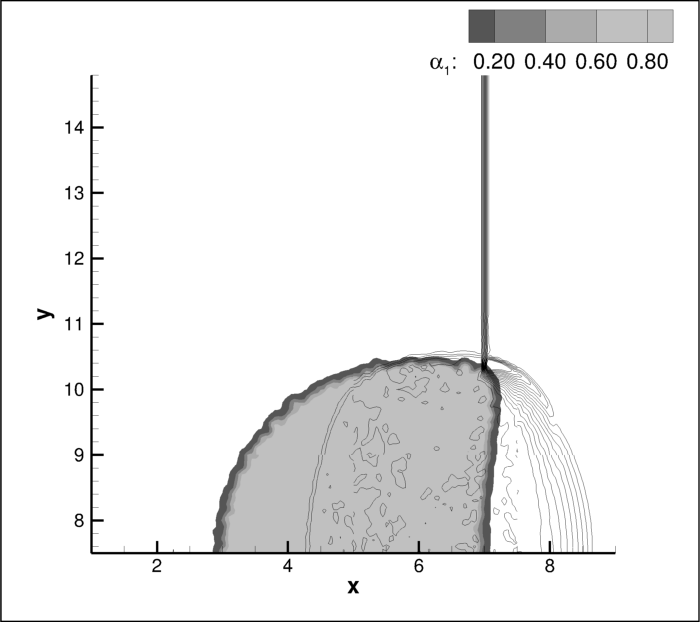}}\hspace*{-0.2cm}
  \subfloat{\includegraphics[height=0.18\paperwidth,trim=0.2cm 0.2cm 0.2cm 0.2cm,clip=true]{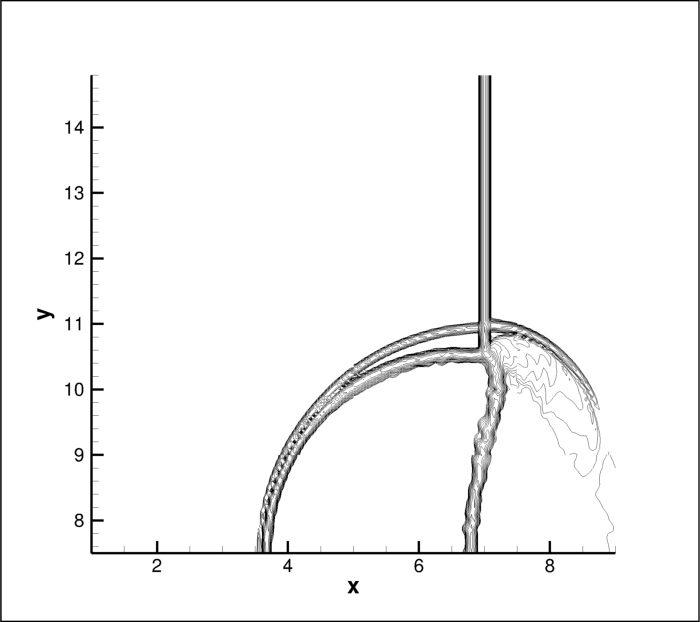}}\hspace*{-0.2cm}
  \subfloat{\includegraphics[height=0.18\paperwidth,trim=0.2cm 0.2cm 0.2cm 0.2cm,clip=true]{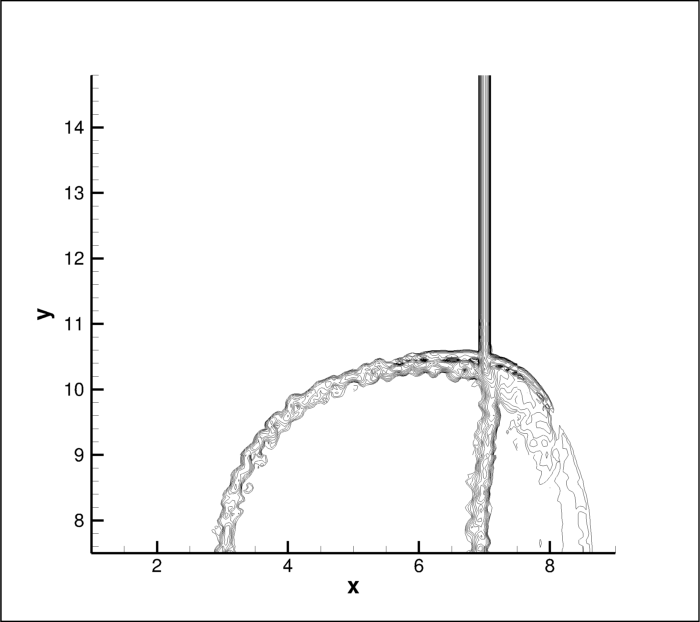}}\\
%
%
 \subfloat[]{\begin{picture}(0,0) \put(-10,50){\rotatebox{90}{$t=2.5\mu s$}} \end{picture}}
 \subfloat{\includegraphics[height=0.18\paperwidth,trim=0.2cm 0.2cm 0.2cm 0.2cm,clip=true]{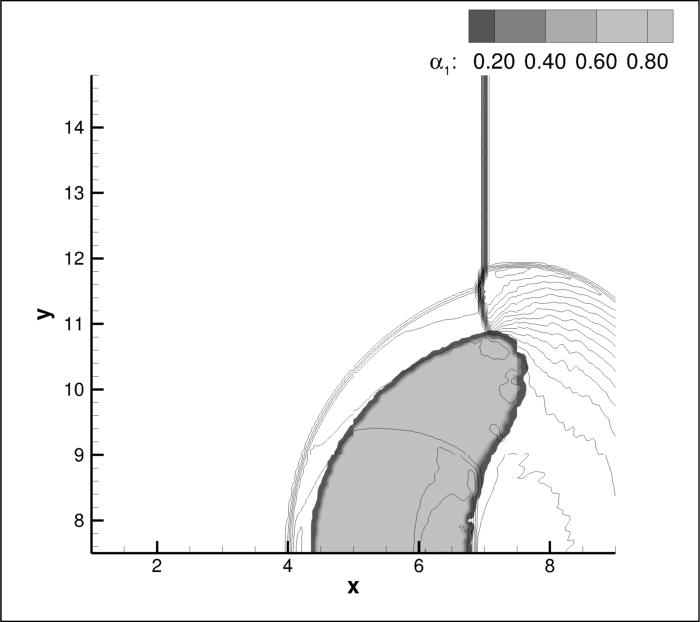}}\hspace*{-0.2cm}
  \subfloat{\includegraphics[height=0.18\paperwidth,trim=0.2cm 0.2cm 0.2cm 0.2cm,clip=true]{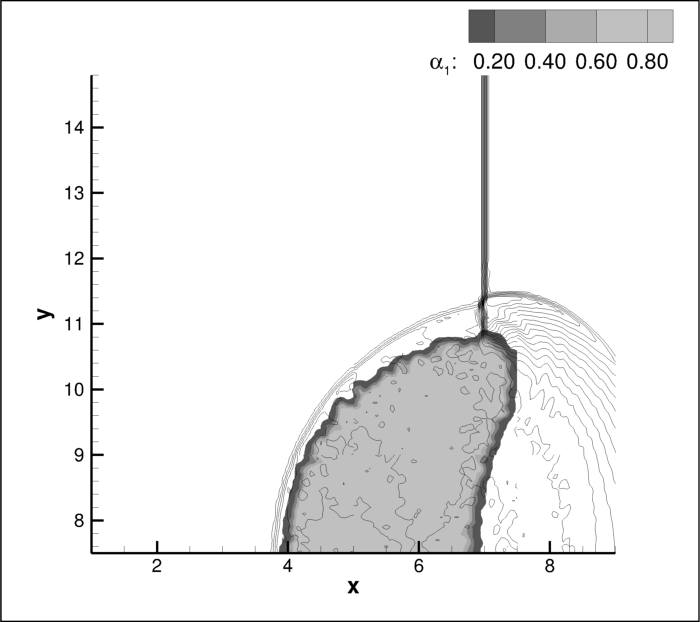}}\hspace*{-0.2cm}
  \subfloat{\includegraphics[height=0.18\paperwidth,trim=0.2cm 0.2cm 0.2cm 0.2cm,clip=true]{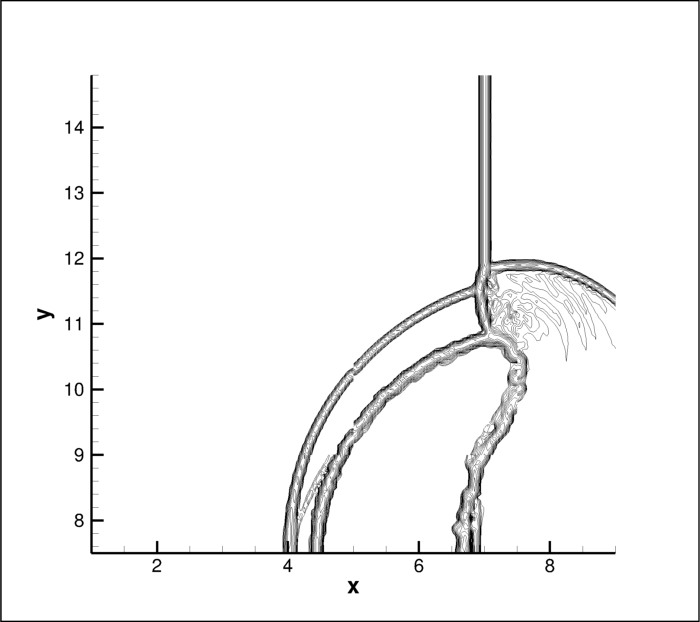}}\hspace*{-0.2cm}
  \subfloat{\includegraphics[height=0.18\paperwidth,trim=0.2cm 0.2cm 0.2cm 0.2cm,clip=true]{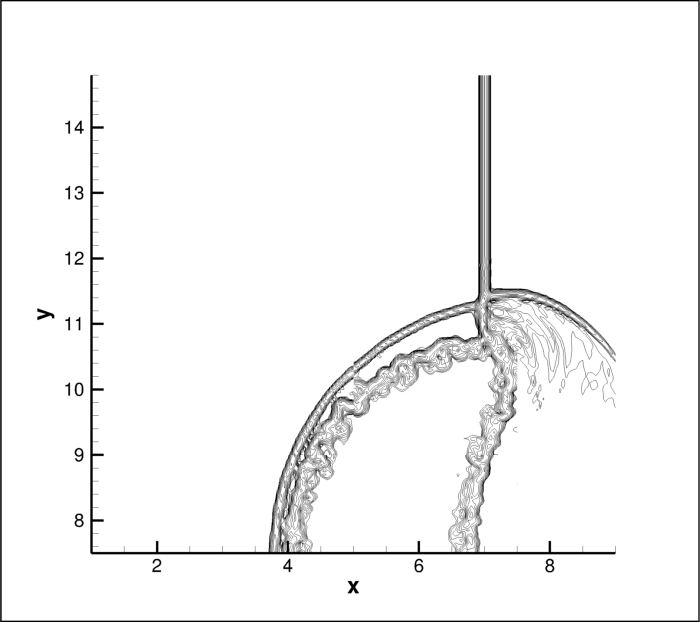}}\\
%
%
 \subfloat[]{\begin{picture}(0,0) \put(-10,50){\rotatebox{90}{$t=4.0\mu s$}} \end{picture}}
 \subfloat{\includegraphics[height=0.18\paperwidth,trim=0.2cm 0.2cm 0.2cm 0.2cm,clip=true]{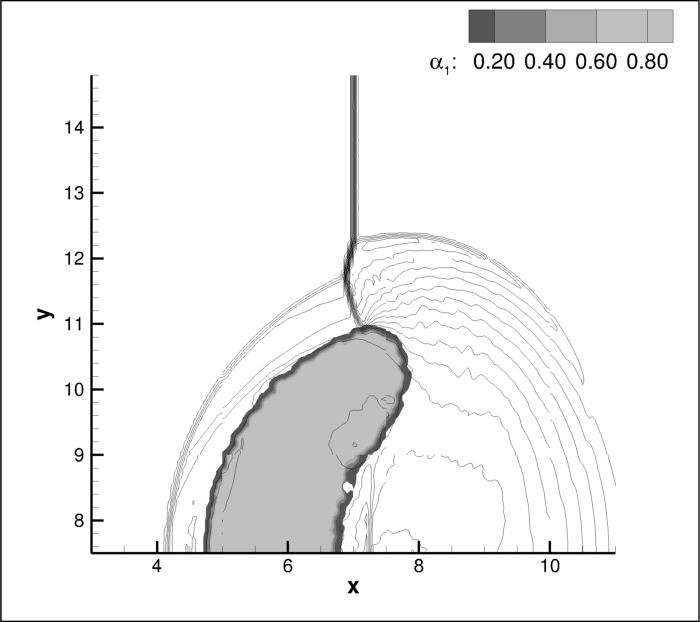}}\hspace*{-0.2cm}
  \subfloat{\includegraphics[height=0.18\paperwidth,trim=0.2cm 0.2cm 0.2cm 0.2cm,clip=true]{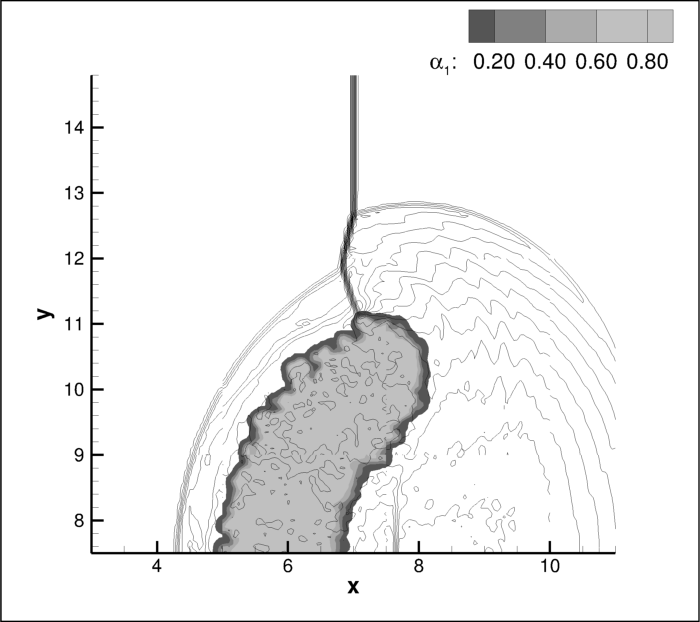}}\hspace*{-0.2cm}
  \subfloat{\includegraphics[height=0.18\paperwidth,trim=0.2cm 0.2cm 0.2cm 0.2cm,clip=true]{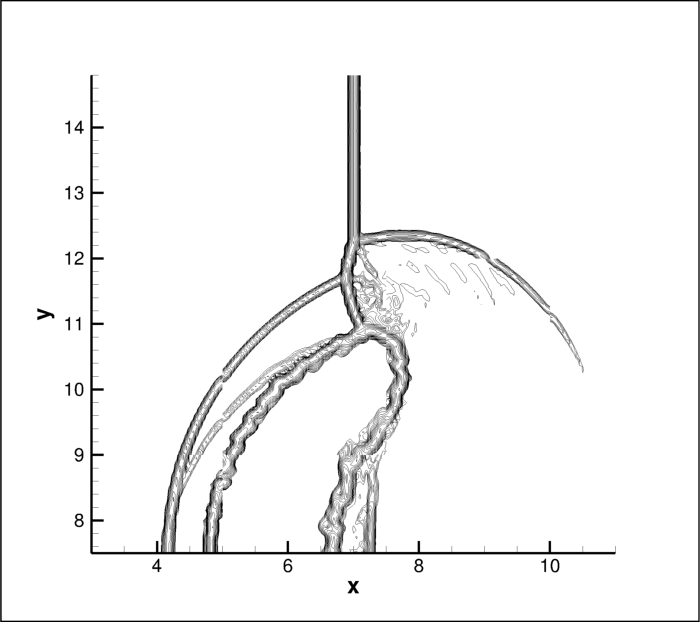}}\hspace*{-0.2cm}
  \subfloat{\includegraphics[height=0.18\paperwidth,trim=0.2cm 0.2cm 0.2cm 0.2cm,clip=true]{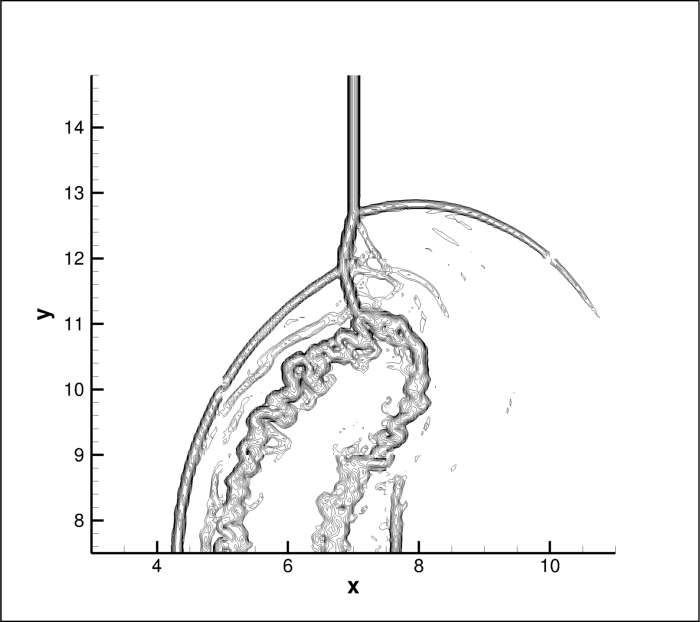}}\\  
 \subfloat[]{\begin{picture}(0,0) \put(-10,50){\rotatebox{90}{$t=8.8\mu s$}} \end{picture}}
 \subfloat{\includegraphics[height=0.18\paperwidth,trim=0.2cm 0.2cm 0.2cm 0.2cm,clip=true]{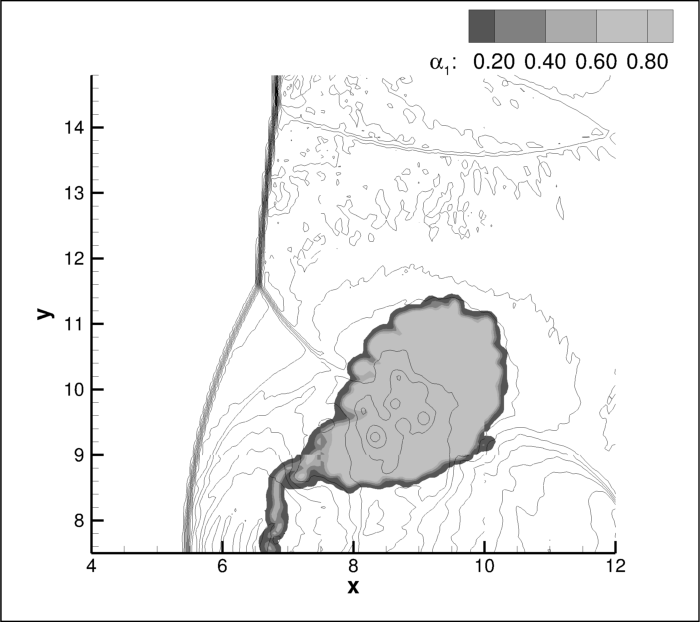}}\hspace*{-0.2cm}
  \subfloat{\includegraphics[height=0.18\paperwidth,trim=0.2cm 0.2cm 0.2cm 0.2cm,clip=true]{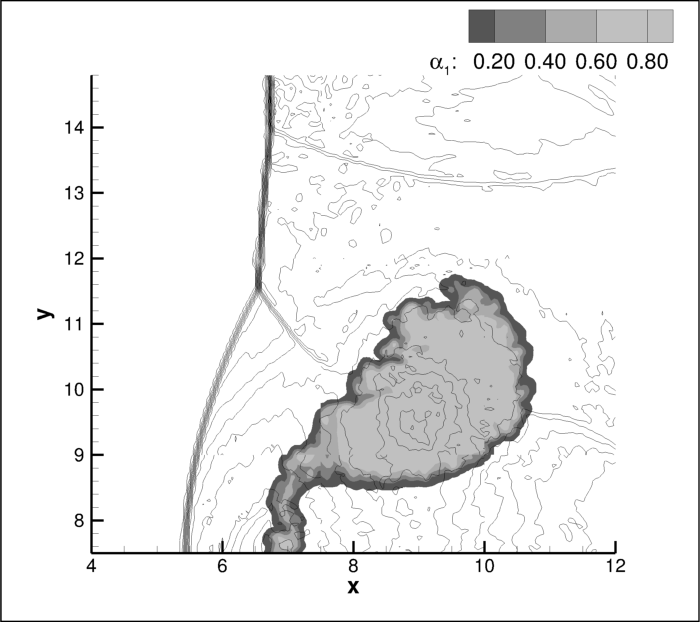}}\hspace*{-0.2cm}
  \subfloat{\includegraphics[height=0.18\paperwidth,trim=0.2cm 0.2cm 0.2cm 0.2cm,clip=true]{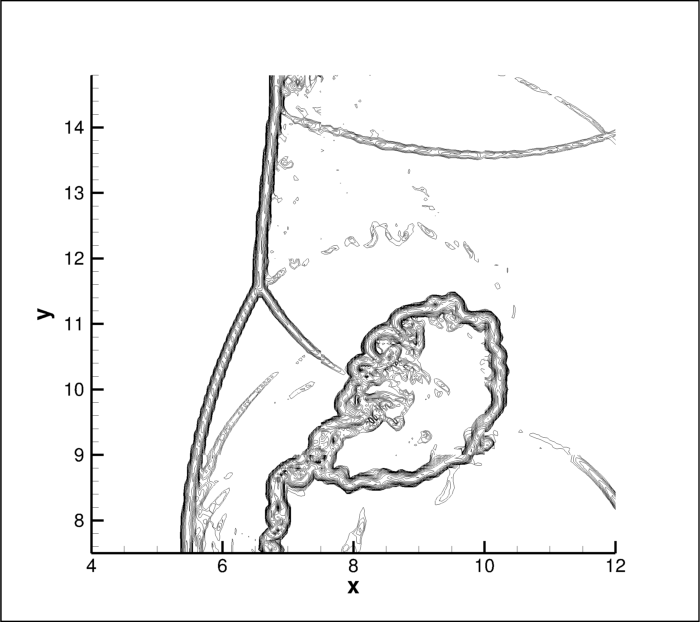}}\hspace*{-0.2cm}
  \subfloat{\includegraphics[height=0.18\paperwidth,trim=0.2cm 0.2cm 0.2cm 0.2cm,clip=true]{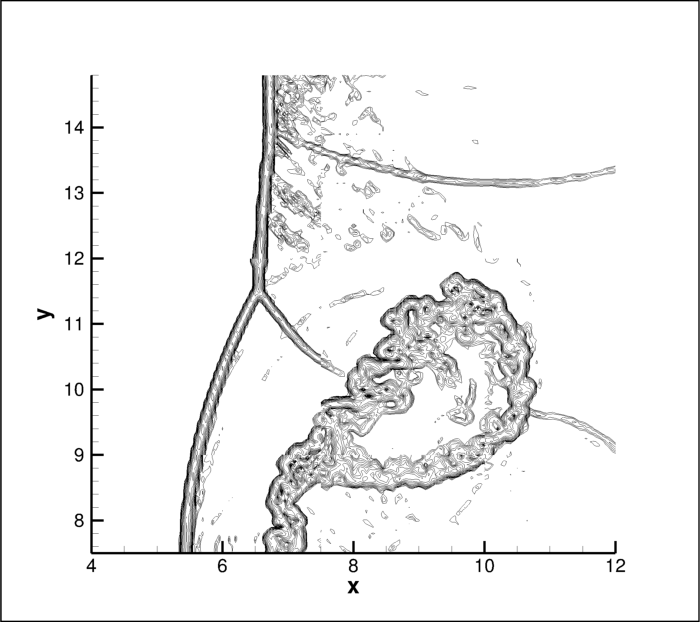}}\\
 \subfloat[]{\begin{picture}(0,0) \put(-10,50){\rotatebox{90}{$t=13.6\mu s$}} \end{picture}}
 \subfloat[CP fluxes]{\includegraphics[height=0.18\paperwidth,trim=0.2cm 0.2cm 0.2cm 0.2cm,clip=true]{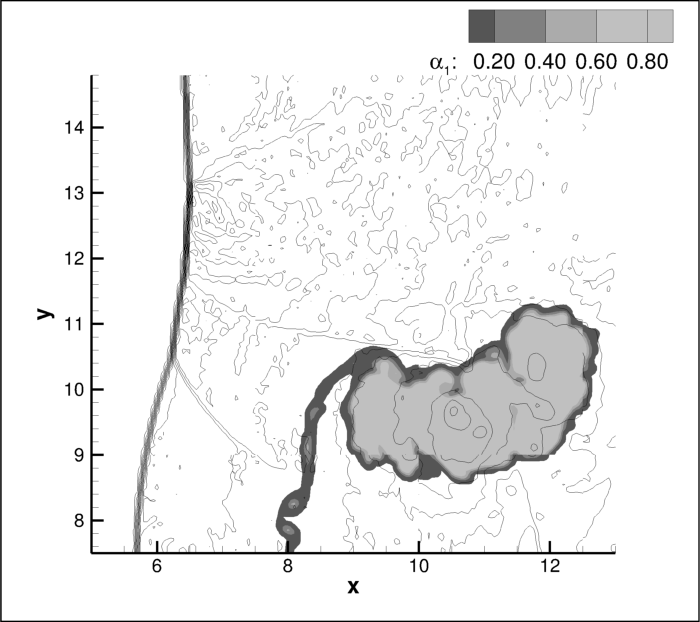}}\hspace*{-0.2cm}
  \subfloat[EC fluxes]{\includegraphics[height=0.18\paperwidth,trim=0.2cm 0.2cm 0.2cm 0.2cm,clip=true]{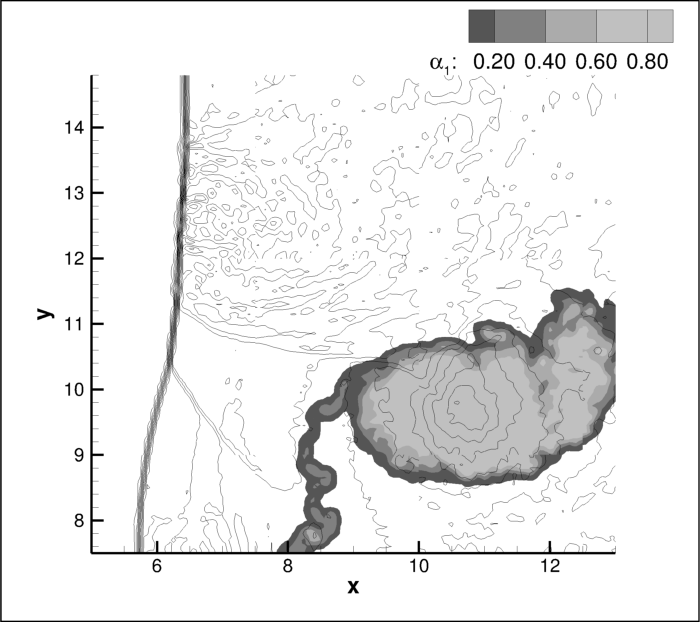}}\hspace*{-0.2cm}
  \subfloat[CP fluxes]{\includegraphics[height=0.18\paperwidth,trim=0.2cm 0.2cm 0.2cm 0.2cm,clip=true]{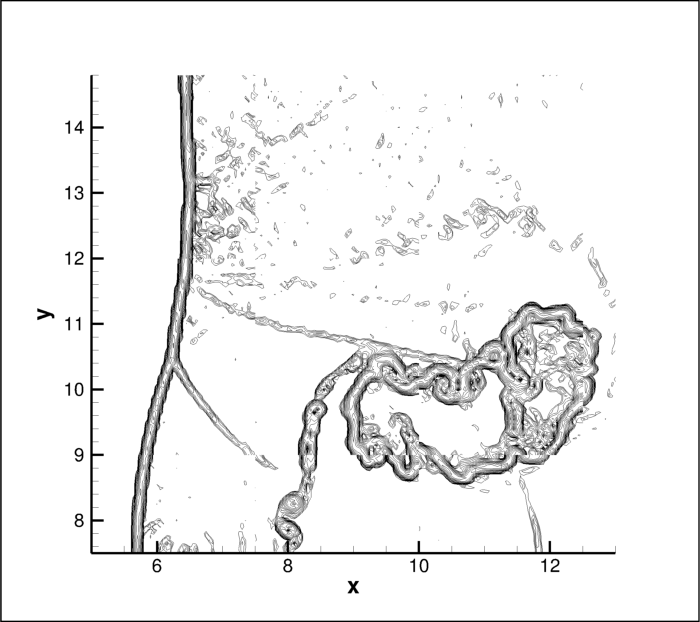}}\hspace*{-0.2cm}
  \subfloat[EC fluxes]{\includegraphics[height=0.18\paperwidth,trim=0.2cm 0.2cm 0.2cm 0.2cm,clip=true]{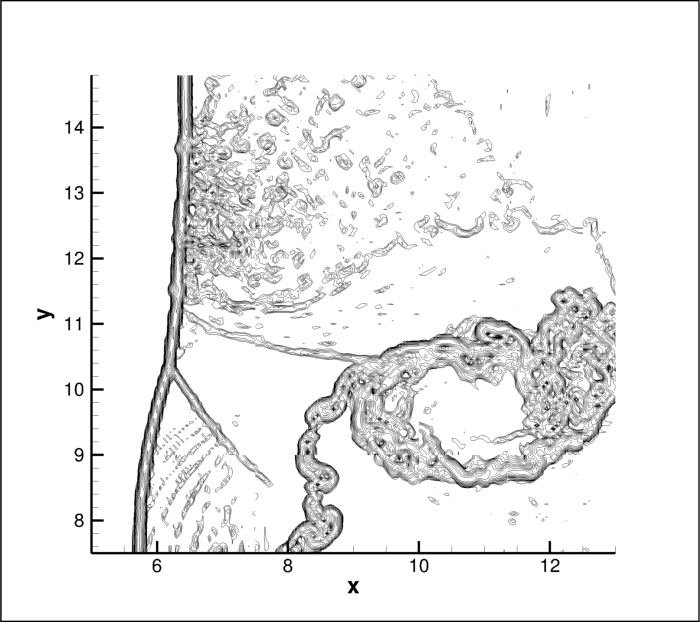}}
	\caption{Interaction of a $M=2$ shock in air with a Hydrogen bubble: fourth-order accurate in space, $p=3$, numerical solutions obtained at different times using either CP, or EC numerical fluxes in the volume integral, on an unstructured mesh with $154622$ elements. For each snapshot the left image shows the helium void fraction contours (color levels) and the pressure contours (lines), while the right image shows the numerical Schlieren of the density $\phi=\exp{|\nabla\rho|/|\nabla\rho_{\max}|}$.}
  \label{result: chap 5 shock-bubble interaction_H2}    
\end{figure}

\Cref{result: chap 5 shock-bubble interaction_H2} shows the deformation of the bubble as the shock passes through it, where we plot contours of the void fraction of the hydrogen bubble, $\alpha_1$, mixture pressure pressure, $\press$, and the numerical Schlieren. Here, we once again observe that the numerical scheme is able to resolve the bubble interface well along with the shock. The oscillations at the interface are due to the Kelvin-Helmholtz instability and they were also observed in \cite{billet2008impact}. Using CP fluxes maintains sharp resolution of the interface while also proving to well capture shocks.

%
%
\section{Concluding remarks}\label{sec: chap 5 summary}

In this work, we propose a high-order, robust and entropy stable discretization of the nonconservative multicomponent SG-gamma model \cite{shyue1998efficient}. The space discretization of this system relies on the DGSEM framework \cite{renac2019entropy,rai2021entropy} based on the modification of the integral over discretization elements where we replace the physical fluxes and nonconservative products by two-point numerical fluctuations. We first extend this framework to multidimensional unstructured grids with curved elements and derive conditions under which the semi-discrete scheme is high-order accurate, free-stream preserving and entropy stable. We then design a robust, entropy stable, material interface preserving, and maximum principle preserving HLLC solver for the SG-gamma model that does not require a root-finding algorithm to evaluate the nonconservative product. The HLLC solver is then used as interface fluctuations in the \red{proposed DG scheme}, while we consider either EC, or CP fluctuations in the integrals over discretization elements. The \red{DGSEM} is shown to either satisfy a semi-discrete entropy inequality with EC fluctuations (when excluding material interfaces), or to preserve material interfaces with CP fluctuations.

We then analyze the fully discrete \red{numerical scheme} with a forward Euler time discretization. We derive conditions on the time step to guarantee that the cell-averaged solution remains in the set of states and satisfies a minimum principle on the entropy and maximum principles on the EOS parameters with either EC, or CP fluctuations. We use a posteriori scaling limiters \cite{zhang2010positivity,zhang2010maximum} to extend these properties to all DOFs within elements, while a strong-stability preserving Runge-Kutta scheme \cite{shu1988efficient} is used for the high-order time integration and keeps properties of the first-order in time scheme.

High-order accurate numerical simulations of flows in one and two space dimensions with discontinuous solutions and complex wave interactions confirm robustness, stability and accuracy of the present scheme with either EC, or CP fluctuations, and the scheme with CP fluxes present better resolution capabilities.

\section*{Declaration of competing interest}
The authors declare that they have no known competing financial interests or personal relationships that could have appeared to influence the work reported in this paper.

\section*{Acknowledgements}
This work was performed as part of P. Rai's PhD thesis. The authors would like to express their gratitude towards ONERA for providing the funds and resources for this research work, and would like to thank the doctoral school at the Ecole Polytechnique for their administrative support.

\bibliographystyle{siam}
\bibliography{main}

\newcommand{\noopsort}[1]{} \newcommand{\printfirst}[2]{#1}
  \newcommand{\singleletter}[1]{#1} \newcommand{\switchargs}[2]{#2#1}
\begin{thebibliography}{10}

\bibitem{abgrall1996prevent}
{\sc R.~Abgrall}, {\em How to prevent pressure oscillations in multicomponent
  flow calculations: a quasi conservative approach}, J. Comput. Phys., 125
  (1996), pp.~150--160.

\bibitem{batten1997choice}
{\sc P.~Batten, N.~Clarke, C.~Lambert, and D.~M. Causon}, {\em On the choice of
  wavespeeds for the {HLLC} {R}iemann solver}, J. Sci. Comput., 18 (1997),
  pp.~1553--1570.

\bibitem{berthon_etal_ES_HLLC_12}
{\sc C.~Berthon, B.~Dubroca, and A.~Sangam}, {\em A local entropy minimum
  principle for deriving entropy preserving schemes}, SIAM J. Numer. Anal., 50
  (2012), pp.~468--491.

\bibitem{billet2008impact}
{\sc G.~Billet, V.~Giovangigli, and G.~De~Gassowski}, {\em Impact of volume
  viscosity on a shock--hydrogen-bubble interaction}, Combust. Theory Model.,
  12 (2008), pp.~221--248.

\bibitem{black1999conservative}
{\sc K.~Black}, {\em A conservative spectral element method for the
  approximation of compressible fluid flow}, Kybernetika, 35 (1999),
  pp.~133--146.

\bibitem{bouchut2004nonlinear}
{\sc F.~Bouchut}, {\em Nonlinear stability of finite Volume Methods for
  hyperbolic conservation laws: And Well-Balanced schemes for sources},
  Springer Science \& Business Media, 2004.

\bibitem{carlier_renac_IDP_22}
{\sc V.~Carlier and F.~Renac}, {\em Invariant domain preserving high-order
  spectral discontinuous approximations of hyperbolic systems},
  arXiv:2203.05452 [math.NA],  (2022).

\bibitem{carpenter2014entropy}
{\sc M.~H. Carpenter, T.~C. Fisher, E.~J. Nielsen, and S.~H. Frankel}, {\em
  Entropy stable spectral collocation schemes for the {N}avier--{S}tokes
  equations: Discontinuous interfaces}, SIAM J. Sci. Comput., 36 (2014),
  pp.~B835--B867.

\bibitem{castro2017well}
{\sc M.~J. Castro, T.~M. de~Luna, and C.~Par{\'e}s}, {\em Well-balanced schemes
  and path-conservative numerical methods}, in Handbook of Numer. Anal.,
  vol.~18, Elsevier, 2017, pp.~131--175.

\bibitem{castro2013entropy}
{\sc M.~J. Castro, U.~S. Fjordholm, S.~Mishra, and C.~Par{\'e}s}, {\em Entropy
  conservative and entropy stable schemes for nonconservative hyperbolic
  systems}, SIAM J. Numer. Anal., 51 (2013), pp.~1371--1391.

\bibitem{castro2006high}
{\sc M.~J. Castro, J.~Gallardo, and C.~Par{\'e}s}, {\em High order finite
  volume schemes based on reconstruction of states for solving hyperbolic
  systems with nonconservative products. applications to shallow-water
  systems}, Math. {C}omput., 75 (2006), pp.~1103--1134.

\bibitem{castro_etal_HLLC_NC_13}
{\sc {Castro D\'{\i}az, Manuel Jes\'us}, {Fern\'andez-Nieto, Enrique Domingo},
  {Morales de Luna, Tom\'as}, {Narbona-Reina, Gladys}, and {Par\'es, Carlos}},
  {\em A hllc scheme for nonconservative hyperbolic problems. application to
  turbidity currents with sediment transport}, ESAIM: M2AN, 47 (2013),
  pp.~1--32.

\bibitem{chen2017entropy}
{\sc T.~Chen and C.-W. Shu}, {\em Entropy stable high order discontinuous
  {G}alerkin methods with suitable quadrature rules for hyperbolic conservation
  laws}, J. Comput. Phys., 345 (2017), pp.~427--461.

\bibitem{cheng2020discontinuous}
{\sc J.~Cheng, F.~Zhang, and T.~Liu}, {\em A discontinuous {G}alerkin method
  for the simulation of compressible gas-gas and gas-water two-medium flows},
  J. Computat. Phys., 403 (2020), p.~109059.

\bibitem{rai2021entropy}
{\sc F.~Coquel, C.~Marmignon, P.~Rai, and F.~Renac}, {\em An entropy stable
  high-order discontinuous {G}alerkin spectral element method for the
  {B}aer-{N}unziato two-phase flow model}, J. Comput. Phys., 431 (2021),
  p.~110135.

\bibitem{coralic2014finite}
{\sc V.~Coralic and T.~Colonius}, {\em Finite-volume {WENO} scheme for viscous
  compressible multicomponent flows}, J. Computat. Phys., 274 (2014),
  pp.~95--121.

\bibitem{de2015new}
{\sc M.~T.~H. de~Frahan, S.~Varadan, and E.~Johnsen}, {\em A new limiting
  procedure for discontinuous {G}alerkin methods applied to compressible
  multiphase flows with shocks and interfaces}, J. Comput. Phys., 280 (2015),
  pp.~489--509.

\bibitem{derigs2017novel}
{\sc D.~Derigs, A.~R. Winters, G.~J. Gassner, and S.~Walch}, {\em A novel
  averaging technique for discrete entropy-stable dissipation operators for
  ideal {MHD}}, J. Comput. Phys., 330 (2017), pp.~624--632.

\bibitem{despres98}
{\sc B.~Despres}, {\em Entropy inequality for high order discontinuous
  {G}alerkin approximation of {E}uler equations}, in Hyperbolic Problems:
  Theory, Numerics, Applications, M.~Fey and R.~Jeltsch, eds., Basel, 1999,
  Birkh{\"a}user Basel, pp.~225--231.

\bibitem{Despres2000}
\leavevmode\vrule height 2pt depth -1.6pt width 23pt, {\em Discontinuous
  galerkin method for the numerical solution of euler equations in axisymmetric
  geometry}, in Discontinuous Galerkin Methods, B.~Cockburn, G.~E. Karniadakis,
  and C.-W. Shu, eds., Berlin, Heidelberg, 2000, Springer Berlin Heidelberg,
  pp.~315--320.

\bibitem{dumbser_balsara_HLLEM_16}
{\sc M.~Dumbser and D.~S. Balsara}, {\em A new efficient formulation of the
  hllem riemann solver for general conservative and non-conservative hyperbolic
  systems}, J. Comput. Phys., 304 (2016), pp.~275--319.

\bibitem{einfeldt1991godunov}
{\sc B.~Einfeldt, C.-D. Munz, P.~L. Roe, and B.~Sj{\"o}green}, {\em On
  godunov-type methods near low densities}, J. Comput. Phys., 92 (1991),
  pp.~273--295.

\bibitem{fisher2013high}
{\sc T.~C. Fisher and M.~H. Carpenter}, {\em High-order entropy stable finite
  difference schemes for nonlinear conservation laws: Finite domains}, J.
  Comput. Phys., 252 (2013), pp.~518--557.

\bibitem{fjordholm2012arbitrarily}
{\sc U.~S. Fjordholm, S.~Mishra, and E.~Tadmor}, {\em Arbitrarily high-order
  accurate entropy stable essentially nonoscillatory schemes for systems of
  conservation laws}, SIAM J. Numer. Anal., 50 (2012), pp.~544--573.

\bibitem{franquet2012runge}
{\sc E.~Franquet and V.~Perrier}, {\em {R}unge--{K}utta discontinuous
  {G}alerkin method for the approximation of {B}aer and {N}unziato type
  multiphase models}, J. Comput. Phys., 231 (2012), pp.~4096--4141.

\bibitem{gallice_ARS_NC_03}
{\sc G.~Gallice}, {\em Positive and entropy stable {G}odunov-type schemes for
  gas dynamics and {MHD} equations in {L}agrangian or {E}ulerian coordinates},
  Numer. Math., 94 (2003), pp.~673--713.

\bibitem{gassner2013SBP-SAT}
{\sc G.~J. Gassner}, {\em A skew-symmetric discontinuous {G}alerkin spectral
  element discretization and its relation to {SBP}-{SAT} finite difference
  methods}, SIAM J. Sci. Comput., 35 (2013), pp.~A1233--A1253.

\bibitem{gassner2016split}
{\sc G.~J. Gassner, A.~R. Winters, and D.~A. Kopriva}, {\em Split form nodal
  discontinuous {G}alerkin schemes with summation-by-parts property for the
  compressible {E}uler equations}, J. Comput. Phys., 327 (2016), pp.~39--66.

\bibitem{giordano2006richtmyer}
{\sc J.~Giordano and Y.~Burtschell}, {\em {R}ichtmyer-{M}eshkov instability
  induced by shock-bubble interaction: Numerical and analytical studies with
  experimental validation}, Phys. Fluids, 18 (2006), p.~036102.

\bibitem{haas1987interaction}
{\sc J.~F. Haas and B.~Sturtevant}, {\em Interaction of weak shock waves with
  cylindrical and spherical gas inhomogeneities}, J. Fluid Mech., 181 (1987),
  pp.~41--76.

\bibitem{harten1983upstream}
{\sc A.~Harten, P.~D. Lax, and B.~V. Leer}, {\em On upstream differencing and
  {G}odunov-type schemes for hyperbolic conservation laws}, SIAM review, 25
  (1983), pp.~35--61.

\bibitem{helluy_seguin_phase_trans_06}
{\sc {Helluy, Philippe} and {Seguin, Nicolas}}, {\em Relaxation models of phase
  transition flows}, ESAIM: M2AN, 40 (2006), pp.~331--352.

\bibitem{hiltebrand2014entropy}
{\sc A.~Hiltebrand and S.~Mishra}, {\em Entropy stable shock capturing
  space--time discontinuous {G}alerkin schemes for systems of conservation
  laws}, Numer. Math., 126 (2014), pp.~103--151.

\bibitem{hiltebrand2018entropy}
{\sc A.~Hiltebrand, S.~Mishra, and C.~Par{\'e}s}, {\em Entropy-stable
  space--time {DG} schemes for non-conservative hyperbolic systems}, ESAIM:
  M2AN, 52 (2018), pp.~995--1022.

\bibitem{hu2006conservative}
{\sc X.~Y. Hu, B.~Khoo, N.~A. Adams, and F.~L. Huang}, {\em A conservative
  interface method for compressible flows}, J. Comput. Phys., 219 (2006),
  pp.~553--578.

\bibitem{ismail2009affordable}
{\sc F.~Ismail and P.~L. Roe}, {\em Affordable, entropy-consistent {E}uler flux
  functions ii: Entropy production at shocks}, J. Comput. Phys., 228 (2009),
  pp.~5410--5436.

\bibitem{jiang1994cell}
{\sc G.~S. Jiang and C.-W. Shu}, {\em On a cell entropy inequality for
  discontinuous {G}alerkin methods}, Math. {C}omput., 62 (1994), pp.~531--538.

\bibitem{johnsen2006implementation}
{\sc E.~Johnsen and T.~Colonius}, {\em Implementation of {WENO} schemes in
  compressible multicomponent flow problems}, J. Comput. Phys., 219 (2006),
  pp.~715--732.

\bibitem{kawai2011high}
{\sc S.~Kawai and H.~Terashima}, {\em A high-resolution scheme for compressible
  multicomponent flows with shock waves}, Int. J. Numer. Methods. Fluids, 66
  (2011), pp.~1207--1225.

\bibitem{kennedy2008reduced}
{\sc C.~A. Kennedy and A.~Gruber}, {\em Reduced aliasing formulations of the
  convective terms within the navier--stokes equations for a compressible
  fluid}, J. Comput. Phys., 227 (2008), pp.~1676--1700.

\bibitem{kopriva_metric_id_06}
{\sc D.~A. Kopriva}, {\em Metric identities and the discontinuous spectral
  element method on curvilinear meshes.}, J. Sci. Comput., 26 (2006),
  pp.~302--327.

\bibitem{kopriva2010quadrature}
{\sc D.~A. Kopriva and G.~Gassner}, {\em On the quadrature and weak form
  choices in collocation type discontinuous {G}alerkin spectral element
  methods}, J. Sci. Comput., 44 (2010), pp.~136--155.

\bibitem{lax1973hyperbolic}
{\sc P.~D. Lax}, {\em Hyperbolic systems of conservation laws and the
  mathematical theory of shock waves}, SIAM, 1973.

\bibitem{liu2003ghost}
{\sc T.~G. Liu, B.~C. Khoo, and K.~S. Yeo}, {\em Ghost fluid method for strong
  shock impacting on material interface}, J. Comput. Phys., 190 (2003),
  pp.~651--681.

\bibitem{renac_etal_ES_noneq_21}
{\sc C.~Marmignon, F.~Naddei, and F.~Renac}, {\em Energy relaxation
  approximation for compressible multicomponent flows in thermal
  nonequilibrium}, Numer. Math., 151 (2022), pp.~151--184.

\bibitem{MenikoffPlohr_RP_89}
{\sc R.~Menikoff and B.~J. Plohr}, {\em The riemann problem for fluid flow of
  real materials}, Rev. Mod. Phys., 61 (1989), pp.~75--130.

\bibitem{pares2006numerical}
{\sc C.~Par{\'e}s}, {\em Numerical methods for nonconservative hyperbolic
  systems: a theoretical framework.}, SIAM J. Numer. Anal., 44 (2006),
  pp.~300--321.

\bibitem{quirk1996dynamics}
{\sc J.~J. Quirk and S.~Karni}, {\em On the dynamics of a shock--bubble
  interaction}, J. Fluid Mech., 318 (1996), pp.~129--163.

\bibitem{rai2021modelling}
{\sc P.~Rai}, {\em Modelling and numerical simulation of compressible
  multicomponent flows}, PhD thesis, Institut Polytechnique de Paris, 2021.

\bibitem{ranocha_18}
{\sc H.~Ranocha}, {\em Comparison of some entropy conservative numerical fluxes
  for the {E}uler equations}, J. Sci. Comput., 76 (2018), pp.~216--242.

\bibitem{renac17b}
{\sc F.~Renac}, {\em A robust high-order discontinuous {G}alerkin method with
  large time steps for the compressible {E}uler equations}, Commun. Math. Sci.,
  15 (2017), pp.~813--837.

\bibitem{renac17a}
\leavevmode\vrule height 2pt depth -1.6pt width 23pt, {\em A robust high-order
  {L}agrange-projection like scheme with large time steps for the isentropic
  {E}uler equations}, Numer. Math., 135 (2017), pp.~493–--519.

\bibitem{renac2019entropy}
\leavevmode\vrule height 2pt depth -1.6pt width 23pt, {\em Entropy stable
  {DGSEM} for nonlinear hyperbolic systems in nonconservative form with
  application to two-phase flows}, J. Comput. Phys., 382 (2019), pp.~1--26.

\bibitem{renac2021multicomp}
\leavevmode\vrule height 2pt depth -1.6pt width 23pt, {\em Entropy stable,
  robust and high-order {DGSEM} for the compressible multicomponent {E}uler
  equations}, J. Comput. Phys.,  (2021), p.~110584.

\bibitem{renac2015aghora}
{\sc F.~Renac, M.~de~la Llave~Plata, E.~Martin, J.~B. Chapelier, and
  V.~Couaillier}, {\em Aghora: A High-Order {DG} Solver for Turbulent Flow
  Simulations}, Springer International Publishing, Cham, 2015, pp.~315--335.

\bibitem{saurel_pantano_ARFM_18}
{\sc R.~Saurel and C.~Pantano}, {\em Diffuse-interface capturing methods for
  compressible two-phase flows}, Annu. Rev. Fluid Mech., 50 (2018),
  pp.~105--130.

\bibitem{shu1988efficient}
{\sc C.-W. Shu and S.~Osher}, {\em Efficient implementation of essentially
  non-oscillatory shock-capturing schemes}, J. Comput. Phys., 77 (1988),
  pp.~439--471.

\bibitem{shyue1998efficient}
{\sc K.-M. Shyue}, {\em An efficient shock-capturing algorithm for compressible
  multicomponent problems}, J. Comput. Phys., 142 (1998), pp.~208--242.

\bibitem{sjogreen2003grid}
{\sc B.~Sj{\"o}green and H.~C. Yee}, {\em Grid convergence of high order
  methods for multiscale complex unsteady viscous compressible flows}, J.
  Comput. Phys., 185 (2003), pp.~1--26.

\bibitem{tadmor1987numerical}
{\sc E.~Tadmor}, {\em The numerical viscosity of entropy stable schemes for
  systems of conservation laws. {I}}, Math. {C}omput., 49 (1987), pp.~91--103.

\bibitem{tokareva2010hllc}
{\sc S.~A. Tokareva and E.~F. Toro}, {\em {HLLC}-type {R}iemann solver for the
  {B}aer--{N}unziato equations of compressible two-phase flow}, J. Comput.
  Phys., 229 (2010), pp.~3573--3604.

\bibitem{TORO_bound_wave_speed_2020}
{\sc E.~Toro, L.~M\"{u}ller, and A.~Siviglia}, {\em Bounds for wave speeds in
  the {R}iemann problem: Direct theoretical estimates}, Comput. Fluids, 209
  (2020), p.~104640.

\bibitem{toro1989riemann}
{\sc E.~F. Toro}, {\em {R}iemann-problem-based techniques for computing
  reactive two-phased flows}, in Numer. Combustion, Springer, 1989,
  pp.~472--481.

\bibitem{toro2013riemann}
\leavevmode\vrule height 2pt depth -1.6pt width 23pt, {\em {R}iemann solvers
  and numerical methods for fluid dynamics: a practical introduction}, Springer
  Science \& Business Media, 2013.

\bibitem{toro1994restoration}
{\sc E.~F. Toro, M.~Spruce, and W.~Speares}, {\em Restoration of the contact
  surface in the hll-riemann solver}, Shock waves, 4 (1994), pp.~25--34.

\bibitem{volpert1967spaces}
{\sc A.~Volpert}, {\em The spaces {BV} and quasilinear equations}, Math. USSR
  Sbornik, 115 (1967), pp.~255--302.

\bibitem{wang2012robust}
{\sc C.~Wang, X.~Zhang, C.-W. Shu, and J.~Ning}, {\em Robust high order
  discontinuous {G}alerkin schemes for two-dimensional gaseous detonations}, J.
  Comput. Phys., 231 (2012), pp.~653--665.

\bibitem{waruszewski2021entropy}
{\sc M.~Waruszewski, J.~E. Kozdon, L.~C. Wilcox, T.~H. Gibson, and F.~X.
  Giraldo}, {\em Entropy stable discontinuous galerkin methods for balance laws
  in non-conservative form: Applications to euler with gravity}, 2021.

\bibitem{wintermeyer2017entropy}
{\sc N.~Wintermeyer, A.~R. Winters, G.~J. Gassner, and D.~A. Kopriva}, {\em An
  entropy stable nodal discontinuous {G}alerkin method for the two dimensional
  shallow water equations on unstructured curvilinear meshes with discontinuous
  bathymetry}, J. Comput. Phys., 340 (2017), pp.~200--242.

\bibitem{zhang2010positivity}
{\sc X.~Zhang and C.~Shu}, {\em On positivity-preserving high order
  discontinuous {G}alerkin schemes for compressible {E}uler equations on
  rectangular meshes}, J. Comput. Phys., 229 (2010), pp.~8918--8934.

\bibitem{zhang2010maximum}
{\sc X.~Zhang and C.-W. Shu}, {\em On maximum-principle-satisfying high order
  schemes for scalar conservation laws}, J. Comput. Phys., 229 (2010),
  pp.~3091--3120.

\end{thebibliography}

\end{document}